

\input amstex
\documentstyle{amsppt}
\loadbold

\magnification=\magstep1
\hsize=6.5truein
\vsize=9truein

\document

\baselineskip=14pt

\font\smallit=cmti10 at 9pt
\font\smallsl=cmsl10 at 9pt
\font\sc=cmcsc10

\def \loongrightarrow {\relbar\joinrel\relbar\joinrel\rightarrow}
\def \llongrightarrow
{\relbar\joinrel\relbar\joinrel\relbar\joinrel\rightarrow}
\def \longtwoheadrightarrow
{\relbar\joinrel\relbar\joinrel\twoheadrightarrow}

\def \gerg {\frak g}

\def \gerk {\frak k}
\def \germ {\frak m}

\def \und1 {\underline{1}}

\def \h {\hbar}
\def \a {{\bold a}}
\def \boldalpha {{\boldsymbol\alpha}}
\def \b {{\bold b}}
\def \x {{\bold x}}
\def \id {\text{\rm id}}
\def \dif {{\text{dif}}}
\def \G {\Cal{G}}
\def \L {\Cal{L}}

\def \calA {\Cal{A}}

\def \calH {\Cal{H}}
\def \calK {\Cal{K}}
\def \calX {\Cal{X}}

\def \N {\Bbb{N}}

\def \C {\Bbb{C}}


\topmatter

\title
  Poisson geometrical symmetries associated to  \\   
  non-commutative formal diffeomorphisms    
\endtitle

\author
       Fabio Gavarini
\endauthor

\rightheadtext{ Poisson symmetries associated to
non-commutative diffeomorphisms }

\affil
  Universit\`a degli Studi di Roma ``Tor Vergata'' ---
Dipartimento di Matematica  \\
  Via della Ricerca Scientifica 1, I-00133 Roma --- ITALY  \\
\endaffil

\address\hskip-\parindent
  Fabio Gavarini  \newline
     \indent   Universit\`a degli Studi di Roma ``Tor Vergata''
---   Dipartimento di Matematica   \newline
     \indent   Via della Ricerca Scientifica 1, I-00133 Roma, ITALY
---   e-mail: gavarini\@{}mat.uniroma2.it
\endaddress

\abstract
   Let  $ \G^\dif $  be the group of all formal power series starting
with  $ x $  with coefficients in a field  $ \Bbbk $  of zero
characteristic (with the composition product), and let  $ \,
F \big[ \G^\dif \big] \, $  be its function algebra.  In [BF] a
non-commutative, non-cocommutative graded Hopf algebra  $ \calH^\dif $
was introduced via a direct process of ``disabelianisation'' of  $ \,
F \big[ \G^\dif \big] \, $,  \, taking the like presentation of the
latter as an algebra but dropping the commutativity constraint.  In
this paper we apply a general method to provide four one-parameters
deformations of  $ \calH^\dif $,  which are quantum groups whose
semiclassical limits are Poisson geometrical symmetries such as
Poisson groups or Lie bialgebras, namely two quantum function
algebras and two quantum universal enveloping algebras.  In
particular the two Poisson groups are extensions of  $ \G^\dif $,
isomorphic as proalgebraic Poisson varieties but not as proalgebraic
groups.   
\endabstract

\endtopmatter

\footnote""{Keywords: \ {\sl Hopf algebras, Quantum Groups}.}   

\footnote""{ 2000 {\it Mathematics Subject Classification:} \
Primary 16W30, 17B37, 20G42; Secondary 81R50. }


\hfill  \hbox{\vbox{ \baselineskip=10pt
     \hbox{\smallit \hskip53pt ``A series of outlaws joined and
formed the Nottingham group, }
 \hbox{\smallit \hskip71pt   whose renowned chieftain was the
famous Robin Hopf'' }
%
%
                   \vskip4pt
       \hbox{\smallsl  \hskip121pt  N.~Barbecue, ``Robin Hopf'' } }
\hskip1truecm }

%
%

\vskip1,7truecm

\centerline {\bf Introduction }

\vskip10pt

   The most general notion of ``symmetry'' in mathematics is encoded
in the notion of Hopf algebra.  Then, among all Hopf algebras (over a
field  $ \Bbbk $),  there are two special families which are of relevant
interest for their geometrical meaning: assuming for simplicity that 
$ \Bbbk $  have zero characteristic, these are the function algebras 
$ F[G] $  of algebraic groups  $ G $  and the universal enveloping
algebras  $ U(\gerg) $  of Lie algebras  $ \gerg \, $.  Function algebras
are exactly those Hopf algebras which are commutative, and enveloping
algebras those which are connected (in the general sense of Hopf
algebra theory) and cocommutative.   
                                         \par   
   Given a Hopf algebra  $ H $,  encoding some generalized symmetry, one can ask whether there are any other Hopf algebras ``close'' to  $ H $,  which are of either one of the above mentioned  {\sl geometrical types}, 
hence encoding geometrical symmetries associated to  $ H $.  The answer is affirmative: namely (see [Ga4]),
it is possible to give functorial recipes to get out of any Hopf algebra  $ H $  two pairs of Hopf algebras
of geometrical type, say  $ \big( F[G_+], U(\gerg_-) \big) $  and  $ \big( F[K_+], U(\gerk_-) \big) $. 
Moreover, the algebraic groups thus obtained are connected  {\sl Poisson\/}  groups,  and the Lie algebras
are  {\sl Lie bialgebras\/};  therefore in both cases  {\sl Poisson geometry\/}  is involved.  In addition,
the two pairs above are related to each other by Poisson duality (see below), thus only either one
of them is truly relevant.  Finally, these four ``geometrical'' Hopf algebras are ``close'' to  $ H $  in
that they are 1-para\-meter deformations (with pairwise isomorphic fibers) of a quotient or a subalgebra
of  $ H $.      
                                            \par   
   The method above to associate Poisson geometrical Hopf algebras to
general Hopf algebras, called ``Crystal Duality Principle'' (CDP in short),
is explained in detail in [Ga4].  It is a special instance of a more general
result, the ``Global Quantum Duality Principle'' (GQDP in short), explained
in [Ga2--3], which in turn is a generalization of the ``Quantum Duality
Principle'' due to Drinfeld (cf.~[Dr], \S 7, and see [Ga1] for a proof).
                                            \par   
   Drinfeld's QDP deals with quantum universal enveloping algebras (QUEAs
in short) and quantum formal series Hopf algebras (QFSHAs in short) over
the ring of formal power series  $ \Bbbk[[\h]] $.  A QUEA is any
topologically free, topological Hopf  $ \Bbbk[[\h]] $--algebra  whose
quotient modulo  $ \h $  is the universal enveloping algebra  $ U(\gerg) $ 
of some Lie algebra  $ \gerg \, $;  \, in this case we denote the QUEA
by  $ U_\h(\gerg) $.  Instead, a  QFSHA is any topological Hopf 
$ \Bbbk[[\h]] $--algebra  of type  $ {\Bbbk[[\h]]}^S $  (as a 
$ \Bbbk[[\h]] $--module,  $ S $  being a set) whose quotient modulo 
$ \h $  is the function algebra  $ F[[G]] $  of some formal algebraic
group  $ G \, $;  \,
   \hbox{then we denote the QFSHA by  $ F_\h[[G]] \, $.}  
                                            \par   
   The QDP claims that the category of all QUEAs and the category of all QFSHAs are equivalent, and
provides an equivalence in either direction.  From QFSHAs to QUEAs it goes as follows: given a QFSHA, say 
$ F_\h[[G]] $,  let  $ J $  be its augmentation ideal (the kernel of its counit map) and set  $ \,
{F_h[[G]]}^\vee := \sum_{n \geq 0} \h^{-n} J^n \, $.  Then  $ \, F_\h[[G]] \mapsto {F_h[[G]]}^\vee
\, $  defines (on objects) a functor from QFSHAs to QUEAs.  To go the other way round, i.e.~from QUEAs to
QFSHAs, one uses a perfectly dual recipe.  Namely, given a QUEA, say  $ U_\h(\gerg) $,  let again  $ J $ 
be its augmentation ideal; for each  $ \, n \in \N \, $,  let  $ \delta_n $  be the composition of the 
$ n $--fold  iterated coproduct followed by the projection onto  $ J^{\otimes n} $  (this makes sense
since  $ \, U_\h(\gerg) = \Bbbk[[\h]] \! \cdot \! 1_{\scriptscriptstyle U_\h(\gerg)} \oplus J \, $):  then 
set  $ \, {U_\h(\gerg)}' := \bigcap_{n \geq 0} {\delta_n}^{\!-1}\big(\h^n {U_\h(\gerg)}^{\otimes n}\big) \, $,  or more
explicitly  $ \, {U_\h(\gerg)}' := \big\{\, \eta \in U_\h(\gerg) \,\big|\, \delta_n(\eta) \in \h^n {U_\h(\gerg)}^{\otimes
n} , \; \forall \, n \in \N \,\big\} \, $.  Then  $ \, U_\h(\gerg) \mapsto {U_\h(\gerg)}' \, $  defines (on
objects) a functor from QUEAs to QFSHAs.  The functors  $ {(\ )}^\vee $ 
and  $ {(\ )}' $  are inverse to each other, hence they
provide the claimed equivalence.
                                            \par   
   Note that the objects (QUEAs and QFSHAs) involved in the QDP are 
{\sl quantum groups\/};  their semiclassical limits then are endowed
with Poisson structures: namely, every  $ U(\gerg) $  is in fact a
co-Poisson Hopf algebra and every  $ F[[G]] $  is a (topological)
Poisson Hopf algebra.  The geometrical structures they describe are
then  {\sl Lie bialgebras\/}  and  {\sl Poisson groups}.  The QDP then
brings further information: namely, the semiclassical limit of the
image of a given quantum group is  {\sl Poisson dual\/}  to the
   \hbox{Poisson geometrical object we start from.  In short}   
  $$  {F_h[[G]]}^\vee  \Big/ \h \, {F_h[[G]]}^\vee = \; U(\gerg^\times) \; ,  \qquad 
\text{i.e.~(roughly)}  \qquad  {F_h[[G]]}^\vee = \; U_\h(\gerg^\times)   \eqno (\hbox{I}.1)  $$
where  $ \gerg^\times $  is the  {\sl cotangent Lie bialgebra\/}  of the Poisson group  $ G \, $, 
\, and 
  $$  {U_\h(\gerg)}' \Big/ \h \, {U_\h(\gerg)}' \, = \, F\big[\big[G^\star\big]\big] \; ,  \qquad 
\text{i.e.~(roughly)}  \qquad  {U_\h(\gerg)}' \, = \, F_\h\big[\big[G^\star\big]\big]  
\eqno (\hbox{I}.2)  $$   
where  $ G^\star $  is a connected Poisson group with cotangent Lie bialgebra  $ \gerg \, $.  So the QDP invol\-ves
both Hopf duality 
   \hbox{(switching enveloping and function algebras) and Poisson duality.}   
                                            \par   
   The generalization from QDP to GQDP stems from a simple observation: the construction of Drinfeld's
functors needs not to start from quantum groups!  Indeed, in order to define either  $ H^\vee $  or 
$ H' $  one only needs that  $ H $  be a torsion-free Hopf algebra over some 1-dimensional doamin  $ R $ 
and  $ \, \h \in R \, $  be  {\sl any non-zero prime\/}  (actually, even less is truly necessary, see
[Ga2--3]).  On the other hand, the outcome still is, in both cases, a ``quantum group'', now meant in a
new sense.  Namely, a QUEA now will be any torsion-free Hopf algebra  $ H $  over  $ R $  such that  $ \, H
\Big/ \h \, H \cong U(\gerg) \, $,  \, for some Lie (bi)algebra  $ \gerg \, $.  Also, instead of
QFSHAs we consider ``quantum function algebras'', QFAs in short: here a QFA will be any torsion-free Hopf
algebra  $ H $  over  $ R $  such that  $ \, H \Big/ \h \, H \cong F[G] \, $  (plus one additional technical
condition) for some connected (Poisson) group  $ G \, $.  In this new framework Drinfeld's recipes give
that  $ H^\vee $  is a QUEA and  $ H' $  is a QFA, whatever is the torsion-free Hopf  $ R $--algebra  $ H $ 
one starts from.  Moreover, when restricted to quantum groups Drinfeld's functors  $ {(\ )}^\vee $  and 
$ {(\ )}' $  again provide equivalences of quantum group categories, respectively from QFAs to QUEAs and
viceversa; then Poisson duality 
   \hbox{is involved once more, like in (I.1--2).}      
                                            \par   
   Therefore, the generalization process from the QDP to the GQDP spreads
over several concerns.  Arithmetically, one can take as  $ (\h) $  any
non-generic point of the spectrum of  $ R \, $,  and define Drinfeld's
functors and specializations accordingly; in particular, the corresponding
quotient field  $ \, \Bbbk_\h := R \big/ \h \, R \, $  might have  {\sl
positive\/}  characteristic.  Geometrically, one considers algebraic
groups rather than formal groups,  i.e.~{\sl global\/}  vs.~{\sl local\/} 
objects.  Algebraically, one drops any topological worry  ($ \h $--adic 
completeness,  etc.), and deals with  {\sl general\/}  Hopf algebras
rather than with quantum groups.  This last point is the one of most
concern to us now, in that it means that we have (functorial) recipes
to get several quantum groups, hence   --- taking semiclassical limits
---   Poisson geometrical symmetries, springing out of the ``generalized
symmetry'' encoded by a torsion-free Hopf algebra  $ H $  over  $ R \, $: 
\, namely, for each non-trivial point of the spectrum of  $ R \, $,  the
quantum groups  $ H^\vee $  and  $ H' $  given by the corresponding
Drinfeld's functors.  Note, however, that  {\it a priori\/}  nothing
prevents any of these  $ H^\vee $  or  $ H' $  or their semiclassical
limits to be (essentially) trivial.    
                                            \par   
   The CDP comes out when looking at Hopf algebras over a field  $ \Bbbk
\, $,  and then applying the GQDP to their scalar extensions  $ \, H[\h]
:= \Bbbk[\h] \otimes_\Bbbk H \, $  with  $ \, R := \Bbbk[\h] \, $  (and 
$ \, \h := \h \, $  itself).  A first application of Drinfeld's functors
to  $ \, H_\h := H[\h] \, $  followed by specialization at  $ \, \h = 0
\, $  provides the pair  $ \big(F[G_+] \, , U(\gerg_-)\big) $  mentioned
above: in a nutshell,  $ \, \big(F[G_+] \, , U(\gerg_-)\big) = \Big(
{H_\h}^{\!\prime}\big|_{\h=0} \, , {H_\h}^{\!\vee}\big|_{\h=0} \Big) \, $, 
\, where hereafter  $ \, X\big|_{\h=0} := X \Big/ \h \, X \, $.  Then
applying once more Drinfeld's functors to  $ {H_\h}^{\!\vee} $  and to 
$ {H_\h}^{\!\prime} $  and specializing at  $ \, \h = 0 \, $  yields the
pair  $ \big(F[K_+] \, , U(\gerk_-)\big) $,  namely  $ \, \big(F[K_+]
\, , U(\gerk_-)\big) = \Big( {\big({H_\h}^{\!\vee}\big)}'\Big|_{\h=0} \, ,
{\big({H_\h}^{\!\prime}\big)}^{\!\vee}\Big|_{\h=0} \,\Big) \, $.  Finally,
the very last part of the GQDP explained before implies that  $ \, K_+ =
G^{\,\star}_- \, $  and  $ \, \gerk_- = \gerg_+^{\,\times} \, $.   
                                            \par   
   While in the second step above one really needs the full strength of
the GQDP, for the first step instead it turns out that the construction
of Drinfeld's functors on  $ H[\h] \, $,  can be fully ``tracked through''
and described at the ``classical level'', i.e.~in terms of  $ H $  alone. 
In addition, the exact relationship among  $ H $  and the pair  $ \big(
F[G_+] \, , U(\gerg_-) \big) $  can be made quite clear, and more
information is available about this pair.  We now sketch it in
some detail.   
                                            \par
   Let  $ \, J \, $  be the augmentation ideal of  $ H $,  let  $ \, \underline{J} :=
{\big\{ J^n \big\}}_{n \in \N} \, $  be the associated (decreasing)  $ J $--adic  filtration,  $ \,
\widehat{H} := G_{\underline{J}}(H) \, $  the associated graded vector space and  $ \, H^\vee := H \Big/
\bigcap_{n \in \N} J^n \, $.  One can prove that  $ \underline{J} $  is a Hopf algebra filtration, hence 
$ \widehat{H} $  is a graded Hopf algebra.  The latter happens to be connected and cocommutative, so  $ \,
\widehat{H} \cong U(\gerg_-) \, $  for some Lie algebra  $ \gerg_- \, $;  in addition, since  $ \widehat{H} $ 
is graded also  $ \gerg_- $  itself is graded as a Lie algebra.  The fact that  $ \widehat{H} $  be
cocommutative allows to define on it a Poisson cobracket 
%
%
which makes  $ \widehat{H} $  into a graded  {\sl co-Poisson\/}  Hopf algebra; eventually, this implies
that  $ \gerg_- $  is a  {\sl Lie bialgebra}.  The outcome is that our  $ U(\gerg_-) $  is just 
$ \widehat{H} $.
                                            \par
   On the other hand, one considers a second (increasing) filtration
defined in a dual manner to  $ \underline{J} \, $,  \, namely  $ \, \underline{D} := {\big\{ D_n := \text{\sl Ker}\,
(\delta_{n+1}) \big\}}_{n \in \N} \, $.  Let now  $ \, \widetilde{H}
:= G_{\underline{D}}(H) \, $  be the associated graded vector space
and  $ \, H' := \bigcup_{n \in \N} D_n \, $.  Again, one shows
that  $ \underline{D} $  is a Hopf algebra filtration, hence
$ \widetilde{H} $  is a graded Hopf algebra.  Moreover, the latter
is commutative, so  $ \, \widetilde{H} = F[G_+] \, $  for some
algebraic group  $ G_+ \, $.  One proves also that  $ \, \widetilde{H}
= F[G_+] \, $  has no non-trivial idempotents, thus  $ G_+ $  is
connected; in addition, since  $ \widehat{H} $  is graded,  $ G_+ $
as a variety is just an affine space.  The fact that  $ \widetilde{H} $ 
be commutative allows to define on it a Poisson bracket 
%
%
which makes  $ \widetilde{H} $  into a graded  {\sl Poisson\/}  Hopf
algebra: this means that  $ G_+ $  is an algebraic  {\sl Poisson group}. 
Thus eventually  $ F[G_+] $  is just  $ \widetilde{H} $.   
                                            \par
   The relationship among  $ H $  and the ``geometrical'' Hopf
algebras  $ \widehat{H} $  and  $ \widetilde{H} $  can be expressed
in terms of ``reduction steps'' and regular 1-parameter deformations,
namely
  $$  \widetilde{H}  \hskip1pt
\underset{ {\Cal{R}^\h_{\underline{D}}(H)}}  \to
{\overset{0 \,\leftarrow\, \h \,\rightarrow\, 1} \to
{\longleftarrow\joinrel\relbar\joinrel\relbar\joinrel%
\relbar\joinrel\llongrightarrow}}  \hskip1pt  H'  \hskip1pt
\,{\lhook\joinrel\relbar\joinrel\relbar\joinrel\rightarrow}\,
\hskip1pt  H  \relbar\joinrel\relbar\joinrel\twoheadrightarrow
H^\vee  \hskip0pt  \underset{ {\Cal{R}^\h_{\underline{J}}(H^\vee)}}
\to  {\overset{1 \,\leftarrow\, \h \,\rightarrow\, 0} \to
{\longleftarrow\joinrel\relbar\joinrel\relbar\joinrel%
\relbar\joinrel\relbar\joinrel\relbar\joinrel\longrightarrow}}
\hskip1pt   \widehat{H}   \eqno (\hbox{I}.3)  $$
where one-way arrows are Hopf algebra morphisms and
two-ways arrows are regular 1-para\-meter deformations
of Hopf algebras, realized through the  {\sl Rees\/}
Hopf algebras  $ {\Cal{R}^\h_{\underline{D}}(H)} $  and
$ {\Cal{R}^\h_{\underline{J}} (H^\vee)} $  associated to
the filtration  $ \underline{D} $  of  $ H $  and to
the filtration  $ \underline{J} $  of  $ H^\vee $.  Hereafter ``regular'' for a deformation means that all
its fibers are pairwise isomorphic as vector spaces.  In classical terms, (I.3) comes directly from the
construction above; on the other hand, in terms of the 
   \hbox{GQDP it comes from the fact that  $ \;
{\Cal{R}^\h_{\underline{D}}(H)} = {H_\h}^{\!\prime} \;
$  and  $ \; {\Cal{R}^\h_{\underline{J}} (H^\vee)} = {H_\h}^{\!\vee} \, $.}
                                            \par   
   As we mentioned above, next step is the ``application'' of (suitable) Drinfeld's functors to the Rees
algebras  $ \; {\Cal{R}^\h_{\underline{D}}(H)} = {H_\h}^{\!\prime} \; $  and  $ \; {\Cal{R}^\h_{\underline{J}}
(H^\vee)} = {H_\h}^{\!\vee} \, $ occurring in (I.3).  The outcome is a second frame of regular 1-parameter
deformations for  $ H' $  and  $ H^\vee $,  namely   
  $$  U\big(\gerg_+^{\,\times}\big) = U(\gerk_-) 
\hskip0pt  \underset{ (H'_\h)^\vee}  \to
{\overset{0 \,\leftarrow\, \h \,\rightarrow\, 1}
\to {\longleftarrow\joinrel\relbar\joinrel\relbar%
\joinrel\relbar\joinrel\llongrightarrow}}  \hskip0pt  H'
\hskip1pt  \,{\lhook\joinrel\relbar\joinrel\rightarrow}\,
\hskip1pt  H  \relbar\joinrel\relbar\joinrel\twoheadrightarrow
H^\vee  \hskip-1pt  \underset{(H^\vee_\h)'}
\to  {\overset{1 \,\leftarrow\, \h \,\rightarrow\, 0} \to
{\longleftarrow\joinrel\relbar\joinrel\relbar\joinrel\relbar%
\joinrel\relbar\joinrel\relbar\joinrel\longrightarrow}}
\hskip0pt   F[K_+] = F\big[G_-^{\,\star}\big]   \eqno (\hbox{I}.4)  $$
which is the analogue of (I.3).  In particular, when  $ \, H^\vee \!
= H = H' \, $  from (I.3) and (I.4) together we find  $ H $  as the
mid-point of four deformation families, whose ``external points'' are
Hopf algebras of ``Poisson geometrical'' type, namely
%
%
  $$  \hskip6pt   U(\gerg_-)  \hskip2pt
\underset{H^\vee_\h}  \to
{\overset{0 \,\leftarrow\, \h \,\rightarrow\, 1}
\to{\longleftarrow\joinrel\relbar\joinrel%
\relbar\joinrel\relbar\joinrel\relbar\joinrel\llongrightarrow}}
\hskip2pt {}  H  \hskip1pt
\underset{\;(H^\vee_\h)'}  \to
{\overset{1 \,\leftarrow\, \h \,\rightarrow\, 0}
\to{\longleftarrow\joinrel\relbar\joinrel%
\relbar\joinrel\relbar\joinrel\relbar\joinrel\llongrightarrow}}
\hskip2pt  F\big[G^{\,\star}_-\big]  $$
 \vskip-22pt
  $$  \hskip4pt   \Big|\Big|   \eqno (\maltese)  $$
 \vskip-15pt
  $$  \hskip5pt   F[G_+]  \hskip3pt
\underset{\;H'_\h}  \to
{\overset{0 \,\leftarrow\, \h \,\rightarrow\, 1}
\to{\longleftarrow\joinrel\relbar\joinrel%
\relbar\joinrel\relbar\joinrel\relbar\joinrel\llongrightarrow}}
\hskip2pt {}  H  \hskip1pt
\underset{\;(H'_\h)^\vee}  \to
{\overset{1 \,\leftarrow\, \h \,\rightarrow\, 0}
\to{\longleftarrow\joinrel\relbar\joinrel%
\relbar\joinrel\relbar\joinrel\relbar\joinrel\llongrightarrow}}
\hskip3pt  U\big(\gerg_+^{\,\times}\big)  $$
\noindent   which gives  {\sl four\/}  different regular 1-parameter
deformations from  $ H $  to Hopf algebras encoding Poisson geometrical objects.
Then each of these four Hopf algebras may be thought of as a semiclassical
   \hbox{geometrical counterpart of the ``generalized symmetry'' encoded by  $ H $.}
                                            \par   
   The purpose of the present paper is to show the effectiveness
of the CDP, applying it to a key example, the Hopf algebra of
non-commutative formal diffeomorphisms of the line.  Indeed, the
interest of the latter, besides its own reasons, grows bigger as
we can see it as a toy model for a broad family of Hopf algebras
of great  concern in mathematical physics, non-commutative
geometry and beyond.  Now I go and present the results of this paper.  
                                            \par   
   Let  $ \G^\dif $  be the set of all formal power series starting
with  $ x $  with coefficients in a field  $ \Bbbk $  of zero
characteristic.  Endowed with the composition product, this is
an infinite dimensional prounipotent proalgebraic group   ---
known as the ``(normalised) Nottingham group'' among group-theorists
and the ``(normalised) group of formal diffeomorphisms of the line''
among mathematical physicists ---   whose tangent Lie algebra is a
special subalgebra of the one-sided Witt algebra.  The function
algebra  $ F \big[ \G^\dif \big] $  is a graded, commutative Hopf
algebra with countably many generators, which admits a neat combinatorial
description.   
                                            \par   
  In [BF] a non-commutative version of  $ F \big[ \G^\dif \big] $ 
is introduced: this is a non-commutative non-cocommutative Hopf
algebra  $ \calH^\dif $  which is presented exactly like  $ F \big[
\G^\dif \big] $  but dropping commutativity, i.e.~taking the
presentation as one of a unital  {\sl associative\/}   ---
and  {\sl not\/}  commutative ---   algebra; in other words, 
$ \calH^\dif $  is the outcome of applying to  $ F \big[
\G^\dif \big] $  a raw ``disabelianization'' process.  In
particular,  $ \, H = \calH^\dif \, $  is  {\sl graded\/} 
and verifies  $ \, H^\vee = H = H' \, $,  \, hence the scheme 
$ (\maltese) $  makes sense and yields four Poisson symmetries
associated to  $ \, \calH^\dif \, $.   
                                             \par
   Note that in each line in  $ (\maltese) $  there is essentially
only  {\sl one\/}  Poisson geometry involved, since Poisson duality
relates mutually opposite sides; thus any classical symmetry on
the same line carries as much information as the other one (but for
global-to-local differences).  Nevertheless, in the case of  $ \, H = \calH^\dif \, $ 
we shall prove that the pieces of information from either line in  $ (\maltese) $  are complementary,
because  $ G_+ $  and  $ G^{\,\star}_- $  happen to be isomorphic as
proalgebraic Poisson varieties but  {\sl not\/}  as groups.  In
particular, we find that the Lie bialgebras  $ \gerg_- $  and 
$ \gerg_+^{\,\times} $  are both isomorphic as Lie algebras to
the free Lie algebra  $ \Cal{L}(\N_+) $  over a countable set,
but they have different, non-isomorphic Lie coalgebra structures.
Moreover,  $ \, G^{\,\star}_- \cong \G^\dif \times \Cal{N} \cong G_+ \, $ 
as Poisson varieties, where  $ \Cal{N} $  is a proaffine Poisson variety
whose coordinate functions are in bijection with a basis of the derived
subalgebra  $ \Cal{L}(\N_+) $;  \, indeed, the latter are obtained by
iterated Poisson brackets of coordinate functions on  $ \G^\dif $,  in
short because both  $ \, F \big[ G^{\,\star}_- \big] \, $  and  $ \, F
\big[ G_+ \big] \, $  are freely generated as Poisson algebras by a copy
of  $ \, F \big[ \G^\dif \big] \, $.  For  $ G^{\,\star}_- $  we have a
more precise result, namely  $ \, G^{\,\star}_- \cong \G^\dif \ltimes
\Cal{N} \, $  (a semidirect product)  {\sl as proalgebraic groups\/}: 
\, thus in a sense  $ G^{\,\star}_- $  is the  {\sl free Poisson
group over  $ \G^\dif $},  which geometrically speaking is obtained by
``pasting'' to  $ \G^\dif $  all 1-parameter subgroups freely obtained
via iterated Poisson brackets of those of  $ \G^\dif \, $;  \, in
particular, these Poisson brackets iteratively yield 1-parameter
subgroups which generate  $ \Cal{N} $.
                                               \par
   We perform the same analysis simultaneously for  $ \G^\dif $,  for
its subgroup of  {\sl odd\/}  formal diffeomorphisms and for all the
groups  $ \G_\nu $  of  {\sl truncated\/}  (at order  $ \nu \in \N_+ $)
formal diffeomorphisms,
         \hbox{whose projective limit is  $ \G^\dif $  itself;
{\it mutatis mutandis},  the results are the like}.
                                                \par
   The case of  $ \calH^\dif $  is just one of many samples of
the same type: indeed, several cases of Hopf algebras built out
of combinatorial data   --- graphs, trees, Feynman diagrams, etc.~---  
have been introduced in (co)homological theories (see e.g.~[LR] and
[Fo1--2], and references therein) and in renormalization studies (see
[CK1--3]); in most cases these algebras   --- or their (graded) duals
---   are commutative polynomial, like  $ F \big[ \G^\dif \big] $, 
and admit non-commutative analogues (thanks to [Fo1--2]), so our
discussion apply almost  {\it verbatim\/}  to them
too, with like results.  Thus the given analysis of
the ``toy model'' Hopf algebra  $ \calH^\dif $  can
be taken as a general pattern for all those cases.

\vskip9pt

\centerline{ \sc acknowledgements }

  The author thanks Alessandra Frabetti and Loic Foissy
for many helpful discussions.

\vskip1,5truecm

\centerline {\bf \S \; 1 \ Notation and terminology }

\vskip10pt

  {\bf 1.1 The classical data.} \, Let  $ \Bbbk $  be a fixed field
of zero characteristic.
                                                     \par
   Consider the set  $ \; \G^\dif := \big\{\, x + \sum_{n \geq 1} a_n
\, x^{n+1} \,\big|\; a_n \in \Bbbk \; \forall \; n \in \N_+ \,\big\}
\; $  of all formal series starting with $ x \, $:  endowed with the
composition product, this is a group, which can be seen as the group
of all ``formal diffeomorphisms''   \hbox{$ \, f \, \colon \, \Bbbk
\longrightarrow \Bbbk \, $}  such that  $ \, f(0) = 0 \, $  and  $ \,
f'(0) = 1 \, $  (i.e.~tangent to the identity), also known as the  {\sl
Nottingham group\/}  (see, e.g., [Ca] and references therein).  In fact,
$ \G^\dif $  is an infinite dimensional (pro)affine algebraic group,
whose function algebra  $ \, F \big[ \G^\dif \big] \, $  is generated
by the coordinate functions  $ \, a_n \, $  ($ \, n \in \N_+ $).  Giving
to each  $ a_n $  the weight\footnote{We say  {\sl weight\/} instead
        of  {\sl degree\/}  because we save the latter term for
        the degree of polynomials.}
 $ \, \partial(a_n) := n \, $,  we have that  $ F \big[ \G^\dif \big] $
is an  $ \N $--{\sl graded\/}  Hopf algebra, with polynomial structure
$ \, F \big[ \G^\dif \big] = \Bbbk [a_1,a_2,\dots,a_n,\dots] \, $  and
Hopf algebra structure given by
  $$  \eqalign{
   \Delta(a_n) \, = \, a_n \otimes 1 + 1 \otimes a_n + {\textstyle
\sum\nolimits_{m=1}^{n-1}} \, a_m \otimes Q^m_{n-m}(a_*) \, ,
\qquad  \epsilon(a_n) \, = \, 0   \hskip15pt  \cr
   S(a_n) \, = \, - a_n - {\textstyle \sum\nolimits_{m=1}^{n-1}}
\, a_m \, S\big(Q^m_{n-m}(a_*)\big) \, = \, - a_n - {\textstyle
\sum\nolimits_{m=1}^{n-1}} \, S(a_m) \, Q^m_{n-m}(a_*) \cr }  $$
where  $ \; Q_t^\ell(a_*) := \sum_{k=1}^t {{\ell + 1} \choose k}
P_t^{(k)}(a_*) \; $  and  $ \; P_t^{(k)}(a_*) := \sum_{\hskip-4pt
\Sb  j_1, \dots, j_k > 0  \\   j_1 + \cdots + j_k = t  \endSb}
\hskip-5pt  a_{j_1} \cdots a_{j_k} \, $  (the symmetric monic
polynomial of weight  $ m $  and degree  $ k $  in the indeterminates
$ a_j $'s)  for all  $ \, m $,  $ k $,  $ \ell \in \N_+ \, $,  \, and
the formula for  $ S(a_n) $  gives the antipode by recursion.  From now
on, to simplify notation we shall write  $ \, \G := \G^\dif \, $  and
$ \, \G_\infty := \G = \G^\dif \, $.  Note also that the tangent Lie
algebra of  $ \G^\dif $  is just the Lie subalgebra  $ \, {W_1}^{\!
\geq 1} = \text{\sl Span}\, \big(\{\, d_n \,|\, n \in \N_+ \,\}\big)
\, $  of the one-sided Witt algebra  $ \, W_1 := \text{\sl Der} \big(
\Bbbk[t] \big) = \text{\sl Span}\, \big( \big\{\, d_n := t^{n+1} {{\,d\,}
\over {\,dt\,}} \,\big|\, n \in \N \cup \{-1\} \,\big\}\big) \, $.
                                                     \par
   In addition, for all  $ \, \nu \in \N_+ \, $  the subset
%
%
$ \, \G^\nu := \big\{\, f \in \G \,\big|\; a_n(f) = 0 , \;
\forall \; n \leq \nu \,\big\} \, $  is a normal subgroup of
$ \G \, $;  \, the corresponding quotient group  $ \, \G_\nu :=
\G \big/ \G^\nu \, $  is unipotent, with dimension  $ \nu $  and
function algebra  $ \, F \big[ \G_\nu \big] \, $  (isomorphic to)
the Hopf subalgebra of  $ F \big[ \G \big] $  generated by  $ \,
a_1 $,  $ \dots $,  $ a_\nu \, $.  In fact, the  $ \G^\nu $'s  form
exactly the lower central series of  $ \G $  (cf.~[Je]).  Moreover,
$ \G $  is (isomorphic to) the inverse (or projective) limit of these
quotient groups  $ \G_\nu $  ($ \nu \in \N_+ $),  hence  $ \G $  is
pro-unipotent; conversely,  $ F[\G] $  is the direct (or inductive)
limit of the direct system of its graded Hopf subalgebras  $ F[\G_\nu]
\, $  ($ \nu \in \N_+ $).  Finally, the set  $ \; \G^{\text{odd}} :=
\big\{\, f \in \G^\dif \,\big|\; a_{2n-1}(f) = 0 \; \forall \;
n \in \N_+ \,\big\} \; $  is another normal subgroup of
$ \G^\dif $  (the group of  {\sl odd\/}
       formal diffeomorphisms\footnote{\hbox{The fixed-point set of
       the group homomorphism  $ \, \Phi : \G \rightarrow \G \, $,
       $ f \mapsto \Phi(f) \, \big( x \mapsto \big(\Phi(f)\big)(x)
       := -f(-x) \big) $}}
 after [CK3]), whose function algebra  $ F \big[ \G^{\text{odd}} \big] $
is (isomorphic to) the quotient Hopf algebra  $ \, F \big[ \G^\dif \big]
\Big/ \Big( \big\{ a_{2n-1} \big\}_{n \in \N_+} \Big) \, $.  The
latter has the following description: denoting again the cosets of the
$ a_{2n} $'s  with the like symbol, we have  $ \, F \big[ \G^{\text{odd}}
\big] = \Bbbk [a_2, a_4, \dots, a_{2n}, \dots] \, $  with Hopf algebra
structure
  $$  \eqalign{
   \Delta(a_{2n}) \, = \, a_{2n} \otimes 1 + 1 \otimes a_{2n}
+ {\textstyle \sum\nolimits_{m=1}^{n-1}} \, a_{2m} \otimes
\bar{Q}^m_{n-m}(a_{2*}) \, ,  \qquad  \epsilon(a_{2n}) \, = \, 0
\hskip35pt  \cr
   S(a_{2n}) \, = \, - a_{2n} - {\textstyle \sum\nolimits_{m=1}^{n-1}}
\, a_{2m} \, S\big(\bar{Q}^m_{n-m}(a_*)\big) = - a_{2n} - {\textstyle
\sum\nolimits_{m=1}^{n-1}} \, S(a_{2m}) \, \bar{Q}^m_{n-m}(a_{2*})
\cr }  $$
where  $ \; \bar{Q}_t^\ell(a_{2*}) := \sum_{k=1}^t {{2\ell + 1}
\choose k} \bar{P}_t^{(k)}(a_{2*}) \; $  and  $ \; \bar{P}_t^{(k)}
(a_{2*}) := \sum_{\hskip-4pt  \Sb  j_1, \dots, j_k > 0  \\
j_1 + \cdots + j_k = t  \endSb}  \hskip-5pt  a_{2{}j_1} \cdots
a_{2{}j_k} \, $  for all  $ \, m $,  $ k $,  $ \ell \in \N_+ \, $.  For
each  $ \, \nu \in \N_+ \, $  we can consider also the normal subgroup
$ \, \G^\nu \cap \G^{\text{odd}} \, $  and the corresponding quotient
$ \, \G_\nu^{\text{odd}} := \G^{\text{odd}} \big/ \big( \G^\nu \cap
\G^{\text{odd}} \big) \, $:  \, then  $ F \big[ \G_\nu^{\text{odd}}
\big] $  is (isomorphic to) the quotient Hopf algebra  $ \, F \big[
\G^{\text{odd}} \big] \! \Big/ \! \Big( \! \big\{ a_{2n-1} \big\}_{(2n-1)
\in \N_\nu} \Big) \, $,  \, in particular it is the Hopf sub-\break
 \noindent   algebra of  $ F \big[ \G^{\text{odd}} \big] $  generated
by  $ \, a_2, \dots, a_{2\,[\nu/2]} \, $.  All the  $ F \big[
\G_\nu^{\text{odd}} \big] $'s  are graded Hopf (sub)al-\break
 \noindent   gebras forming a direct system with direct limit
$ F \big[ \G^{\text{odd}} \big] $;  \, conversely, the
$ \G_\nu^{\text{odd}} $'s  form an inverse system with inverse
limit  $ \G^{\text{odd}} $.  In the sequel we write  $ \, \G^+ :=
\G^{\text{odd}} \, $  and  $ \, \G^+_\nu := \G_\nu^{\text{odd}} \, $.
                                                     \par
   For each  $ \, \nu \in \N_+ \, $,  \, set  $ \, \N_\nu := \{1, \dots,
\nu\} \, $;  set also  $ \, \N_\infty := \N_+ \, $.  For each  $ \, \nu
\in \N_+ \cup \{\infty\} \, $,  \, let  $ \, \L_\nu = \L(\N_\nu) \, $
be the free Lie algebra over  $ \Bbbk $  generated by  $ {\{x_n\}}_{n
\in \N_\nu} $  and let  $ \, U_\nu = U(\L_\nu) \, $  be its universal
enveloping algebra; let also  $ \, V_\nu = V(\N_\nu) \, $  be the
$ \Bbbk $--vector  space with basis  $ \, {\{ x_n \}}_{n \in \N_\nu}
\, $,  \, and let  $ \, T_\nu = T(V_\nu) \, $  be its associated tensor
algebra.  Then there are canonical identifications  $ \, U(\L_\nu)
= T(V_\nu) = \Bbbk \big\langle \{\, x_n \,|\, n \in \N_\nu \,\}
\big\rangle \, $,  \, the latter being the unital  $ \Bbbk $--algebra
of non-commutative polynomials in the set of indeterminates  $ \,
{\{x_n\}}_{n \in \N_\nu} \, $,  \, and $ \L_\nu $  is just the Lie
subalgebra of  $ U_\nu = T_\nu $  generated by  $ \, {\{x_n\}}_{n \in
\N_\nu} \, $.  Moreover,  $ \L_\nu $  has a basis  $ B_\nu $  made of
Lie monomials in the  $ x_n $'s  ($ n \in \N_\nu $),  like  $ [x_{n_1},
x_{n_2}] $,  $ [[x_{n_1}, x_{n_2}], x_{n_3}] $,  $ [[[x_{n_1}, x_{n_2}],
x_{n_3}], x_{n_4}] $,  etc.: details can be found e.g.~in [Re], Ch.~4--5.
In the sequel I shall use these identifications with no further mention.
We consider on  $ \, U(\L_\nu) \, $  the standard Hopf algebra structure
given by  $ \; \Delta(x) = x \otimes 1 + 1 \otimes x \, $,  $ \,
\epsilon(x) = 0 \, $,  $ \, S(x) = -x \; $  for all  $ \, x \in
\L_\nu \, $,  \, which is also determined by the same formulas for
$ \, x \in {\{x_n\}}_{n \in \N_\nu} \, $  alone.  By construction
$ \, \nu \leq \mu \, $  implies  $ \, \L_\nu \subseteq \L_\mu \, $,
\, whence the  $ \L_\nu $'s  form a direct system (of Lie algebras)
whose direct limit is exactly  $ \L_\infty \, $;  \, similarly,
$ U(\L_\infty) $  is the direct limit of all the  $ U(\L_\nu) $'s.
Finally, with  $ \Bbb{B}_\nu $  we shall mean the obvious PBW-like
basis of  $ U(\L_\nu) $  w.r.t.~some fixed total order  $ \preceq $
of  $ B_\nu $,  namely  $ \, \Bbb{B}_\nu := \big\{\, x_{\underline{b}}
\;\big|\, \underline{b} = b_1 \cdots b_k \, ; \, b_1, \dots, b_k
\in B_\nu \, ; \, b_1 \preceq \cdots \preceq b_k \,\big\} \, $.
   The same construction applies to make out ``odd'' objects, based
on  $ {\{x_n\}}_{n \in \N^+_\nu} $,  with  $ \, \N^+_\nu := \N_\nu
\cap 2\,\N \, $  ($ \nu \! \in \! \N \cup \{\infty\} $),  instead of 
$ {\{x_n\}}_{n \in \N_\nu} $,  $ \, \L^+_\nu = \L(\N^+_\nu) \, $, 
$ \, U^+_\nu = U(\L^+_\nu) \, $,  \, $ \, V^+_\nu = V(\N^+_\nu) \, $, 
$ \, T^+_\nu = T(V^+_\nu) \, $,  \, with the obvious canonical
identifications  $ \, U(\L^+_\nu) = T(V^+_\nu) = \Bbbk \big\langle
\{\, x_n \,|\, n \in \N^+_\nu \,\} \big\rangle \, $;  \, moreover, 
$ \L^+_\nu $  has a basis  $ B^+_\nu $  made of Lie monomials in the 
$ x_n $'s  ($ n \in \N^+_\nu $),  etc.  The  $ \L^+_\nu $'s  form
a direct system whose direct limit is  $ \L^+_\infty \, $,  \, and 
$ U(\L^+_\infty) $  is the direct limit of all the  $ U(\L^+_\nu) $'s.

\vskip4pt

   {\it  $ \underline{\text{Warning}} \, $:}  \, in the sequel, we
shall often deal with subsets  $ {\{\hbox{\bf y}_b\}}_{b \in B_\nu} $
(of some algebra) in bijection with  $ B_\nu \, $,  the fixed basis
of  $ \L_\nu \, $.  Then we shall write things like  $ \, \hbox{\bf
y}_\lambda \, $  with  $ \, \lambda \in \L_\nu \, $:  \, this means
we extend the bijection  $ \, {\{\hbox{\bf y}_b\}}_{b \in B_\nu} \cong
B_\nu \, $  to  $ \, \text{\sl Span}\,\big( \{ \hbox{\bf y}_b \}_{b \in
B_\nu} \big) \cong \L_\nu \, $  by linearity, so that  $ \; \hbox{\bf
y}_\lambda \cong \sum_{b \in B_\nu} c_b \, b \; $  iff  $ \; \lambda
= \sum_{b \in B_\nu} c_b \, b \; $  ($ \, c_b \in \Bbbk \, $).  The
same kind of convention will be applied with  $ B^+_\nu $  instead
of  $ B_\nu $  and  $ \L^+_\nu $  instead of  $ \L_\nu \, $.

\vskip7pt

   {\bf 1.2 The noncommutative Hopf algebra of formal diffeomorphisms.}
\, For all  $ \, \nu \in \N_+ \cup \{\infty\} \, $,  \, let  $ \,
\calH_\nu \, $  be the Hopf  $ \Bbbk $--algebra  given as follows:
as a  $ \Bbbk $--algebra it is simply  $ \, \calH_\nu := \Bbbk
\big\langle \{\, \a_n \,|\, n \in \N_\nu \,\} \big\rangle \, $
(the  $ \Bbbk $--algebra  of non-commutative polynomials in the
set of indeterminates  $ \, {\{\a_n\}}_{n \in \N_\nu} \, $),  and its
Hopf algebra structure is given by (for all  $ \, n \in \N_\nu \, $)
  $$  \hbox{ $ \eqalign{
   \Delta(\a_n) \, = \, \a_n \otimes 1 + 1 \otimes \a_n + {\textstyle
\sum_{m=1}^{n-1}} \, \a_m \otimes Q^m_{n-m}(\a_*) \, ,  \qquad
\hskip15pt  \epsilon(\a_n) = 0  \hskip15pt  \cr
   S(\a_n) \, = \, -\a_n - {\textstyle \sum_{m=1}^{n-1}} \,
\a_m \, S\big(Q^m_{n-m}(\a_*)\big) \, = \, -\a_n - {\textstyle
\sum_{m=1}^{n-1}} \, S(\a_m) \, Q^m_{n-m}(\a_*)  \cr } $ }
\eqno  (1.1)  $$
(notation like in \S 1.1) where the latter formula yields the
antipode by recursion.  Moreover,  $ \calH_\nu $  is in fact an  {\sl
$ \N $--graded  Hopf algebra},  once generators have been given degree
--- in the sequel called  {\sl weight\/}  ---   by the rule  $ \,
\partial(\a_n) := n \, $  (for all  $ \, n \in \N_\nu \, $).  By
construction the various  $ \calH_\nu $'s  (for all  $ \, \nu \in
\N_+ \, $)  form a direct system, whose direct limit is  $ \calH_\infty
\, $:  \, the latter was
    originally introduced\footnote{However, the formulas in [BF]
    give the  {\sl opposite\/}  coproduct, hence change the antipode
    accordingly; we made the present choice to make these formulas
    ``fit well'' with those for  $ F \big[ \G^{\text{dif}} \big] $
    (see below).}
in [BF], \S 5.1 (with  $ \, \Bbbk = \C \, $),  under the name
$ \, \calH^\dif \, $.
                                                     \par
   Similarly, for all  $ \, \nu \in \N_+ \cup \{\infty\} \, $  we
set  $ \, \calK_\nu := \Bbbk \big\langle \{\, \a_n \,|\, n \in
\N^+_\nu \,\} \big\rangle \, $  (where  $ \, \N^+_\nu := \N_\nu
\cap (2\,\N) \, $):  this bears a Hopf algebra structure given
by (for all  $ \, 2\,n \in \N^+_\nu \, $)
  $$  \hbox{ $ \eqalign{
   \Delta(\a_{2n}) \, = \, \a_{2n} \otimes 1 + 1 \otimes \a_{2n} +
{\textstyle \sum_{m=1}^{n-1}} \, \a_{2m} \otimes \bar{Q}^m_{n-m}(\a_{2*})
\, ,  \qquad  \hskip11pt  \epsilon(\a_{2n}) = 0
\hskip19pt  \cr
   S(\a_{2n}) \, = \, -\a_{2n} - {\textstyle \sum_{m=1}^{n-1}} \,
\a_{2m} \, S\big(\bar{Q}^m_{n-m}(\a_{2*})\big) \, = \, - \a_{2n}
- {\textstyle \sum_{m=1}^{n-1}} \, S(\a_{2m}) \,
Q^m_{n-m}(\a_{2*})  \cr } $ }  $$
(notation of \S 1.1).  Indeed, this is an  $ \N $--graded  Hopf algebra
where generators have degree   --- called  {\sl weight\/}  ---   given
by  $ \, \partial(\a_n) := n \, $  (for all  $ \, n \in \N^+_\nu \, $).
All the  $ \calK_\nu $'s  form a direct system with direct limit
$ \calK_\infty \, $.  Finally, for each  $ \, \nu \in \N^+_\nu \, $
there is a graded Hopf algebra epimorphism  $ \, \calH_\nu
\longtwoheadrightarrow \calK_\nu \, $  given by  $ \, \a_{2n}
\mapsto \a_{2n} \, $,  $ \, \a_{2m+1} \mapsto 0 \, $  for
all  $ \, 2n, 2m+1 \in \N_\nu \, $.
                                             \par
   Definitions and \S 1.1 imply that
  $$  {\big(\calH_\nu\big)}_{\text{\it ab}} \, := \,
\calH_\nu \Big/ \big( \big[ \calH_\nu, \calH_\nu \,\big] \big)
\; \cong \; F \big[ \G_\nu \big] \, ,   \qquad   \hbox{via}
\qquad   \a_n \mapsto a_n  \quad  \forall \; n \in \N_\nu  $$
as  $ \N $--graded  Hopf algebras: in other words, the abelianization
of  $ \calH_\nu $  is nothing but  $ F \big[ \G_\nu \big] $.  Thus in
a sense one can think at  $ \calH_\nu $  as a  {\sl non-commutative
version\/}  (indeed, the ``coarsest'' one) of  $ F \big[ \G_\nu \big] $,
hence as a ``quantization'' of  $ \G_\nu $  itself: however, this is
{\sl not\/}  a quantization in the usual sense, because  $ F \big[
\G_\nu \big] $  is attained through abelianization,  {\sl not\/}  via
specialization of some deformation parameter.  Similarly we have also
  $$  {\big(\calK_\nu\big)}_{\text{\it ab}} \, := \,
\calK_\nu \Big/ \big( \big[ \calK_\nu, \calK_\nu \,\big] \big)
\; \cong \; F \big[ \G^+_\nu \big] \, ,   \qquad   \hbox{via}
\qquad   \a_{2n} \mapsto a_{2n}  \quad  \forall \; 2n \in \N^+_\nu  $$
as  $ \N $--graded  Hopf algebras: in other words, the abelianization
of  $ \calK_\nu $  is just  $ F \big[ \G^+_\nu \big] $.
                                                     \par
   In the following I make the analysis explicit for  $ \calH_\nu \, $,
\, the case  $ \calK_\nu $  being the like (details are left to the
reader); I drop the subscript  $ \nu $,  which stands fixed, and write
$ \, \calH := \calH_\nu \, $.

\vskip7pt

   {\bf 1.3 Deformations.} \, Let  $ \h $  be an indeterminate.  In
this paper we shall consider several Hopf algebras over  $ \Bbbk[\h] $,
which can also be seen as 1-parameter depending families of Hopf
algebras over  $ \Bbbk $,  the parameter being  $ \Bbbk \, $;  \,
each  $ \Bbbk $--algebra  in such a family can then be thought
of as a 1-parameter deformation of any other object in the
same family. As a matter of notation, if  $ H $  is such
a Hopf  $ \Bbbk[\h] $--algebra  I call  {\sl fibre\/}  of
$ H $  (though of as a deformation) any Hopf  $ \Bbbk $--algebra
of type  $ \, H \Big/ p(\h) \, H \; $  for some irreducible  $ \,
p(\h) \in \Bbbk[\h] \, $;  \, in particular  $ \, H{\Big|}_{\h=c}
\!\! := H \! \Big/ (\h-c) \, H \, $,  \, for any  $ \, c \in \Bbbk \, $,
\, is called  {\sl specialization\/}  of  $ H $  at  $ \, \h = c \, $.
                                        \par
   We start from  $ \; \calH_\h := \calH[\h] \equiv \Bbbk[\h]
\otimes_{\Bbbk} \! \calH \, $:  \, this is indeed a Hopf
$ \Bbbk[\h] $--algebra,  namely  $ \, \calH_\h = \Bbbk[\h]
\big\langle \{\, \a_n \,|\, n \in \N_\nu \,\} \big\rangle \, $
with Hopf structure given by (1.1) again.  Set also  $ \, \calH(\h)
:= \Bbbk(\h) \otimes_{\Bbbk[\h]} \! \calH_\h = \Bbbk(\h) \otimes_{\Bbbk}
\! \calH = \Bbbk(\h) \big\langle \{\, \a_n \,|\, n \in \N_\nu \,\}
\big\rangle $,  \, a Hopf  $ \Bbbk(\h) $--algebra ruled by
(1.1) too.

\vskip1,5truecm

\centerline {\bf \S \; 2 \ The Rees deformation
$ \, {\calH_\h}^{\!\vee} \, $. }

\vskip10pt

  {\bf 2.1 The goal.} \, The crystal duality principle (cf.~[Ga2], \S
5, or [Ga4]) yields a recipe to produce a 1-parameter deformation  $ \,
{\calH_\h}^{\!\vee} \, $  of  $ \calH $  which is a  {\sl quantized
universal enveloping algebra\/} (QUEA in the sequel): namely,  $ \,
{\calH_\h}^{\!\vee} \, $  is a Hopf  $ \Bbbk[\h] $--algebra  such that
$ \, {\calH_\h}^{\!\vee}\Big|_{\h=1} \!\! = \calH \, $  and  $ \,
{\calH_\h}^{\!\vee}\Big|_{\h=0} \!\! = U(\gerg_-) \, $,  \, the
universal enveloping algebra of a graded Lie bialgebra  $ \gerg_-
\, $.  Thus $ {\calH_\h}^{\!\vee} $  is a  {\sl quantization\/}
of  $ U(\gerg_-) $,  and the quantum symmetry  $ \calH $  is a
deformation of the classical Poisson symmetry  $ U(\gerg_-) $.
By definition  $ {\calH_\h}^{\!\vee} $  is the  {\sl Rees algebra\/}
associated to a distinguished  {\sl decreasing\/}  Hopf algebra
filtration of  $ \calH $,  so that  $ U(\gerg_-) $  is just the graded
Hopf algebra associated to this filtration.  The purpose of this section
is to describe explicitly  $ {\calH_\h}^{\!\vee} $  and its semiclassical
limit  $ U(\gerg_-) $,  hence also  $ \gerg_- $  itself.  This will also
provide a direct, independent proof of all the above mentioned results
about  $ {\calH_\h}^{\!\vee} $  and  $ U(\gerg_-) $  themselves.

\vskip7pt

  {\bf 2.2 The Rees algebra  $ \, {\calH_\h}^{\!\vee} \, $.} \, Let
$ \, J := \text{\sl Ker}\,\big(\epsilon_{\scriptscriptstyle \calH} \,
\colon \calH \longrightarrow \Bbbk \, \big) \, $  be the augmentation
ideal of  $ \calH $,  and let  $ \, \underline{J} := {\big\{ J^n
\big\}}_{n \in \N} \, $  be the  $ J $--adic  filtration in  $ \calH
\, $.  It is easy to show (see [Ga4]) that  $ \underline{J} $  is a
Hopf algebra filtration of  $ \calH \, $;  since  $ \calH $  is graded
connected we have  $ \, J = \calH_+ := \oplus_{n \in \N} \calH_{(n)}
\, $  (where  $ \calH_{(n)} $  is the  $ n $-th  homogeneous component of
$ \calH $),  whence  $ \, \bigcap_{n \in \N} J^n = \{0\} \, $  and  $ \,
\calH^\vee := \calH \Big/ \! \bigcap_{n \in \N} J^n \cong \calH \, $.
We let the  {\sl Rees algebra\/}  associated to  $ \underline{J} $  be
  $$  {\calH_\h}^{\!\vee} \, := \, \Bbbk[\h] \, \cdot {\textstyle
\sum\limits_{n \geq 0}} \; \h^{-n} J^n \, = {\textstyle \sum\limits_{n
\geq 0}} \; \Bbbk[\h] \, \h^{-n} \cdot J^n \, = \, {\textstyle
\sum\limits_{n \geq 0}} \, \Bbbk[\h] \, \big( \h^{-1} \cdot J \,\big)^n
\;  \big( \, \subseteq \calH(\h) \, \big) \, .   \eqno (2.1)  $$
Letting  $ \, J_\h := \text{\sl Ker}\,\big(\epsilon_{\scriptscriptstyle
\calH_\h} \, \colon \calH_\h \! \longrightarrow \Bbbk[\h] \, \big) =
\Bbbk[\h] \cdot J \, $  (the augmentation ideal of  $ \calH_\h \, $)
one has
  $$  {\calH_\h}^{\!\vee} \, = \, {\textstyle \sum_{n \geq 0}}
\; \h^{-n} {J_\h}^{\!n} \, = \, {\textstyle \sum_{n \geq 0}}
\, \big( \h^{-1} J_\h \big)^n  \qquad  \big( \, \subseteq
\calH(\h) \, \big) \, .  $$
   \indent   For all  $ \, n \in \N_\nu \, $,  \, set  $ \, \x_n :=
\h^{-1} \a_n \, $;  \, clearly  $ \, {\calH_\h}^{\!\vee} \, $  is
the  $ \Bbbk[\h] $--subalgebra of  $ \, \calH(\h) \, $  generated by
$ \, J^\vee := \h^{-1} J \, $,  \, hence by  $ \, {\{\x_n\}}_{n \in
\N_\nu} \, $,  \, so  $ \; \displaystyle{ {\calH_\h}^{\!\vee}
\, = \, \Bbbk[\h] \big\langle \{\, \x_n \,|\, n \in \N_\nu \,\}
\big\rangle } \, $.  Moreover,
  $$  \hskip-7pt   \hbox{ $ \eqalign{
   \Delta(\hbox{\bf x}_n)  &  \, = \, \hbox{\bf x}_n \! \otimes \! 1
+ 1 \! \otimes \! \hbox{\bf x}_n + {\textstyle \sum_{m=1}^{n-1}} \,
{\textstyle \sum_{k=1}^m} \; \h^k {\textstyle {{n-m+1} \choose k}}
\, \hbox{\bf x}_{n-m} \otimes P^{(k)}_m(\hbox{\bf x}_*) \, ,
\hskip6pt  \epsilon(\x_n) \, = \, 0  \cr
   S(\hbox{\bf x}_n)  &  \, = \, - \hbox{\bf x}_n -
{\textstyle \sum_{m=1}^{n-1}} \, {\textstyle \sum_{k=1}^m}
\; \h^k {\textstyle {{n-m+1} \choose k}} \, \hbox{\bf x}_{n-m}
\, S \big( P^{(k)}_m(\hbox{\bf x}_*) \big) \, = \,  \cr
   &  \hskip61pt   \, = \, - \hbox{\bf x}_n -
{\textstyle \sum_{m=1}^{n-1}} \, {\textstyle \sum_{k=1}^m} \;
\h^k {\textstyle {{n-m+1} \choose k}} \, S(\hbox{\bf x}_{n-m}) \,
P^{(k)}_m(\hbox{\bf x}_*)  \cr } $ }   \hfill \hskip4,5pt (2.2)  $$
%
for all  $ \, n \in \N_\nu \, $,  due to (1.1).  From this one
sees by hands that the following holds:

\vskip7pt

\proclaim{Proposition 2.1} \, Formulas (2.2) make  $ \, {\calH_\h}^{\!
\vee} = \Bbbk[\h] \big\langle \{\, \x_n, \,|\, n \in \N_\nu \,\}
\big\rangle \, $  into a graded Hopf\/ $ \, \Bbbk[\h] $--algebra,
embedded into  $ \, \calH(\h) := \Bbbk(\h) \otimes_{\Bbbk} \calH \, $
as a graded Hopf subalgebra.  Moreover,  $ {\calH_\h}^{\!\vee} $  is
a deformation of  $ \, \calH $,
   \hbox{for its specialization at  $ \, \h = 1 \, $
is isomorphic to  $ \calH $,  i.e.}
  $$  {\calH_\h}^{\!\vee}{\Big|}_{\h=1} := \, {\calH_\h}^{\!\vee}
\Big/ (\h\!-\!1) \, {\calH_\h}^{\!\vee} \, \cong \, \calH  \quad
\text{via}  \quad  \x_n \!\!\mod (\h\!-\!1) \, {\calH_\h}^{\!\vee}
\, \mapsto \, \a_n  \quad  (\,\forall \;\;  n \in \N_\nu \,)  $$
as graded Hopf algebras over\/  $ \Bbbk \, $.   \qed
\endproclaim

%
%
 \eject

   {\sl  $ \underline{\hbox{\it Remark}} $:}  \; the previous result
shows that  $ \calH_\h $  {\sl is a deformation of  $ \, \calH $,  which
is recovered as specialization (of  $ \calH_\h $)  at}  $ \, \h = 1 \, $.
Next result instead shows that  $ \calH_\h $  {\sl is also a deformation
of  $ U(\L_\nu) $,  recovered as specialization at}  $ \, \, \h = 0
\, $.  Altogether, this gives the top-left horizontal arrow in the
frame  $ (\maltese) $  in the Introduction for  $ \, H = \calH :=
\calH_\nu \, $,  \, with  $ \, \gerg_- = \L_\nu \, $.

\vskip7pt

\proclaim{Theorem 2.1} \,  $ {\calH_\h}^{\!\vee} $  is a QUEA at  $ \,
\h = 0 \, $.  Namely, the specialization limit of  $ {\calH_\h}^{\!
\vee} $  at  $ \, \h \! = \! 0 \, $  is  $ \; {\calH_\h}^{\!\vee}
{\Big|}_{\h=0} := \, {\calH_\h}^{\!\vee} \Big/ \h \, {\calH_\h}^{\!\vee}
\, \cong \, \, U(\L_\nu) \; $  via  $ \; \x_n \mod \h \, {\calH_\h}^{\!
\vee} \mapsto x_n \; $  for all  $ \; n \in \N_\nu \, $,  \, thus
inducing on  $ U(\L_\nu) $  the structure of co-Poisson Hopf algebra
uniquely provided by the Lie bialgebra structure on  $ \L_\nu $  given
by  $ \, \delta(x_n) = \sum_{\ell=1}^{n-1} (\ell+1) \, x_\ell \wedge
x_{n-\ell} \, $  (for
      all  $ \, n \in \N_\nu $)\,\footnote{Hereafter, I use notation
                   $ \, a \wedge b := a \otimes b - b \otimes a \, $.}.
In particular in the diagram  $ \, (\maltese) $  for  $ \, H =
\calH \; (= \calH_\nu) \, $  we have  $ \, \gerg_- = \L_\nu \, $.
                                   \hfill\break
   \indent   Finally, the grading  $ d $  given by  $ \, d(x_n) := 1
\;\, (n \in \N_+) \, $  makes  $ \, {\calH_\h}^{\!\vee}{\Big|}_{\h=0}
\!\! \cong \, U(\L_\nu) \, $  into a graded co-Poisson Hopf algebra;
similarly, the grading  $ \, \partial $  given by  $ \, \partial(x_n)
:= n \;\, (n \in \N_+) \, $  makes  $ \, {\calH_\h}^{\!\vee}
{\Big|}_{\h=0} \!\! \cong \, U(\L_\nu) \, $  into a graded Hopf
algebra and  $ \L_\nu $  into a graded Lie bialgebra.
\endproclaim

\demo{Proof}  First observe that since  $ \, {\calH_\h}^{\!\vee} =
\Bbbk[\h] \, \big\langle \{\, \x_n \,|\, n \in \N_\nu \,\} \big\rangle
\, $  and  $ \, U(\L_\nu) = T(V_\nu) = \Bbbk \, \big\langle \{\, x_n
\,|\, n \in \N_\nu \,\} \big\rangle \, $  mapping $ \; \x_n \mod \h \,
{\calH_\h}^{\!\vee} \mapsto x_n \; $  ($ \, \forall \; n \in \N_\nu
\, $)  does really define an isomorphism  {\sl of algebras}  $ \;
\Phi \, \colon \, {\calH_\h}^{\!\vee} \Big/ \h \, {\calH_\h}^{\!\vee}
\, \cong \, U(\L_\nu) \, $.  Second, formulas (2.2) give
%
%
%
 \vskip-15pt
  $$  \displaylines{
   \Delta(\x_n) \equiv \x_n \otimes 1 + 1 \otimes \x_n  \mod \h \,
\big( {\calH_\h}^{\!\vee} \otimes {\calH_\h}^{\!\vee} \big)  \cr
   \epsilon(\x_n) \equiv 0  \mod \h \, \Bbbk[\h] \, ,  \qquad
S(\x_n) \equiv - \x_n  \mod \h \, {\calH_\h}^{\!\vee}  \cr }  $$
 \vskip-2pt
\noindent
for all  $ \, n \in \N_\nu \, $;  \, comparing with the standard Hopf
structure of  $ U(\L_\nu) $  this shows that  $ \Phi $  is an isomorphism
{\sl of Hopf algebras\/}  too.  Finally, as  $ {\calH_\h}^{\!\vee}
{\Big|}_{\h=0} $  is cocommutative, a Poisson co-bracket is defined
on it by the standard recipe used in quantum group theory, namely
  $$  \eqalign{
   \delta(x_n) \;  &  := \; \big( \h^{-1} \, \big( \Delta(\x_n) -
\Delta^{\text{op}}(\x_n) \big) \big) \mod \h \, \big( {\calH_\h}^{\!
\vee} \otimes {\calH_\h}^{\!\vee} \,\big) \; =   \hskip75pt   \cr
                   &  \phantom{:}= \; {\textstyle \sum_{m=1}^{n-1}} \,
{\textstyle {{n-m+1} \choose 1}} \, x_{n-m} \wedge P_m^{(1)}(x_*) \;
= \; {\textstyle \sum_{\ell=1}^{n-1}} \, (\ell+1) \, x_\ell \wedge
x_{n-\ell}  \qquad  \forall \;\; n \in \N_\nu \; .   \quad
\hskip-3pt  \square  \cr }  $$
\enddemo

\vskip1,5truecm

\centerline {\bf \S \; 3 \ The Drinfeld's deformation
$ \, \big( {\calH_\h}^{\!\vee} \big)' \, $. }

\vskip10pt

  {\bf 3.1 The goal.} \, The second step in the crystal duality
principle is to build a second deformation basing upon the Rees
deformation  $ {\calH_\h}^{\!\vee} $.  This will be a new Hopf
$ \Bbbk[\h] $--algebra  $ \big( {\calH_\h}^{\!\vee} \big)' $,
contained in  $ {\calH_\h}^{\!\vee} $,  which for  $ \, \h = 1 \, $
specializes to  $ \calH $  and for  $ \, \h = 0 \, $  specializes to
$ F[K_+] $,  \, the function algebra of some connected Poisson group
$ K_+ \, $;  \, in other words,  $ \, \big( {\calH_\h}^{\!\vee} \big)'
\Big|_{\h=1} \!\! = \calH \, $  and  $ \, \big( {\calH_\h}^{\!\vee}
\big)'\Big|_{\h=0} \!\! = F[K_+] \, $,  \, the latter meaning that
$ \big( {\calH_\h}^{\!\vee} \big)' $  is a  {\sl quantized function
algebra\/}  (QFA in the sequel).  Therefore  $ \big( {\calH_\h}^{\!\vee}
\big)' $  is a  {\sl quantization\/}  of  $ F[K_+] \, $,  and the quantum
symmetry  $ \calH $  is a deformation of the classical Poisson symmetry
$ F[K_+] \, $.
                                       \par
   In addition, the general theory also describes the relationship
between  $ K_+ $  and the Lie bialgebra  $ \, \gerg_- = \L_\nu \, $  in
\S 2.1, which is  $ \, \L_\nu = \text{\it coLie}(K_+) \, $,  \, so that
we can write  $ \, K_+ = {G_{\!\L_\nu}\phantom{|}}^{\hskip-8pt \star}
\, $.  Comparing with \S 2.1, one eventually concludes that the quantum
symmetry encoded by  $ \calH $  is intermediate between the two classical,
Poisson symmetries ruled by  $ {G_{\!\L_\nu}\phantom{|}}^{\hskip-8pt
\star} $  and  $ \L_\nu \, $.
                                       \par
   In this section I describe explicitly  $ \big( {\calH_\h}^{\!\vee}
\big)' $  and its semiclassical limit  $ F[G_-] \, $,  \, hence  $ G_- $
itself too.  This yields a direct proof of all above mentioned results
about  $ \big( {\calH_\h}^{\!\vee} \big)' $  and  $ G_- \, $.

\vskip7pt

  {\bf 3.2 Drinfeld's  $ \delta_\bullet $--{\,}maps.} \, Let  $ H $  be
any Hopf algebra (over a ring  $ R \, $).  For every  $ \, n \in \N \, $,
\, define the iterated coproduct  $ \; \Delta^n \colon H \longrightarrow
H^{\otimes n} \; $  by $ \, \Delta^0 := \epsilon \, $,  $ \, \Delta^1
:= \id_{\scriptscriptstyle C} $,  \, and finally  $ \, \Delta^n :=
\big( \Delta \otimes \id_{\scriptscriptstyle C}^{\,\otimes (n-2)}
\big) \circ \Delta^{n-1} \, $  if  $ \, n > 2 \, $.  For any ordered
subset  $ \, \Phi = \{i_1, \dots, i_k\} \subseteq \{1, \dots, n\}
\, $  with  $ \, i_1 < \dots < i_k \, $,  \, define the linear map
$ \; j_{\scriptscriptstyle \Phi} : H^{\otimes k} \longrightarrow
H^{\otimes n} \; $  by  $ \; j_{\scriptscriptstyle \Phi} (a_1 \otimes
\cdots \otimes a_k) := b_1 \otimes \cdots \otimes b_n \; $  with  $ \,
b_i := \und1 \, $  if  $ \, i \notin \Phi \, $  and  $ \, b_{i_m}
:= a_m \, $  for  $ \, 1 \leq m \leq k \, $;  \, then set  $ \;
\Delta_\Phi := j_{\scriptscriptstyle \Phi} \circ \Delta^k \, $,
$ \, \Delta_\emptyset := \Delta^0 \, $,  and  $ \; \delta_\Phi
:= \sum_{\Psi \subset \Phi} {(-1)}^{n-|\Psi|} \Delta_\Psi \, $,
$ \; \delta_\emptyset := \epsilon \, $.  The inverse formula
$ \; \Delta_\Phi = \sum_{\Psi \subseteq \Phi} \delta_\Psi \, $
also holds.  We shall also use the shorthand notation  $ \,
\delta_0 := \delta_\emptyset \, $,  $ \, \delta_n := \delta_{\{1,
2, \dots, n\}} \, $  for  $ \, n \in \N_+ \, $.  The following
properties of the maps  $ \delta_\Phi $  will be used:
                                             \par
   {\it (a)} \;  $ \displaystyle{ \delta_n \, = \, \big(
\id_{\scriptscriptstyle C} - u \circ \epsilon \,\big)^{\otimes n}
\circ \Delta^n } \; $ for all  $ \, n \in \N_+ \, $,  \, where
$ \, u : R \longrightarrow H \, $  is the unit map;
                                             \par
   {\it (b)} \;  the maps  $ \, \delta_n \, $  are coassociative,
that is  $ \; \displaystyle{ \Big( \id_{\scriptscriptstyle C}^{\,
\otimes s} \otimes \delta_\ell \otimes \id_{\scriptscriptstyle C}^{\,
\otimes (n-1-s)} \Big) \circ \delta_n \, = \, \delta_{n+\ell-1} } \; $
for all  $ \, n, \ell, s \in \N \, $,  $ \, 0 \leq s \leq n-1 \, $,
\, and similarly in general for the maps  $ \delta_\Phi \, $;
                                             \par
   {\it (c)} \;  $ \displaystyle{ \delta_\Phi(ab) \; = \; {\textstyle
\sum\nolimits_{\Lambda \cup Y = \Phi}} \, \delta_\Lambda(a) \,
\delta_Y(b) } \; $  for all finite subset  $ \, \Phi \subseteq \N \, $
and all  $ \, a, b \in H \, $;
                                             \par
   {\it (d)} \;  $ \displaystyle{ \delta_\Phi(ab - ba) \; = \;
{\textstyle \sum\nolimits_{\Sb  \Lambda \cup Y = \Phi  \\
\Lambda \cap Y \not= \emptyset  \endSb}} \,
\big( \delta_\Lambda(a) \, \delta_Y(b) -
\delta_Y(b) \, \delta_\Lambda(a) \big) } \; $
for all  $ \, \Phi \not= \emptyset \, $  and  $ \, a, b \in H \, $.

\vskip7pt

  {\bf 3.3 Drinfeld's algebra  $ \big( {\calH_\h}^{\!\vee} \big)' $.}
\, Using Drinfeld's  $ \delta_\bullet $--{\,}maps  of \S 3.2, we define
  $$  \big( {\calH_\h}^{\!\vee} \big)' := \Big\{\, \eta \in
{\calH_\h}^{\!\vee} \;\Big\vert\;\, \delta_n(\eta) \in \h^n \big(
{\calH_\h}^{\!\vee} \big)^{\otimes n} \; \forall\; n \in \N \,\Big\}
\qquad  \big( \, \subseteq {\calH_\h}^{\!\vee} \,\big) \; .
\eqno (3.1)  $$
   \indent   Now I describe  $ \big( {\calH_\h}^{\!\vee} \big)' $  and
its specializations at  $ \, \h = 1 \, $  and  $ \, \h = 0 \, $,  \,
in several steps.

\vskip3pt

   {\sl Step I:} \,  {\it A direct check shows that  $ \, \tilde{\x}_n
:= \h \, \x_n = \a_n \in \big( {\calH_\h}^{\!\vee} \big)' $,  \, for
all  $ \, n \in \N_\nu \, $}.  Indeed, we have of course  $ \, \delta_0
(\tilde{\x}_n) = \epsilon(\tilde{\x}_n) \in \h^0 \, {\calH_\h}^{\!\vee}
\, $  and  $ \, \delta_1(\tilde{\x}_n) = \tilde{\x}_n - \epsilon
(\tilde{\x}_n) \in \h^1 \, {\calH_\h}^{\!\vee} \, $.  Moreover,
$ \, \delta_2(\tilde{\x}_n) = \sum_{m=1}^{n-1} \tilde{\x}_{n-m}
\otimes Q_m^{n-m} (\tilde{\x}_*) = \sum_{m=1}^{n-1} \sum_{k=1}^m
\h^{k+1} {{n-m+1} \choose k} \hbox{\bf x}_{n-m} \otimes P^{(k)}_m
(\hbox{\bf x}_*) \in \h^2 \Big( {\calH_\h}^{\!\vee} \otimes
{\calH_\h}^{\!\vee} \Big) \, $.  Since in general  $ \, \delta_\ell
= \big( \delta_{\ell-1} \otimes \id \big) \circ \delta_2 \, $  for
all  $ \, \ell \in \N_+ \, $,  \, we have
  $$  \delta_\ell(\tilde{\x}_n) \, = \, \big( \delta_{\ell-1} \otimes
\id  \big) \big( \delta_2 (\tilde{\x}_n) \big) \, = \, {\textstyle
\sum_{m=1}^{n-1} \sum_{k=1}^m} \; \h^k {\textstyle {{n-m+1} \choose k}} \,
\delta_{\ell-1}(\hbox{\bf x}_{n-m}) \otimes P^{(k)}_m(\hbox{\bf x}_*)  $$  
whence induction gives  $ \, \delta_\ell(\tilde{\x}_n) \in \h^\ell
\, \big( {\calH_\h}^{\!\vee} \big)^{\otimes \ell} \, $  for all  $ \,
\ell \in \N \, $,  \, thus  $ \, \tilde{\x}_n \in \big( {\calH_\h}^{\!
\vee} \big)' \, $,  \, q.e.d.

\vskip3pt

   {\sl Step II:} \, Using property  {\it (c)\/}  in \S 3.2 one
easily checks that  $ \big( {\calH_\h}^{\!\vee} \big)' $  is a
$ \Bbbk[\h] $--subalgebra  of  $ {\calH_\h}^{\!\vee} \, $  (see
[Ga2--3], Proposition 3.5 for details).  In particular, by  {\sl
Step I\/}  and the very definitions this implies that  $ \big(
{\calH_\h}^{\!\vee} \big)' $  contains  $ \calH_\h \, $.

\vskip3pt

   {\sl Step III:} \, Using property  {\it (d)\/}  in \S 3.2 one easily
sees that  $ \big( {\calH_\h}^{\!\vee} \big)'\Big|_{\h=0} $  is
commutative (cf.~[Ga2--3], Theorem 3.8 for details): this means  $ \,
[a,b] \equiv 0 \mod \h \, \big( {\calH_\h}^{\!\vee} \big)' \, $,  \,
that is  $ \, [a,b] \in \h \, \big( {\calH_\h}^{\!\vee} \big)' \, $
hence also  $ \, \h^{-1} [a,b] \in \big( {\calH_\h}^{\!\vee} \big)' \, $,
\, for all  $ \, a $,  $ b \in \big( {\calH_\h}^{\!\vee} \big)' \, $.  In
particular, we get  $ \, \widetilde{[\,\x_n,\x_m]} := \h \, [\,\x_n,\x_m]
= \h^{-1} [\,\tilde{\x}_n,\tilde{\x}_m] \in \big( {\calH_\h}^{\!\vee}
\big)' \, $  for all  $ \, n $, $ m \in \N_\nu \, $,  \, whence iterating
(and recalling  $ \L_\nu $ is generated by the  $ \x_n $'s)  {\it we get
$ \, \tilde{\x} := \h \, \x \in \big( {\calH_\h}^{\!\vee} \big)' \, $
for every  $ \, \x \in \L_\nu \, $}.  Hereafter we identify the free
Lie algebra  $ \L_\nu $  with its image via the natural embedding  $ \,
\L_\nu \lhook\joinrel\longrightarrow U(\L_\nu) = \Bbbk \big\langle
{\{x_n\}}_{n \in \N_\nu} \big\rangle \lhook\joinrel\longrightarrow
\Bbbk[\h] \big\langle {\{\x_n\}}_{n \in \N_\nu} \big\rangle =
{\calH_\h}^{\!\vee} \, $  given by  $ \, x_n \mapsto \x_n \, $
($ \, n \in \N_\nu $) \, .

\vskip3pt

   {\sl Step IV:} \, The previous step showed that, if we
embed  $ \, \L_\nu \lhook\joinrel\longrightarrow U(\L_\nu)
\lhook\joinrel\longrightarrow {\calH_\h}^{\!\vee} \, $  via
$ \, x_n \mapsto \x_n \, $  ($ n \in \N_\nu $)  we find  $ \,
\widetilde{\L_\nu} := \h \, \L_\nu \subseteq \big( {\calH_\h}^{\!
\vee} \big)' \, $.  {\it Let  $ \, \big\langle \widetilde{\L_\nu}
\,\big\rangle \, $  be the  $ \Bbbk[\h] $--subalgebra  of  $ \big(
{\calH_\h}^{\!\vee} \big)' $  generated by  $ \widetilde{\L_\nu} \, $:
then  $ \, \big\langle \widetilde{\L_\nu} \,\big\rangle \subseteq \big(
{\calH_\h}^{\!\vee} \big)' \, $},  \, because  $ \big( {\calH_\h}^{\!
\vee} \big)' $  is a subalgebra.  In particular, if  $ \, \b_b \in
{\calH_\h}^{\!\vee} \, $  is the image of any  $ \, b \in B_\nu \, $
(cf.~\S 1.1) we have  $ \, \widetilde\b_b := \h \, \b_b \in \big(
{\calH_\h}^{\!\vee} \big)' \, $.

\vskip3pt

   {\sl Step V:} \,  {\it Conversely to  {\sl Step IV},  we have
$ \, \big\langle \widetilde{\L_\nu} \,\big\rangle \supseteq \big(
{\calH_\h}^{\!\vee} \big)' \, $}.  In fact, let  $ \, \eta \in \big(
{\calH_\h}^{\!\vee} \big)' \, $;  \, then there are unique  $ \, d
\in \N \, $,  $ \, \eta_+ \in {\calH_\h}^{\!\vee} \setminus \h \,
{\calH_\h}^{\!\vee} \, $  such that  $ \, \eta = \h^d \eta_+ \, $;
\, set also  $ \, \bar{y} := y \mod \h \, {H_\h}^{\!\vee} \in {H_\h}^{\!
\vee} \Big/ \h \, {H_\h}^{\!\vee} \, $  for all  $ \, y \in {H_\h}^{\!\vee}
\, $.  As  $ \, {\calH_\h}^{\!\vee} = \Bbbk[\h] \big\langle \{\, \x_n \,|\,
n \in \N_\nu \,\} \big\rangle \, $  there is a unique  $ \h $--adic
expansion  of  $ \eta_+ $,  namely  $ \; \eta_+ = \eta_0 + \h \, \eta_1
+ \cdots + \h^s \, \eta_s = \sum_{k=0}^s \h^k \, \eta_k \; $  with all  $ \;
\eta_k \in \Bbbk \big\langle \{\, \x_n \,|\, n \in \N_\nu \,\} \big\rangle
\; $  and  $ \; \eta_0 \not= 0 \, $.  Then  $ \, \bar\eta_+ = \bar\eta_0
:= \eta_0 \mod \h \, {\calH_\h}^{\!\vee} \, $,  \, with  $ \, \bar\eta_+
= \bar\eta_0 \in {\calH_\h}^{\!\vee}\Big|_{\h=0} \! = U(\L_\nu) \, $  by
Theorem 2.1.  On the other hand,  $ \, \eta \in \big( {\calH_\h}^{\!\vee}
\big)' \, $  implies  $ \, \delta_{d+1}(\eta) \in \h^{d+1} \big(
{\calH_\h}^{\!\vee} \big)^{\otimes (d+1)} $,  \, whence  $ \,
\delta_{d+1}(\eta_+) = \h^{-d} \delta_{d+1}(\eta) \in \h \, \big(
{\calH_\h}^{\!\vee} \big)^{\otimes (d+1)} \, $  so that  $ \,
\delta_{d+1}\big(\bar\eta_0\big) = 0 \, $;  \, the latter implies
that the degree  $ \partial(\bar\eta_0) $  of  $ \bar\eta_0 $  for
the standard filtration of  $ U(\L_\nu) $  is at most  $ d \, $
(cf.~[Ga2--3],  Lemma 4.2{\it (d)\/}  for a proof).  By the PBW
theorem,  $ \partial(\bar\eta_0) $  is also the degree of  $ \bar\eta_0 $
as a polynomial in the  $ \bar{\text{x}}_b $'s,  hence also of  $ \eta_0 $
as a polynomial in the  $ \text{x}_b $'s  ($ b \in B_\nu $):  \, then
$ \, \h^d \, \eta_0 \in \big\langle \widetilde{\L_\nu} \,\big\rangle
\subseteq \big( {\calH_\h}^{\!\vee} \big)' \, $  (using  {\it Step
III\/}),  hence we find
  $$  \eta_{(1)} \, := \, \h^{d+1} \, \big( \eta_1 + \h \, \eta_2
+ \cdots + \h^{s-1} \, \eta_s \big) = \eta - \h^d \, \eta_0 \in
\big( {\calH_\h}^{\!\vee} \big)' \, .  $$
Thus we can apply our argument again, with  $ \eta_{(1)} $  instead
of  $ \eta $.  Iterating we find  $ \, \partial(\bar\eta_k) \leq
d+k \, $,  \, whence  $ \, \h^{d+k} \, \eta_k \in \big\langle
\widetilde{\L_\nu} \,\big\rangle \, $  $ \Big( \subseteq \big(
{\calH_\h}^{\!\vee} \big)' \Big) $  for all  $ k \, $,
\hbox{thus  $ \, \eta = \sum_{k=0}^s \h^{d+k}
   \, \eta_k \, \in \big\langle \widetilde{\L_\nu} \,\big\rangle \, $,
   \, q.e.d.}

\vskip3pt

   An entirely similar analysis clearly works with  $ \calK_\h $
taking the role of  $ \calH_\h $,  with similar results ({\sl
mutatis mutandis\/}).  On the upshot, we get the following
description:
%
%

\vskip7pt

\proclaim {Theorem 3.1} \, (a) \, With notation of  {\sl Step III\/}
in \S 3.3 (and  $ \, [a,c\,] := a \, c - c \, a \, $),  we have
  $$  \big( {\calH_\h}^{\!\vee} \big)' \; = \; \Big\langle
\widetilde{\L_\nu} \;\Big\rangle \; = \; \Bbbk[\h] \, \Big\langle
{\big\{\, \widetilde\b_b \,\big\}}_{b \in B_\nu} \Big\rangle \Bigg/
\bigg( \Big\{\, \Big[ \widetilde\b_{b_1}, \widetilde\b_{b_2} \Big]
- \h \; \widetilde{\big[\,\b_{b_1},\b_{b_2}\big]} \;\Big|\;
\forall \; b_1, b_2 \in B_\nu \;\Big\} \bigg) \; .  $$
  \indent   (b) \,  $ \big( {\calH_\h}^{\!\vee} \big)' $  is a graded
Hopf\/  $ \, \Bbbk[\h] $--subalgebra  of  $ \, {\calH_\h}^{\!\vee}
\, $,  \, and  $ \, \calH $  is naturally embedded into  $ \, \big(
{\calH_\h}^{\!\vee} \big)' $  as a graded Hopf subalgebra via  $ \,\;
\calH \lhook\joinrel\loongrightarrow \big( {\calH_\h}^{\!\vee}
\big)' \, $,  $ \, \a_n \mapsto \tilde\x_n \; $  (for all
$ \, n \in \N_\nu $).
                                        \hfill\break
  \indent   (c) \,  $ \big( {\calH_\h}^{\!\vee} \big)'{\Big|}_{\h=0} :=
\big( {\calH_\h}^{\!\vee} \big)' \Big/ \h \, \big( {\calH_\h}^{\!\vee}
\big)' = F \big[ {G_{\!\L_\nu}\phantom{|}}^{\hskip-8pt \star} \big]
\, $,  \, where  $ \, {G_{\!\L_\nu}\phantom{|}}^{\hskip-8pt \star}
\, $  is an infinite dimensional connected Poisson algebraic group
with cotangent Lie bialgebra isomorphic to  $ \L_\nu $  (with the
graded Lie bialgebra structure of Theorem 2.1).  Indeed,  $ \big(
{\calH_\h}^{\!\vee} \big)'{\Big|}_{\h=0} $  is the free Poisson
(commutative) algebra over  $ \N_\nu \, $,  generated by all the\/
$ \tilde\x_n{\big|}_{\h=0} $  ($ \, n \in \N_\nu \, $)  with Hopf
structure given by (1.1) with\/  $ \tilde\x_* $  instead of\/
$ \a_* \, $.  Thus  $ \big( {\calH_\h}^{\!\vee} \big)'{\Big|}_{\h=0} $
is the polynomial algebra  $ \, \Bbbk \big[ {\{\, \beta_b \,\}}_{b \in
B_\nu} \big] \, $  generated by a set of indeterminates  $ \, {\{\,
\beta_b \,\}}_{b \in B_\nu} \, $  in bijection with the basis
$ B_\nu $  of\/  $ \L_\nu \, $,  so  $ \; {G_{\!\L_\nu}
\phantom{|}}^{\hskip-8pt \star} \cong \Bbb{A}_\Bbbk^{B_\nu}
\, $  (a (pro)affine\/  $ \Bbbk $--space)  as algebraic varieties.
Finally,  $ \, F \big[ {G_{\!\L_\nu}\phantom{|}}^{\hskip-8pt \star}
\big] = \big( {\calH_\h}^{\!\vee} \big)'{\Big|}_{\h=0} \cong \Bbbk
\big[ {\{\, \beta_b \,\}}_{b \in B_\nu} \big] \, $  bears the natural
{\sl algebra grading}  $ \, d $  of polynomial algebras and the  {\sl
Hopf algebra grading}  inherited from  $ \big( {\calH_\h}^{\!\vee}
\big)' $,  respectively given by  $ \, d\big(\widetilde\b_b\big) =
1 \, $  and  $ \, \partial\big(\widetilde\b_b\big) = \sum_{i=1}^k
n_i \, $  for all  $ \, b = [[\cdots[[x_{n_1},x_{n_2}],x_{n_3}],
\cdots],x_{n_k}] \in B_\nu \, $.
                                        \hfill\break
  \indent   (d) \,  $ F \big[ \G_\nu \big] $  is naturally embedded
into  $ \, \big( {\calH_\h}^{\!\vee} \big)'{\Big|}_{\h=0} = F \big[
{G_{\!\L_\nu}\phantom{|}}^{\hskip-8pt \star} \big] \, $  as a graded
Hopf subalgebra via  $ \; \mu \, \colon \, F \big[ \G_\nu \big]
\lhook\joinrel\loongrightarrow \big( {\calH_\h}^{\!\vee} \big)'
{\Big|}_{\h=0} = F \big[ {G_{\!\L_\nu} \phantom{|}}^{\hskip-8pt
\star} \big] \, $,  $ \, a_n \mapsto \Big( \, \tilde\x_n \! \mod
\h \, \big( {\calH_\h}^{\!\vee} \big)' \Big) \; $  (for all  $ \,
n \in \N_\nu \, $);  moreover,  $ F \big[ \G_\nu \big] $  freely
generates  $ F \big[ {G_{\!\L_\nu}\phantom{|}}^{\hskip-8pt \star}
\big] $  as a\/  {\rm Poisson}  algebra.   Thus there is an
algebraic group epimorphism  $ \, \mu_* \, \colon \,
{G_{\!\L_\nu} \phantom{|}}^{\hskip-8pt \star} \,
{\relbar\joinrel\relbar\joinrel\twoheadrightarrow}
\; \G_\nu \, $,  \, that is  $ \, {G_{\!\L_\nu}
\phantom{|}}^{\hskip-8pt \star} \, $  is
an extension of  $ \, \G_\nu \, $.
                                        \hfill\break
  \indent   (e) \,  Mapping  $ \; \Big( \, \tilde\x_n \! \mod
\h \, \big( {\calH_\h}^{\!\vee} \big)' \Big) \mapsto a_n \; $
(for all  $ \, n \in \N_\nu $)  gives a well-defined graded
Hopf algebra epimorphism  $ \; \pi \, \colon \, F \big[
{G_{\!\L_\nu}\phantom{|}}^{\hskip-8pt \star} \big] \,
{\relbar\joinrel\relbar\joinrel\twoheadrightarrow} \, F
\big[ \G_\nu \big] \, $.  Thus there is an algebraic
group mono\-morphism  $ \, \pi_* \, \colon \, \G_\nu
\lhook\joinrel\loongrightarrow\, {G_{\!\L_\nu}
\phantom{|}}^{\hskip-8pt \star} \, $,  \, that is
$ \, \G_\nu \, $  is an algebraic subgroup of  $ \,
{G_{\!\L_\nu} \phantom{|}}^{\hskip-8pt \star} \, $.
                                        \hfill\break
  \indent   (f) \,  The map  $ \mu $  is a section of  $ \pi $,
hence  $ \pi_* $  is a section of  $ \mu_* \, $.  Thus  $ \,
{G_{\!\L_\nu}\phantom{|}}^{\hskip-8pt \star} \, $  is a
semidirect product of algebraic groups, namely  $ \; {G_{\!\L_\nu}
\phantom{|}}^{\hskip-8pt \star} = \, \G_\nu \ltimes \Cal{N}_\nu \; $
where  $ \, \Cal{N}_\nu := \text{\sl Ker}\,(\mu_*) \trianglelefteq
{G_{\!\L_\nu}\phantom{|}}^{\hskip-8pt \star} \, $.
                                        \hfill\break
  \indent   (g) \,  The analogues of statements (a)--(f) hold
with  $ \calK $  instead of  $ \, \calH \, $,  \, with  $ X^+ $ 
instead of  $ X $  for all  $ \, X = \L_\nu, B_\nu, \N_\nu, \mu,
\pi, \Cal{N}_\nu \, $,  and with
$ {G_{\!\L_\nu^+}\phantom{|}}^{\hskip-8pt \star} $
instead of  $ {G_{\!\L_\nu}\phantom{|}}^{\hskip-8pt \star} $.
\endproclaim

\demo{Proof} {\it (a)} \, This part follows directly from
{\sl Step IV\/ {\rm and}  Step V\/}  in \S 3.3.
                                              \par
   {\it (b)} \, To show that  $ \big( {\calH_\h}^{\!\vee} \big)' $
is a graded Hopf subalgebra we use its presentation in  {\it (a)}.
But first recall that, by  {\sl Step II},  $ \, \calH \, $  embeds
into  $ \big( {\calH_\h}^{\!\vee} \big)' $  via an embedding which
is compatible with the Hopf operations (it is a restriction of the
identity on  $ \calH(\h) \, $):  then this will be a Hopf
algebra monomorphism,  {\sl up to proving that  $ \big( {\calH_\h}^{\!
\vee} \big)' $  is a Hopf subalgebra\/}  (of  $ {\calH_\h}^{\!\vee}
\, $).
                                              \par
   Now,  $ \epsilon_{{\calH_\h}^{\!\vee}} $  obviously restricts to
give a counit for  $ \big( {\calH_\h}^{\!\vee} \big)' $.  Second, we
show that  $ \, \Delta \Big( \big( {\calH_\h}^{\!\vee} \big)' \Big)
\subseteq \big( {\calH_\h}^{\!\vee} \big)' \otimes \big( {\calH_\h}^{\!
\vee} \big)' \, $,  \, so  $ \Delta $  restricts to a coproduct for
$ \big( {\calH_\h}^{\!\vee} \big)' $.  Indeed, each  $ \, b \in B_\nu
\, $  is a Lie monomial, say  $ \, b = [[[ \dots [x_{n_1}, x_{n_2}],
x_{n_3}], \dots], x_{n_k}] \, $  for some  $ \, k $,  $ n_1 $,
$ \dots $,  $ n_k \in \N_\nu \, $,  where  $ k $  is its  Lie
degree: by induction on  $ k \, $  we'll prove  $ \, \Delta \big(
\widetilde\b_b \big) \in \big( {\calH_\h}^{\!\vee} \big)' \otimes
\big( {\calH_\h}^{\!\vee} \big)' \, $  (with  $ \, \widetilde\b_b
:= \h \, \b_b = \h \, [[[ \dots [\x_{n_1}, \x_{n_2}], \x_{n_3}],
\dots], \x_{n_k}] \, $).
                                              \par
   If  $ \, k = 1 \, $  then  $ \, b = x_n \, $  for some  $ \, n
\in \N_\nu \, $.  Then  $ \, \widetilde\b_b = \h \, \x_n = \a_n
\, $  and
  $$  \Delta\Big(\widetilde\b_b\Big) \! = \Delta(\a_n) =
\a_n \otimes 1 + 1 \otimes \a_n + \! {\textstyle \sum\limits_{m=1}^{n-1}}
\a_{n-m} \otimes Q_m^{n-m}(\a_*) \, \in \calH^{\text{dif}} \otimes
\calH^{\text{dif}} \subseteq \big( {\calH_\h}^{\!\vee} \big)' \otimes
\big( {\calH_\h}^{\!\vee} \big)' .  $$
   \indent   If  $ \, k > 1 \, $  then  $ \, b = [b^-,x_n] \, $  for
some  $ \, n \in \N_\nu \, $  and some  $ \, b^- \in B_\nu \, $
expressed by a Lie monomial of degree  $ \, k-1 \, $.  Then  $ \,
\widetilde\b_b = \h \, [\b^-,\x_n] = \Big[ \widetilde\b^-,\x_n \Big]
\, $  and
  $$  \displaylines{
   \Delta\Big(\widetilde\b_b\Big) \; = \; \Delta\Big( \Big[\widetilde\b^-,
\x_n \Big] \Big) \; = \; \Big[ \Delta\Big(\widetilde\b^-\Big),
\Delta(\x_n) \Big] \; = \; \h^{-1} \, \Big[ \Delta\Big(
\widetilde\b^-\Big), \Delta(\a_n) \Big] \; =   \hfill  \cr
   \hfill   = \; \h^{-1} \, \left[\; {\textstyle \sum\nolimits_{\left(
\, \widetilde\b^- \!\right)}} \widetilde\b^-_{(1)} \otimes
\widetilde\b^-_{(2)} \; , \; \a_n \otimes 1 + 1 \otimes \a_n
+ {\textstyle \sum\nolimits_{m=1}^{n-1}} \, \a_{n-m} \otimes
Q_m^{n-m}(\a_*) \;\right] \; =  \cr
%
%
   \hfill   + \; {\textstyle \sum_{\left(\, \widetilde\b^- \!\right)}}
{\textstyle \sum\limits_{m=1}^{n-1}} \bigg( \h^{-1} \! \left[\,
\widetilde\b^-_{(1)} \, , \, \a_{n-m} \right] \! \otimes
\widetilde\b^-_{(2)} \, Q_m^{n-m}(\a_*) \; + \; \widetilde\b^-_{(1)}
\, \a_{n-m} \otimes \h^{-1} \! \left[\, \widetilde\b^-_{(2)} \, ,
\, Q_m^{n-m}(\a_*) \right] \hskip-2pt \bigg)  \cr }  $$
where we used the standard  $ \Sigma $--notation  for  $ \, \Delta
\Big(\widetilde\b^-\Big) = \sum\nolimits_{\left(\, \widetilde\b^-
\!\right)} \widetilde\b^-_{(1)} \otimes \widetilde\b^-_{(2)} \, $.
By inductive
hypothesis we have  $ \, \widetilde\b^-_{(1)} $,  $ \widetilde\b^-_{(2)}
\in \big( {\calH_\h}^{\!\vee} \big)' \, $;  \, then since also  $ \,
\a_\ell \in \big( {\calH_\h}^{\!\vee} \big)' \, $  for all  $ \ell $
and since  $ \big( {\calH_\h}^{\!\vee} \big)' $  is commutative
modulo $ \h $  we have
  $$  \h^{-1} \left[\, \widetilde\b^-_{(1)} \, , \, \a_n \,\right]
\, ,  \; \h^{-1} \left[\, \widetilde\b^-_{(2)} \, , \, \a_n
\,\right] \, ,  \; \h^{-1} \left[\, \widetilde\b^-_{(1)} \, ,
\, \a_{n-m} \,\right] \, ,  \; \h^{-1} \left[\, \widetilde\b^-_{(2)}
\, , \, Q_m^{n-m}(\a_*) \,\right] \, \in \, \big( {\calH_\h}^{\!\vee}
\big)'  $$
for all  $ n $  and  $ (n-m) $  above: so the previous formula
   \hbox{gives  $ \, \Delta\big(\,\widetilde\b_b\big)
\in \big( {\calH_\h}^{\!\vee} \big)' \! \otimes
\big( {\calH_\h}^{\!\vee} \big)' \, $,  \, q.e.d.}
                                              \par
   Finally, the antipode.  Take the Lie monomial  $ \, b = [[[ \dots
[x_{n_1},x_{n_2}], x_{n_3}], \dots], x_{n_k}] \in B_\nu \, $,  \, so
$ \, \widetilde\b_b = \h \, \b_b = \h \, [[[ \dots [\x_{n_1}, \x_{n_2}],
\x_{n_3}], \dots], \x_{n_k}] \, $.  We prove that  $ \, S \big( \,
\widetilde\b_b \big) \in \big( {\calH_\h}^{\!\vee} \big)' \, $  by
induction on the degree  $ k \, $.  If  $ \, k = 1 \, $  then  $ \,
b = x_n \, $  for some  $ n \, $,  so  $ \, \widetilde\b_b = \h \,
\x_n = \a_n \, $  and
  $$  S\big(\,\widetilde\b_b\big) \, = \, S(\a_n) \, = \,
- \a_n - {\textstyle \sum\nolimits_{m=1}^{n-1}} \, \a_{n-m}
\, S\big(Q_m^{n-m}(\a_*)\big) \; \in \; \calH^{\text{dif}}
\subseteq \big( {\calH_\h}^{\!\vee} \big)' \; ,  \quad
\text{q.e.d.}  $$
   \indent   If  $ \, k > 1 \, $  then  $ \, b = [b^-,x_n] \, $  for
some  $ \, n \in \N_\nu \, $  and some  $ \, b^- \in B_\nu \, $  which
is a Lie monomial of degree  $ \, k-1 \, $.  Then  $ \, \widetilde\b_b
= \h \, [\,\b^-,\x_n] = \Big[ \widetilde\b^-, \x_n \Big] = \h^{-1} \,
\Big[ \widetilde\b^-, \a_n \Big] \, $  and so
  $$  S\big(\,\widetilde\b_b\big) \, = \, S\Big( \Big[\,\widetilde\b^-,
\x_n \Big] \Big) \, = \, \h^{-1} \Big[ S(\a_n), S \big( \,
\widetilde\b^- \big) \Big] \; \in \; \h^{-1} \Big[ \big(
{\calH_\h}^{\!\vee} \big)', \big( {\calH_\h}^{\!\vee} \big)'
\Big] \; \subseteq \; \big( {\calH_\h}^{\!\vee} \big)'  $$
using the fact  $ \, S(\a_n) = S\big(\widetilde\x_n\big) =
S\big(\,\widetilde\b_{x_n}\big) \in \big( {\calH_\h}^{\!\vee} \big)'
\, $  (by the case $ \, k \! = \! 1 \, $)  along with the inductive
assumption  $ \, S\big(\,\widetilde\b^-\big) \in \big( {\calH_\h}^{\!
\vee} \big)' \, $  and the commutativity of  $ \big( {\calH_\h}^{\!\vee}
\big)' $  modulo  $ \h \, $.
                                              \par
   {\it (c)} \, As a consequence of  {\it (a)},  the  $ \Bbbk $--algebra
$ \, \big( {\calH_\h}^{\!\vee} \big)'{\Big|}_{\h=0} $  is a  {\sl
polynomial algebra},  namely  $ \; \big( {\calH_\h}^{\!\vee} \big)'
{\Big|}_{\h=0} \! = \, \Bbbk \big[ \{\,\beta_b\,\}_{b \in B}
\big] \; $  with  $ \; \beta_b := \widetilde\b_b \mod \h \, \big(
{\calH_\h}^{\!\vee} \big)' \, $  for all  $ \, b \in B_\nu \, $.  So
$ \big( {\calH_\h}^{\!\vee} \big)'{\Big|}_{\h=0} $  is the algebra of
regular functions  $ F[\varGamma] $ of some (affine) algebraic variety
$ \varGamma \, $;  \, as $ \big( {\calH_\h}^{\!\vee} \big)' $  is a Hopf
algebra the same is true for  $ \, \big( {\calH_\h}^{\!\vee} \big)'
{\Big|}_{\h=0} = F[\varGamma] \, $,  \, so  $ \varGamma $  is an
algebraic group; and since  $ \, F[\varGamma] = \big( {\calH_\h}^{\!
\vee} \big)'{\Big|}_{\h=0} $  is a specialization limit of  $ \big(
{\calH_\h}^{\!\vee} \big)' $,  it is endowed with the Poisson bracket
$ \; \big\{ a|_{\h=0} , b|_{\h=0} \big\} := \big( \h^{-1} [a,b]
\big){\big|}_{\h=0} \; $  which makes  $ \varGamma $  into a
{\sl Poisson\/}  group too.
                                              \par
   We compute the cotangent Lie bialgebra of  $ \varGamma $.  First,
$ \, \germ_e := \text{\sl Ker}\, \big( \epsilon_{F[\varGamma]}
\big) = \Big( {\big\{ \beta_b \big\}}_{b \in B_\nu} \Big) \, $
(the ideal generated by the $ \beta_b $'s)  by construction, so
$ \, {\germ_e}^{\hskip-3pt 2} = \Big( \! {\big\{ \beta_{b_1}
\beta_{b_2} \big\}}_{b_1, b_2 \in B_\nu} \Big) $.  Therefore the
cotangent Lie bialgebra  $ \, Q\big(F[\varGamma]\big) := \germ_e \Big/
{\germ_e}^{\hskip-3pt 2} \, $  as a  $ \Bbbk $--vector space has basis
$ \, {\big\{\overline\beta_b \big\}}_{b \in B_\nu} \, $  where
$ \, \overline\beta_b := \beta_b \mod {{\germ_e}^{\hskip-3pt 2}}
\, $  for all  $ \, b \in B_\nu \, $.  For its Lie bracket
we have (cf.~Remark 1.5)
  $$  \displaylines{
   \big[ \, \overline\beta_{b_1}, \overline\beta_{b_2} \big] \;
:= \; \big\{ \beta_{b_1}, \beta_{b_2} \big\}  \hskip-3pt  \mod
{{\germ_e}}^{\hskip-3pt 2} \; = \; \Big( \h^{-1} \big[\,
\widetilde\b_{b_1}, \widetilde\b_{b_2} \big]  \hskip-3pt  \mod
\h \, \big( {\calH_\h}^{\!\vee} \big)' \Big)  \hskip-3pt  \mod
{{\germ_e}^{\hskip-3pt 2}} \; =   \hfill  \cr
   = \; \Big( \h^{-1} \h^2 \, \big[\,\b_{b_1},\b_{b_2}\big]
\hskip-3pt  \mod \! \h \, \big( {\calH_\h}^{\!\vee} \big)' \Big)
\hskip-3pt  \mod {{\germ_e}^{\hskip-3pt 2}} \; = \, \Big( \h \,
\b_{[b_1,b_2]} \hskip-3pt  \mod \! \h \, \big( {\calH_\h}^{\!\vee}
\big)' \Big)  \hskip-3pt  \mod  {{\germ_e}^{\hskip-3pt 2}} \; =  \cr
   \hfill   = \; \Big( \widetilde\b_{[b_1,b_2]}  \hskip-3pt  \mod
\h \, \big( {\calH_\h}^{\!\vee} \big)' \Big)  \hskip-3pt  \mod
{{\germ_e}^{\hskip-3pt 2}} \; = \; \beta_{[b_1,b_2]}  \hskip-3pt
\mod {{\germ_e}^{\hskip-3pt 2}} \; = \; \overline\beta_{[b_1,b_2]}
\; ,  \cr }  $$
thus the  $ \Bbbk $--linear  map  $ \; \Psi \, \colon \, \L_\nu
\longrightarrow \germ_e \big/ {\germ_e}^{\hskip-3pt 2} \; $
defined by  $ \, b \mapsto \overline\beta_b \, $  for all
$ \, b \in B_\nu \, $  is a  {\sl Lie algebra isomorphism}.  As for
the Lie cobracket, using the general identity  $ \, \delta = \Delta -
\Delta^{\text{op}} \mod \big( {{\germ_e}^{\hskip-3pt 2}} \otimes
F[\varGamma] + F[\varGamma] \otimes {{\germ_e}^{\hskip-3pt 2}} \big)
\, $  (written  $ \hskip-3pt \mod \widehat{{\germ_e}^{\hskip-3pt 2}}
\, $  for short) we get, for all  $ \, n \in \N_\nu \, $,
  $$  \displaylines{
   \delta\big(\,\overline\beta_{x_n}\big) = \big( \Delta \! - \!
\Delta^{\text{op}} \big) (\beta_{x_n})  \hskip-5pt  \mod  \hskip-2pt
\widehat{{\germ_e}^{\hskip-3pt 2}} =  \hskip-1pt  \Big(  \hskip-2pt
\big(  \hskip-1pt  \Delta \! - \! \Delta^{\text{op}}  \hskip-1pt
\big) (\tilde\x_n)  \hskip-5pt \mod \hskip-1,6pt  \h \Big( \! \big(
{\calH_\h}^{\!\vee} \big)' \! \otimes \big( {\calH_\h}^{\!\vee}
\big)' \Big)  \hskip-2pt  \Big)  \hskip-5pt  \mod  \hskip-2pt
\widehat{{\germ_e}^{\hskip-3pt 2}} =   \hfill  \cr
   = \; \left( \left( \a_n \wedge 1 \, + \, 1 \wedge \a_n \, + \,
{\textstyle \sum_{m=1}^{n-1}} \, \a_{n-m} \wedge Q_m^{n-m}(\a_*) \right)
\hskip-3pt  \mod \h \, \big( {\calH_\h}' \otimes {\calH_\h}' \big)
\right)  \hskip-3pt  \mod \widehat{{\germ_e}^{\hskip-3pt 2}} \; =  \cr
%
%
   \hfill   = \! {\textstyle \sum\nolimits_{m=1}^{n-1}}
\beta_{x_{n-m}} \wedge Q_m^{n-m}(\beta_{x_*})  \hskip-3pt
\mod \widehat{{\germ_e}^{\hskip-3pt 2}} = \! {\textstyle
\sum\nolimits_{m=1}^{n-1}} {\textstyle \sum\nolimits_{k=1}^m}
{\textstyle {n-m+1 \choose k}} \, \beta_{x_{n-m}} \wedge
P_m^{(k)}(\beta_{x_*}) \hskip-3pt  \mod
\widehat{{\germ_e}^{\hskip-3pt 2}} =  \cr
   \quad   = \; \left( {\textstyle \sum_{m=1}^{n-1}} {\textstyle
{{n-m+1} \choose 1}} \, \beta_{x_{n-m}} \wedge P_m^{(1)}(\beta_{x_*})
\right)  \hskip-3pt  \mod \widehat{{\germ_e}^{\hskip-3pt 2}} \; =
%
%
%
\; {\textstyle \sum_{\ell=1}^{n-1}} (\ell+1) \; \overline\beta_{x_\ell}
\wedge \overline\beta_{x_{n-\ell}}  \cr }  $$
because   --- among other things ---   one has  $ \; P_m^{(k)}
(\beta_{x_*}) \in {{\germ_e}^{\hskip-3pt 2}} \; $  for all
$ \, k > 1 \, $:  \, therefore
  $$  \delta\big(\,\overline\beta_{x_n}\big) \; = \; {\textstyle
\sum\nolimits_{\ell=1}^{n-1}} (\ell+1) \; \overline\beta_{x_\ell}
\wedge \overline\beta_{x_{n-\ell}}  \qquad  \forall  \hskip7pt
n \in \N_\nu \; .   \eqno (3.2)  $$
Since  $ \L_\nu $  is generated (as a Lie algebra) by the  $ x_n $'s,
the last formula shows that the map  $ \; \Psi \, \colon \, \L_\nu
\longrightarrow \germ_e \big/ {\germ_e}^{\hskip-3pt 2} \; $  given
above is also an isomorphism  {\sl of Lie bialgebras},  \, q.e.d.
                                              \par
   Finally, the statements about gradings of  $ \, \big(
{\calH_\h}^{\!\vee} \big)'{\Big|}_{\h=0} $  should be trivially
clear.
                                              \par
   {\it (d)} \, The part about Hopf algebras is a direct consequence
of  {\it (a)\/}  and  {\it (b)},  noting that the  $ \tilde\x_n $'s
commute modulo  $ \, \h \, \big( {\calH_\h}^{\!\vee} \big)' \, $,
\, since  $ \big( {\calH_\h}^{\!\vee} \big)'{\Big|}_{\h=0} $  is
commutative.
Taking spectra (i.e.~sets of characters of each
Hopf algebras) we get an algebraic group morphism
$ \, \mu_* \, \colon \, {G_{\!\L_\nu}\phantom{|}}^{\hskip-8pt \star}
\,{\relbar\joinrel\relbar\joinrel\rightarrow}\; \G_\nu \, $,  \, which
in fact is  {\sl onto\/}  because, as these algebras are polynomial,
each character of  $ \, F \big[ \G_\nu \big] \, $  does extend to a
character  of  $ \, F \big[ {G_{\!\L_\nu}\phantom{|}}^{\hskip-8pt \star}
\big] \, $,  so the former arises from restriction of the latter.
                                              \par
   {\it (e)} \, Due to the explicit description of  $ F \big[
{G_{\!\L_\nu} \phantom{|}}^{\hskip-8pt \star} \big] $  coming
from  {\it (a)\/}  and  {\it (b)},  mapping  $ \; \Big( \,
\tilde\x_n \! \mod \h \, \big( {\calH_\h}^{\!\vee} \big)' \Big)
\mapsto a_n \; $  (for all  $ \, n \in \N_\nu $)  clearly yields
a Hopf algebra epimorphism  $ \; \pi \, \colon \, F \big[
{G_{\!\L_\nu} \phantom{|}}^{\hskip-6pt \star} \big] \,
{\relbar\joinrel\relbar\joinrel\twoheadrightarrow} \,
F \big[ \G_\nu \big] \, $.  Taking spectra gives an
algebraic group monomorphism  $ \, \pi_* \, \colon \,
\G_\nu \,\lhook\joinrel\loongrightarrow\, {G_{\!\L_\nu}
\phantom{|}}^{\hskip-8pt \star} \, $  as required.
                                              \par
   {\it (f)} \, The map  $ \mu $  is a section of  $ \pi $  by
construction.  Then clearly  $ \pi_* $  is a section of  $ \mu_* \, $,
which implies  $ \; {G_{\!\L_\nu}\phantom{|}}^{\hskip-8pt \star}
= \, \G_\nu \ltimes \Cal{N}_\nu \; $  (with  $ \, \Cal{N}_\nu
:= \text{\sl Ker}\,(\mu_*) \trianglelefteq {G_{\!\L_\nu}
\phantom{|}}^{\hskip-8pt \star} \, $)  by general theory.
                                              \par
   {\it (g)} \, This ought to be clear from the whole discussion,
for all arguments apply again   --- {\sl mutatis mutandis\/} ---
when starting with  $ \calK $  instead of  $ \, \calH \, $;  \,
details are left to the reader.   \qed
\enddemo

\vskip3pt

   {\sl $ \underline{\hbox{\it Remark}} $:}  \; Roughly speaking,
we can say that the extension  $ \, F \big[ \G_\nu \big]
\,\lhook\joinrel\loongrightarrow\, F \big[ {G_{\!\L_\nu}
\phantom{|}}^{\hskip-8pt \star} \big] \, $  is performed simply
by adding to  $ F \big[ \G_\nu \big] $  a  {\sl free\/}  Poisson
structure, which happens to be compatible with the Hopf structure.
Then the Poisson bracket starting from the ``elementary'' coordinates
$ a_n $  (for  $ \, n \in \N_\nu $)  freely generates  {\sl new\/}
coordinates  $ \{a_{n_1},a_{n_2}\} $,  $ \big\{\! \{a_{n_1},
a_{n_2}\}, a_{n_3} \big\} $,  etc., thus enlarging  $ F \big[
\G_\nu \big] $  and generating  $ F \big[ {G_{\!\L_\nu}
\phantom{|}}^{\hskip-8pt \star} \big] $.  At the group level,
this means that  $ \G_\nu $  {\sl freely Poisson-generates\/}
the Poisson group  $ {G_{\!\L_\nu}\phantom{|}}^{\hskip-8pt \star}
\, $:  \, new 1-parameter subgroups, build up in a  {\sl Poisson-free\/}
manner from those attached to the  $ a_n $'s,  are freely ``pasted''
to  $ \G_\nu \, $,  expanding it and building up  $ {G_{\!\L_\nu}
\phantom{|}}^{\hskip-8pt \star} \, $.  Then the epimorphism  $ \,
{G_{\!\L_\nu}\phantom{|}}^{\hskip-8pt \star} \, {\buildrel \mu_*
\over {\relbar\joinrel\relbar\joinrel\twoheadrightarrow}} \; \G_\nu
\, $  is just a forgetful map: it kills the new 1-parameter subgroups
and is injective (hence an isomorphism) on the subgroup generated by
the old ones.  On the other hand, definitions imply that  $ \, F \big[
{G_{\!\L_\nu}\phantom{|}}^{\hskip-8pt \star} \big] \Big/ \Big( \big\{
F \big[ {G_{\!\L_\nu}\phantom{|}}^{\hskip-8pt \star} \big], F \big[
{G_{\!\L_\nu} \phantom{|}}^{\hskip-8pt \star}\big] \big\} \Big) \cong
F \big[ \G_\nu \big] \, $,  \, and with this identification  $ \, F
\big[ {G_{\!\L_\nu}\phantom{|}}^{\hskip-8pt \star} \big] \, {\buildrel
\pi \over {\relbar\joinrel\relbar\joinrel\twoheadrightarrow}} \; F
\big[ \G_\nu \big] \, $  is just the canonical map, which mods out
all Poisson brakets  $ \{f_1,f_2\} $,  \, for  $ \, f_1, f_2 \in
F \big[ {G_{\!\L_\nu}\phantom{|}}^{\hskip-8pt \star} \big] \, $.

\vskip7pt

   {\bf 3.4 Specialization limits.} \, So far, we have already
pointed out (by Proposition 2.1, Theorem 2.1,  Theorem 3.1{\it
(c)\/})  the following specialization limits of  $ {\calH_\h}^{\!
\vee} $  and  $ \big( {\calH_\h}^{\!\vee} \big)' \, $:
  $$  {\calH_\h}^{\!\vee} \;{\buildrel \h \rightarrow 1 \over
\llongrightarrow}\; \calH \; ,  \qquad  {\calH_\h}^{\!\vee}
\;{\buildrel \h \rightarrow 0 \over \llongrightarrow}\; U(\L_\nu) \; ,
\qquad  \big( {\calH_\h}^{\!\vee} \big)' \;{\buildrel \h \rightarrow 0
\over \llongrightarrow}\; F \big[ {G_{\!\L_\nu} \phantom{|}}^{\hskip-8pt
\star} \big]  $$
as graded Hopf  $ \Bbbk $--algebras,  with some (co-)Poisson structures
in the last two cases.  As for the specialization limit of  $ \big(
{\calH_\h}^{\!\vee} \big)' $  at  $ \, \h = 1 \, $,  \, Theorem 3.1
implies that it is $ \calH \, $.  Indeed, by  Theorem 3.1{\it (b)\/}
$ \, \calH $  embeds into  $ \big( {\calH_\h}^{\!\vee} \big)' $  via
$ \; \a_n \mapsto \tilde\x_n \; $  (for all  $ \, n \in \N_\nu $):
then
  $$  [\a_n,\a_m] \, = \, \big[ \tilde\x_n,\tilde\x_m \big] \, = \,
\h \, \widetilde{[\x_n,\x_m]} \, \equiv \, \widetilde{[\x_n,\x_m]}
\mod (\h \! - \! 1) \, \big( {\calH_\h}^{\!\vee} \big)'   \eqno
\big(\, \forall \; n, m \in \N_\nu \big)  $$
whence, due to the presentation of  $ \big( {\calH_\h}^{\!\vee}
\big)' $  by generators and relations in  Theorem 3.1{\it (a)\/},
  $$  \big( {\calH_\h}^{\!\vee} \big)'{\Big|}_{\h=1} \; := \;
\big( {\calH_\h}^{\!\vee} \big)' \Big/ (\h \! - \! 1) \,
\big( {\calH_\h}^{\!\vee} \big)' \; = \; \Bbbk \big\langle
\overline{\tilde\x}_1, \overline{\tilde\x}_2, \dots,
\overline{\tilde\x}_n, \ldots \big\rangle \; = \;
\Bbbk \big\langle \overline\a_1, \overline\a_2, \dots,
\overline\a_n, \ldots \big\rangle  $$
(where  $ \, \overline{\bold{c}} := \bold{c} \mod (\h\!-\!1) \,
\big( {\calH_\h}^{\!\vee} \big)' \, $)  \, as  $ \Bbbk $--algebras,
and the Hopf structure is exactly the one of  $ \calH $  because it
is given by the like formulas on generators.  In a nutshell, we have
$ \; \big( {\calH_\h}^{\!\vee} \big)' \,{\buildrel {\h \rightarrow 1}
\over \llongrightarrow}\; \calH \; $  as Hopf  $ \Bbbk $--algebras.
This completes the top part of the diagram  ($ \maltese $) in the
Introduction, for  $ \, H = \calH \; (:= \calH_\nu) \, $,  \,
because  $ \, \calH^\vee := \calH \big/ \cap_{n \in \N} J^n =
\calH \, $  by \S 2.2: namely,
%
%
%
  $$  U(\L_\nu) \underset{{\calH_\h}^{\!\!\vee}}
\to {\overset{0 \leftarrow \h \rightarrow 1} \to
{\longleftarrow\joinrel\relbar\joinrel\relbar\joinrel\llongrightarrow}}
\, \calH \underset{({\calH_\h}^{\!\!\vee})'}  \to
{\overset{1 \leftarrow \h \rightarrow 0} \to
{\longleftarrow\joinrel\relbar\joinrel\relbar\joinrel\llongrightarrow}}
F \big[ {G_{\!\L_\nu}\phantom{|}}^{\hskip-8pt \star} \big] \quad .  $$
%
%

\vskip1,5truecm

\centerline {\bf \S \; 4 \ The Rees deformation
$ \, {\calH_\h}^{\!\prime} \, $. }

\vskip10pt

  {\bf 4.1 The goal.} \, The crystal duality principle (cf.~[Ga2],
[Ga4]) yields also a recipe to produce a 1-parameter deformation
$ \, {\calH_\h}^{\!\prime} \, $  of  $ \calH $  which is a  {\sl
quantized function algebra\/} (QFA in the sequel): namely,  $ \,
{\calH_\h}^{\!\prime} \, $  is a Hopf  $ \Bbbk[\h] $--algebra  such
that  $ \, {\calH_\h}^{\!\prime}\Big|_{\h=1} \!\! = \calH \, $  and
$ \, {\calH_\h}^{\!\prime}\Big|_{\h=0} \!\! = F[G_+] \, $,  \, the
function algebra of a connected algebraic Poisson group  $ G_+ \, $.
Thus
 \eject   
\noindent   $ {\calH_\h}^{\!\prime} $  is a  {\sl quantization\/} of
$ F[G_+] $,  and the quantum symmetry  $ \calH $  is a deformation
of the classical Poisson symmetry  $ F[G_+] $.  By definition
$ {\calH_\h}^{\!\prime} $  is the  {\sl Rees algebra\/}  associated
to a distinguished  {\sl increasing\/}  Hopf algebra filtration of
$ \calH $,  and  $ F[G_+] $  is simply the graded Hopf algebra
associated to this filtration.  The purpose of this section is
to describe explicitly  $ {\calH_\h}^{\!\prime} $  and its
semiclassical limit  $ F[G_+] $,  hence also  $ G_+ $  itself.
This will also provide a direct, independent proof of all the
above mentioned results about  $ {\calH_\h}^{\!\prime} $  and
$ F[G_+] $  themselves.

\vskip7pt

   {\bf 4.2 The Rees algebra  $ \, {\calH_\h}^{\!\prime} \, $.} \,
Let's consider Drinfeld's  $ \delta_\bullet $--maps,  as in \S 3.2,
for the Hopf algebra  $ \calH \, $.  Using them, we define the  {\sl
$ \delta_\bullet $--filtration\/}  $ \; \underline{D} := {\big\{ D_n
\big\}}_{n \in \N} \; $  of  $ \calH $  by  $ \; D_n := \text{\it Ker}\,
(\delta_{n+1}) \, $,  \; for all  $ \, n \in \N \, $.
  It is easy to show (cf.~[Ga4]) that  $ \underline{D} $  is a Hopf
algebra filtration of  $ \calH \, $;  \, moreover, since  $ \calH $
is graded connected,
%
%
we have  $ \, \calH = \bigcup_{n \in \N} D_n =: \calH' \, $.  We define
the  {\sl Rees algebra\/}  associated to  $ \underline{D} $  as
  $$  {\calH_\h}^{\!\prime} \, := \; \Bbbk[\h] \cdot {\textstyle
\sum\limits_{n \geq 0}} \; \h^{+n} D_n \, = \, {\textstyle
\sum\limits_{n \geq 0}} \; \Bbbk[\h] \h^{+n} \cdot D_n  \qquad
\big( \, \subseteq \calH_\h := \calH[\h] \, \big) \; .   \eqno (4.1)  $$
A trivial check shows that the following intrinsic characterization
(inside  $ \calH_\h $)  also holds:
  $$  {\calH_\h}^{\!\prime} \, = \, \big\{\, \eta \in \calH_\h
\;\big\vert\;\, \delta_n(\eta) \in \h^n {\calH_\h}^{\!\otimes n} ,
\;\; \forall\; n \in \N \,\big\}  \qquad  \big( \, \subseteq \calH_\h
\, \big) \; .  $$
   \indent   We shall describe  $ {\calH_\h}^{\!\prime} $  explicitly,
and we'll compute its specialization at  $ \, \h = 0 \, $  and at  $ \,
\h = 1 \, $:  in particular we'll show that it is really a QFA and a
deformation of  $ \calH \, $,  as claimed.
                                               \par
   By (4.1), all we need is to compute the filtration  $ \, \underline{D}
= {\big\{ D_n \big\}}_{n \in \N} \, $;  \, the idea is to describe it in
combinatorial terms, based on the non-commutative polynomial nature
of  $ \calH \, $.

\vskip7pt

{\bf 4.3 Gradings and filtrations:} \, Let  $ \, \partial_- $  be
the unique Lie algebra grading of  $ \L_\nu $  given by  $ \;
\partial_-(x_n) := n - 1 + \delta_{n,1} \; $  (for all
$ \, n \in \N_\nu \, $).  Let also  $ d $  be the standard Lie algebra
grading associated with the central lower series of  $ \L_\nu \, $,
\, i.e.~the one defined by  $ \; d \big( [ \cdots [[x_{s_1}, x_{s_2}],
\dots x_{s_k}] \big) = k - 1 \; $  on any Lie monomial of  $ \L_\nu \, $.
As both  $ \partial_- $  and  $ d $  are Lie algebra gradings,  $ \,
(\partial_- - d) \, $  is a Lie algebra grading too.  Let  $ \,
{\big\{ F_n \big\}}_{n \in \N} \, $  be the Lie algebra filtration
associated with  $ (\partial_- - d \,) $;  then the down-shifted
filtration  $ \, \underline{T} := {\big\{ \, T_n := F_{n-1} \,
\big\}}_{n \in \N} \, $  is again a Lie algebra filtration of
$ \L_\nu \, $.  There is a unique algebra filtration on  $ U(\L_\nu) $
extending  $ \underline{T} $:  we denote it  $ \, \underline{\varTheta}
= {\big\{ \varTheta_n \big\}}_{n \in \N} \, $,  \, and set also  $ \,
\varTheta_{-1} := \{0\} \, $.  Finally, for each  $ \, y \in U(\L_\nu)
\setminus \{0\} \, $  there is a unique  $ \, \tau(y) \in \N \, $  with
$ \, y \in \varTheta_{\tau(y)} \setminus \varTheta_{\tau(y)-1} \, $;  \,
in particular  $ \, \tau(b) = \partial_-(b) - d(b) \, $,  $ \, \tau(b
\, b') = \tau(b) + \tau(b') \; $  and  $ \; \tau\big([b,b']\big) \! =
\tau(b) + \tau(b') - 1 \; $  for  $ \, b, b' \! \in \! B_\nu \, $.
                                               \par
   We can explicitly describe  $ \underline{\varTheta} $.  Indeed, let
us fix any total order  $ \preceq $  on the basis  $ B_\nu $  of \S
1.1: then  $ \; \calX := \Big\{\, \underline{b} := b_1 \cdots b_k
\,\Big|\; k \in \N \, , \; b_1, \dots, b_k \in B_\nu \, , \; b_1
\preceq \cdots \preceq b_k \,\Big\} \; $  is a  $ \Bbbk $--basis  of
$ U(\L_\nu) $,  by the PBW theorem.  It follows that  $ \varTheta $
induces a set-theoretic filtration  $ \, \underline{\calX} = {\big\{
\calX_n \big\}}_{n \in \N} \, $  of  $ \calX $  with  $ \, \calX_n :=
\calX \cap \varTheta_n = \Big\{\, \underline{b} := b_1 \cdots b_k
\,\Big|\; k \in \N \, , \; b_1, \dots, b_k \in B_\nu \, , \; b_1
\preceq \cdots \preceq b_k \, , \; \tau(\underline{b}\,) = \tau(b_1)
+ \cdots + \tau(b_k) \leq n \,\Big\} \, $,  \; and also that  $ \;
\varTheta_n = \text{\sl Span}\,\big(\calX_n\big) \; $  for all
$ \, n \in \N \, $.
                                               \par
   Let us define  $ \, \boldalpha_1 := \a_1 \, $  and  $ \,
\boldalpha_n := \a_n - {\a_1}^{\!n} \, $  for all  $ \, n \in
\N_\nu \setminus \{1\} \, $.  This ``change of variables''   ---
which switch from the  $ \a_n $'s  to their  {\sl differentials},
in a sense ---   is the key to achieve a complete description
of  $ \underline{D} \, $,  \, via a close comparison between
$ \calH $  and  $ U(\L_\nu) \, $.
                                               \par
   By definition  $ \, \calH = \calH_\nu \, $  is the free associative
algebra over  $ \{\a_n\}_{n \in \N_\nu} \, $,  \, hence (by definition
of the  $ \boldalpha $'s)  also over  $ \{\boldalpha_n\}_{n \in
\N_\nu} $;  so we have an algebra isomorphism  $ \; \Phi \,
\colon \, \calH \,{\buildrel \cong \over
{\lhook\joinrel\relbar\joinrel\twoheadrightarrow}}\,
U(\L_\nu) \; $  given by  $ \, \boldalpha_n \mapsto x_n \, $
($ \, \forall \; n \in \N_\nu \, $).  Via  $ \Phi $  we pull back
all data and results about gradings, filtrations, PBW bases and so
on mentioned above for  $ U(\L_\nu) \, $;  in particular we set  $ \,
\boldalpha_{\underline{b}} := \Phi(x_{\underline{b}}) = \boldalpha_{b_1}
\cdots \boldalpha_{b_k} \, $  ($ \, b_1, \dots, b_k \in B_\nu $),  $ \,
\calA_n := \Phi(\calX_n) \, $  ($ n \in \N $),  $ \, \calA := \Phi(\calX)
= \bigcup_{n \in \N} \calA_n \, $.   For gradings on  $ \calH $  we stick
to the like notation, i.e.~$ \partial_- $,  $ d $  and  $ \tau \, $,  \,
and similarly for  $ \underline{\varTheta} \, $.
                                                    \par
   Finally, for all  $ \, a \in \calH \setminus \{0\} \, $  \, we set
$ \, \kappa\,(a) := k \, $  iff  $ \, a \in D_k \setminus D_{k-1} \, $
(with  $ \, D_{-1} := \{0\} \, $).
                                                    \par
   {\sl Our goal is to prove an identity of filtrations, namely  $ \,
\underline{D} = \underline{\varTheta} \, $},  \, or equivalently
$ \, \kappa = \tau \, $.  In fact, this would give to the Hopf
filtration  $ \underline{D} $,  which is defined intrinsically
in Hopf algebraic terms, an explicit  {\sl combinatorial\/}
description, namely the one of  $ \varTheta $  explained above.

\vskip7pt

\proclaim{Lemma 4.1}  \hbox{$ Q^\ell_t(\a_*) \in
\varTheta_t \! \setminus \varTheta_{t-1} \, ,  \;
Z^\ell_t(\boldalpha_*) := \! \Big( Q^\ell_t(\a_*) - \!
{\textstyle {{\ell + t} \choose t}} \, {\a_1}^{\!t} \Big)
\! \in \varTheta_{t-1} \, (\, \ell, t \in \N \, , t \!
\geq \! 1 ) $.}
\endproclaim

%
%

\demo{Proof} \,  When  $ \, t = 1 \, $  definitions give
$ \, Q^\ell_1(\a_*) = (\ell+1) \, \a_1 \in \varTheta_1 \, $  and
so  $ \, Z^\ell_1(\boldalpha_*) = (\ell+1) \, \a_1 - {{\ell + 1}
\choose 1} \, \a_1 = 0 \in \varTheta_0 \, $,  \, for all  $ \,
\ell \in \N \, $.  Similarly, when  $ \, \ell = 0 \, $  we have
$ \, Q^0_t(\a_*) = \a_t \in \varTheta_t \, $  and so  $ \,
Z^0_t(\boldalpha_*) = \a_t - {1 \choose 1} \, {\a_1}^{\!t}
= \boldalpha_t \in \varTheta_{t-1} \, $  (by definition),
\, for all  $ \, t \in \N_+ \, $.
                                                  \par
   When  $ \, \ell > 0 \, $  and  $ \, t > 1 \, $,  \, we can prove
the claim using two independent methods.
                                                  \par
   {\it $ \underline{\text{First method}} $:} \,  The very
definitions imply that the following recurrence formula holds:
  $$  Q^\ell_t(\a_*) \, = \, Q^{\ell-1}_t(\a_*) \, + \, {\textstyle
\sum\nolimits_{s=1}^{t-1}} \, Q^{\ell-1}_{t-s}(\a_*) \cdot \a_s +
\a_t  \qquad  \forall \quad \ell \geq 1 \, ,  \; t \geq 2 \; .  $$
From this formula and from the identities  $ \, \a_1 = \boldalpha_1 \, $,
$ \, \a_s = \boldalpha_s + {\boldalpha_1}^{\!s} \, $  ($ s \in \N_+ $),
we argue
  $$  \displaylines{
   Z^\ell_t(\boldalpha_*) \; := \; Q^\ell_t(\a_*) \, - \,
{\textstyle {{\ell+t} \choose t}} \, {\a_1}^{\!t} \; =
\; Q^{\ell-1}_t(\a_*) \, + \, {\textstyle \sum_{s=1}^{t-1}}
Q^{\ell-1}_{t-s}(\a_*) \, \a_s \, + \, \a_t \, - \, {\textstyle
{{\ell+t} \choose t}} \, {\a_1}^{\!t} \; =   \hfill  \cr
   {} \hfill   = \; Z^{\ell-1}_t(\a_*) \, + \, {\textstyle {{\ell-1+t}
\choose t}} \, {\a_1}^{\!t} \, + \, {\textstyle \sum_{s=1}^{t-1}}
\left( Z^{\ell-1}_{t-s}(\a_*) \, + \, {\textstyle {{\ell-1+t-s}
\choose {t-s}}} \, {\a_1}^{\!t-s} \right) \, \a_s + \a_t -
{\textstyle {{\ell+t} \choose t}} \, {\a_1}^{\!t} \; =  \cr
%
%
   = \; Z^{\ell-1}_t(\a_*) \, + \, {\textstyle \sum_{s=1}^{t-1}}
Z^{\ell-1}_{t-s}(\a_*) \, \big( \boldalpha_s + {\boldalpha_1}^{\!s}
\big) \, + \, {\textstyle \sum_{s=1}^{t-1}} {\textstyle
{{\ell-1+t-s} \choose {t-s}}} \, {\boldalpha_1}^{\!t-s} \,
\boldalpha_s \, + \, \boldalpha_t \, + \,   \hfill  \cr
   {} \hfill   + \, {\textstyle \sum_{s=1}^{t-1}} {\textstyle
{{\ell-1+t-s} \choose {t-s}}} \, {\boldalpha_1}^{\!t-s} \,
{\boldalpha_1}^{\!s} \, + \, {\boldalpha_1}^{\!t} \,
+ \, {\textstyle {{\ell-1+t} \choose t}}
\, {\boldalpha_1}^{\!t} \, - \,
{\textstyle {{\ell+t} \choose t}} \, {\boldalpha_1}^{\!t} \; =  \cr
   = \; Z^{\ell-1}_t(\a_*) \, + \, {\textstyle \sum_{s=1}^{t-1}}
Z^{\ell-1}_{t-s}(\a_*) \, \big( \boldalpha_s + {\boldalpha_1}^{\!s}
\big) \, + \, {\textstyle \sum_{s=1}^{t-1}} {\textstyle
{{\ell-1+t-s} \choose {t-s}}} \, {\boldalpha_1}^{\!t-s} \,
\boldalpha_s \, + \, \Big( {\textstyle \sum_{r=0}^t} {\textstyle
{{\ell-1+r} \choose {\ell-1}}} \, - \,   \hfill  \cr
   {} \hfill   \, - \, {\textstyle {{\ell+t} \choose t}} \! \Big) \,
{\boldalpha_1}^{\!t} \, + \, \boldalpha_t = \, Z^{\ell-1}_t(\a_*) \,
+ \, {\textstyle \sum_{s=1}^{t-1}} Z^{\ell-1}_{t-s}(\a_*) \, \big(
\boldalpha_s + {\boldalpha_1}^{\!s} \big) \, + \, {\textstyle
\sum_{s=1}^{t-1}} {\textstyle {{\ell-1+t-s} \choose {t-s}}} \,
{\boldalpha_1}^{\!t-s} \, \boldalpha_s \, + \, \boldalpha_t  \cr
%
%
 }  $$
because of the classical identity  $ \; {{\ell+t} \choose \ell} =
\sum_{r=0}^t {{\ell-1+r} \choose {\ell-1}} \; $.  Then induction
upon  $ \ell \, $  and the very definitions allow to argue that
all summands in the final sum belong to  $ \varTheta_{t-1} $,
hence  $ \, Z^\ell_t(\boldalpha_*) \in \varTheta_{t-1} \, $
as well.  Finally, this implies  $ \; Q^\ell_t(\a_*) \, =
\, Z^\ell_t (\boldalpha_*) \, + \, {{\ell+t} \choose t} \,
{\boldalpha_1}^{\!t} \in \varTheta_t \setminus \varTheta_{t-1} \, $.
                                                  \par
   {\it $ \underline{\text{Second method}} $:}  $ \; Q_t^\ell(\a_*)
:= \sum_{s=1}^t {{\ell + 1} \choose s} \, P_t^{(s)}(\a_*) =
\sum_{s=1}^t {{\ell + 1} \choose s} \sum_{\hskip-4pt  \Sb  j_1,
\dots, j_s > 0  \\   j_1 + \cdots + j_s \, = \, t  \endSb}
\hskip-5pt  \a_{j_1} \cdots \a_{j_s} \, $,  \;  by definition;
then expanding the  $ \a_j $'s
%
%
(for  $ \, j > 1 \, $) as above
we find that  $ \, Q_t^\ell(\a_*) = Q_t^\ell \big(
\boldalpha_* + {\boldalpha_1}^{\!*} \big) \, $  is a linear combination
of monomials  $ \, \boldalpha_{(j_1)} \cdots \boldalpha_{(j_s)} \, $
with  $ \, j_1, \dots, j_s > 0 \, $,  $ \, j_1 + \cdots + j_s =
t \, $,  $ \, \boldalpha_{(j_r)} \in \big\{ \boldalpha_{j_r},
{\boldalpha_1}^{\!j_r} \big\} \, $  for all  $ r \, $.  Let
$ Q_- $  be the linear combination of those monomials such
that  $ \, (\boldalpha_{(j_1)}, \boldalpha_{(j_2)}, \dots,
\boldalpha_{(j_s)}\big) \not= \big( {\boldalpha_1}^{\!j_1},
{\boldalpha_1}^{\!j_2}, \dots, {\boldalpha_1}^{\!j_s}\big) \, $;
\, the remaining monomials enjoy  $ \, \boldalpha_{j_1} \cdot
\boldalpha_{j_2} \cdots \boldalpha_{j_s} = {\boldalpha_1}^{\!
j_1 + \cdots + j_s} = {\boldalpha_1}^{\!t} \, $,  \, so their
linear combination giving  $ \, Q_+ := Q^\ell_t(\a_*) - Q_- \, $
is a multiple of  $ {\boldalpha_1}^{\!t} $,  say  $ \, Q_+ = N
\, {\boldalpha_1}^{\!t} \, $.  Now we compute this  $ N \, $.
                                             \par
   By construction,  $ N $  is nothing but  $ \, N = Q^\ell_t(1_*)
= Q^\ell_t(1, 1, \dots, 1, \dots) \, $  where the latter is the
value of  $ Q^\ell_t $  when all indeterminates are set equal
to  $ 1 $;  thus we compute  $ Q^\ell_t(1_*) \, $.
                                             \par
   Recall that the  $ Q^\ell_t $'s  enter in the definition of the
coproduct of  $ F \big[ \G^\dif \big] $:  \, the latter is dual
to the (composition) product of series in  $ \G^\dif $,  thus
if  $ \{a_n\}_{n \in \N_+} $  and  $ \{b_n\}_{n \in \N_+} $
are two countable sets of commutative indeterminates then
  $$  \displaylines{
   \Big( \, x \, + {\textstyle \sum_{n=1}^{+\infty}} \, a_n \,
x^{n+1} \Big) \circ \Big( \, x \, + {\textstyle \sum_{m=1}^{+\infty}}
\, b_m \, x^{m+1} \Big) \; :=   \hfill  \cr
   {} \hfill   = \; \bigg( \! \Big( \, x \, +
{\textstyle \sum\nolimits_{m=1}^{+\infty}} \, b_m \, x^{m+1} \Big)
+ {\textstyle \sum\nolimits_{n=1}^{+\infty}} \, a_n \, \Big( \, x \, +
{\textstyle \sum\nolimits_{m=1}^{+\infty}} \, b_m \, x^{m+1} \Big)^{n+1}
\bigg) \; = \; x \, + {\textstyle \sum\nolimits_{k=0}^{+\infty}} \,
c_k \, x^{k+1}  \cr }  $$
with  $ \; c_k = Q^0_k(b_*) + \sum_{r=1}^k a_r \cdot Q^r_{k-r}(b_*)
\, $  (cf.~\S 1.1).  Specializing  $ \, a_\ell = 1 \, $  and  $ \,
a_r = 0 \, $  for all  $ \, r \not= \ell \, $  we get  $ \, c_{t+\ell}
= Q^0_{t+\ell}(b_*) + Q^\ell_t(b_*) = b_{t+\ell} + Q^\ell_t(b_*) \, $.
In particular setting  $ \, b_* = 1_* \, $  we have that  $ \, 1 +
Q^\ell_t(1_*) \, $  is the coefficient  $ \, c_{\ell+t} \, $  of
$ \, x^{\ell + t + 1} \, $  in the series
  $$  \displaylines{
   \big( \, x \, + x^{\ell+1} \big) \circ \Big( \, x \,
+ {\textstyle \sum_{m=1}^{+\infty}} \, x^{m+1} \Big) \;
= \; \big( \, x \, + x^{\ell+1} \big) \circ
\big( x \cdot {(1-x)}^{-1} \big) \; =   \hfill  \cr
   = \; x \cdot {(1-x)}^{-1} + {\big( x \cdot {(1-x)}^{-1} \big)}^{\ell+1}
\; = \; {\textstyle \sum_{m=0}^{+\infty}} \, x^{m+1} \, + \,
x^{\ell+1} \Big( {\textstyle \sum_{m=0}^{+\infty}} \, x^m \Big)^{\ell+1}
=  \cr
   {} \hfill   = \; {\textstyle \sum_{m=0}^{+\infty}} \, x^{m+1}
\, + \, x^{\ell+1} \, {\textstyle \sum_{n=0}^{+\infty} {{\ell+n}
\choose \ell}} \, x^n \; = \; {\textstyle \sum_{s=0}^{\ell-1}}
\, x^{s+1} \, + \, {\textstyle \sum_{s=\ell}^{+\infty}
\left( 1 + {s \choose \ell} \right)} \, x^{s+1}  \quad ;  \cr }  $$
therefore  $ \, 1 + Q^\ell_t(1_*) = c_{\ell+t} = 1 + {{\ell+t} \choose
\ell} \, $,  \, whence  $ \, Q^\ell_t(1_*) = {{\ell+t} \choose
\ell} \, $.  As an alternative approach, one can prove that
$ \, Q^\ell_t(1_*) = {{\ell+t} \choose \ell} \, $  by induction
using the recurrence formula  $ \; Q^\ell_t(\x_*) \, = \,
Q^{\ell-1}_t(\x_*) + \sum_{s=1}^{t-1} Q^{\ell-1}_{t-s}(\x_*)
\, \x_s + \x_t \; $  and the identity  $ \; {{\ell+t} \choose
\ell} = \sum_{s=0}^t {{\ell+t-1} \choose {\ell-1}} \; $.
                                               \par
   The outcome is  $ \, N = Q^\ell_t(1_*) = {{\ell+t} \choose \ell}
\, $  (for all  $ t, \ell \, $),  thus  $ \; Q^\ell_t(\a_*) - {{\ell+t}
\choose \ell} \, \a_t \, = \, Q_- + Q_+ - {{\ell+t} \choose \ell} \,
\a_t \, = \, Q_- + N \, \a_t - {{\ell+t} \choose \ell} \, \a_t \, = \,
Q_- \, $.  Now, by definition  $ \, \tau(\boldalpha_{j_r}) = j_r - 1
\, $  and  $ \, \tau\big({\boldalpha_1}^{\!j_r} \big) = j_r \, $.
Therefore if  $ \, \boldalpha_{(j_r)} \in \big\{ \boldalpha_{j_r},
{\boldalpha_1}^{\!j_r} \big\} \, $  (for all  $ \, r = 1, \dots, s
\, $)  and  $ \, (\boldalpha_{(j_1)}, \boldalpha_{(j_2)}, \dots,
\boldalpha_{(j_s)}) \not= \big( {\boldalpha_1}^{\!j_1},
{\boldalpha_1}^{\!j_2}, \dots, {\boldalpha_1}^{\!j_s} \big) \, $,
\, then  $ \, \tau \big( \boldalpha_{(j_1)} \cdots \boldalpha_{(j_s)}
\big) \leq j_1 + \cdots + j_s - 1 = t - 1 \, $.  Then by construction
$ \, \tau(Q_-) \leq t - 1 \, $,  \, whence, since  $ \,
Z^\ell_t(\boldalpha_*) := Q^\ell_t(\a_*) - {{\ell+t} \choose \ell}
\, \a_t = Q_- \, $,  \, we get also  $ \, \tau \big( Z^\ell_t
(\boldalpha_*) \big) \leq t - 1 \, $,  \, i.e.~$ \, Z^\ell_t
(\boldalpha_*) \in \varTheta_{t-1} \, $,  \, so  $ \;
Q^\ell_t(\a_*) \, = \, Z^\ell_t(\boldalpha_*) \, + \,
{{\ell+t} \choose t} \, {\boldalpha_1}^{\!t} \in
\varTheta_t \setminus \varTheta_{t-1} \, $.   \qed
\enddemo

\vskip7pt

\proclaim{Proposition 4.1}  $ \, \underline{\varTheta} $  is
a Hopf algebra filtration of  $ \, \calH \, $.
\endproclaim

\demo{Proof} \, By construction (cf.~\S 4.3)
$ \underline{\varTheta} $  is an algebra filtration; so to
check it is  {\sl Hopf\/}  too we are left only to show that
$ \; (\star) \, \Delta(\varTheta_n) \subseteq \sum_{r+s=n}
\varTheta_r \otimes \varTheta_s \; $  (for all  $ \, n \in \N \, $),
\, for then  $ \, S(\varTheta_n) \subseteq \varTheta_n \, $  (for
all  $ \, n \, $)  will follow from that by recurrence (and Hopf
algebra axioms).
                                                  \par
   By definition  $ \, \varTheta_0 = \Bbbk \cdot 1_{\scriptscriptstyle
\calH} \, $;  then  $ \, \Delta(1_{\scriptscriptstyle \calH}) =
1_{\scriptscriptstyle \calH} \otimes 1_{\scriptscriptstyle \calH} \, $
proves  $ (\star) $  for  $ \, n = 0 \, $.  For  $ \, n = 1 \, $,  \,
by definition  $ \varTheta_1 $  is the direct sum of  $ \varTheta_0 $
with the (free) Lie (sub)algebra (of  $ \calH \, $)  generated by
$ \{\boldalpha_1,\boldalpha_2\} $.  Since  $ \, \Delta(\boldalpha_1)
= \boldalpha_1 \otimes 1 + 1 \otimes \boldalpha_1 \, $  and  $ \,
\Delta(\boldalpha_2) = \boldalpha_2 \otimes 1 + 1 \otimes \boldalpha_2
\, $  and
  $$  \Delta\big([x,y]\big) \, = \, \big[\Delta(x),\Delta(y)\big]
\, = \, {\textstyle \sum_{(x),(y)}} \big( [x_{(1)},y_{(1)}] \otimes
x_{(2)} y_{(2)} + x_{(1)} y_{(1)} \otimes [x_{(2)},y_{(2)}] \big)  $$
(for all  $ \, x, y \in \calH \, $)  we argue  $ (\star) $  for
$ \, n = 1 \, $  too.  Moreover, for every  $ \, n > 1 \, $
(setting $ \, Q^n_0(\a_*) = 1 = \a_0 \, $  for short) we have
$ \; \Delta(\boldalpha_n) = \Delta(\a_n) - \Delta \big(
{\a_1}^{\!n} \big) = {\textstyle \sum_{k=0}^n} \, \a_k
\otimes Q^k_{n-k}(\a_*) - {\textstyle \sum_{k=0}^n \,
{n \choose k}} \, {\a_1}^k \otimes \, {\a_1}^{\!n-k} \;
= {\textstyle \sum_{k=2}^n} \, \boldalpha_k \otimes \,
Q^k_{n-k}(\a_*) + {\textstyle \sum_{k=0}^{n-1}} \,
{\boldalpha_1}^{\!k} \otimes \, Z^k_{n-k}(\boldalpha_*) \, $,
\; and therefore  $ \, \Delta(\boldalpha_n) \in \sum_{r+s=n-1}
\varTheta_r \otimes \varTheta_s \, $  thanks to Lemma 4.1
(and to  $ \, \boldalpha_m \in \varTheta_{m-1} \, $  for
$ \, m > 1 \, $).
                                                   \par
   Finally, as  $ \, \Delta\big([x,y]\big) \! = \! \big[ \Delta(x),
\Delta(y) \big] \! = \! \sum_{(x),(y)} \! \big( [x_{(1)},y_{(1)}]
\otimes x_{(2)} y_{(2)} + x_{(1)} y_{(1)} \otimes [x_{(2)},y_{(2)}]
\big) $  and similarly  $ \, \Delta(x\,y) = \Delta(x) \Delta(y) =
\sum_{(x),(y)} x_{(1)} y_{(1)} \otimes x_{(2)} y_{(2)} \, $  (for
$ \, x, y \in \calH \, $),  we have that  $ \Delta $  does not
increase  $ \, (\partial_- - d\,) \, $:  \, as  $ \varTheta $  is
exactly the (algebra) filtration induced by  $ (\partial_- - d\,)
\, $,  \, it is a Hopf algebra filtration as well.   \qed
\enddemo

\vskip7pt

\proclaim{Lemma 4.2}  (notation of \S 4.3)
                                    \hfill\break
   \indent   (a)  $ \; \kappa\,(a) \leq \partial(a) \; $  for
every  $ \, a \in \calH \setminus \{0\} \, $  which is
$ \, \partial(a) $--homogeneous.
                                    \hfill\break
   \indent   (b)  $ \; \kappa\,(a\,a') \leq \kappa\,(a) + \kappa\,(a')
\; $  and  $ \; \kappa\,\big([a,a']\big) < \kappa\,(a) + \kappa\,(a')
\; $  for all  $ \, a, a' \in \calH \setminus \{0\} \, $.
                                    \hfill\break
   \indent   (c)  $ \; \kappa\,(\boldalpha_n) = \partial_-(\boldalpha_n)
= \tau(\boldalpha_n) \; $  for all  $ \, n \in \N_\nu \, $.
                                    \hfill\break
   \indent   (d)  $ \; \kappa \, \big( [\boldalpha_r,\boldalpha_s] \big)
= \partial_-(\boldalpha_r) + \partial_-(\boldalpha_s) - 1 = \tau \big(
[\boldalpha_r, \boldalpha_s] \big) \; $  for all  $ \, r $,  $ s \in
\N_\nu \, $  with  $ \, r \not= s \, $.
                                    \hfill\break
   \indent   (e)  $ \; \kappa\,(\boldalpha_b) = \partial_-(\boldalpha_b)
- d(\boldalpha_b) + 1 = \tau(\boldalpha_b) \; $  for every  $ \, b
%
%
 \in B_\nu \, $.
                                    \hfill\break
   \indent   (f)  $ \; \kappa\,(\boldalpha_{b_1} \boldalpha_{b_2}
\cdots \boldalpha_{b_\ell}) = \tau(\boldalpha_{b_1} \boldalpha_{b_2}
\cdots \boldalpha_{b_\ell}) \; $  for all  $ \, b_1, b_2, \dots,
b_\ell \in B_\nu \, $.
                                       \hfill\break
   \indent   (g)  $ \; \kappa \big([\boldalpha_{b_1},
\boldalpha_{b_2}]\big) = \kappa\,(\boldalpha_{b_1})
+ \kappa\,(\boldalpha_{b_2}) - 1 = \tau \big([\boldalpha_{b_1},
\boldalpha_{b_2}]\big) \, $,  \; for all  $ \, b_1, b_2 \in B_\nu \, $.
\endproclaim

\demo{Proof} {\it (a)} \, Let  $ \, a \in \calH \setminus \{0\} \, $
be  $ \partial(a) $--homogeneous.  Since  $ \calH $  is graded, we
have  $ \, \partial \big( \delta_\ell(a) \big) = \partial(a) \, $  for
all  $ \ell \, $;  \, moreover,  $ \, \delta_\ell(a) \in \! J^{\otimes
\ell} \, $  (with  $ \, J := \hbox{\sl Ker}\,(\epsilon_{\calH}) \, $)
by definition, and  $ \, \partial(y) > 0 \, $  for each
$ \partial $--homogeneous  $ \, y \in J \setminus \{0\} \, $.
Then  $ \, \delta_\ell(a) = 0 \, $  for all  $ \, \ell >
\partial(a) \, $,  \, whence the claim.
                                             \par
   {\it (b)} \, Let  $ \, a \in D_m \, $,  $ \, b \in D_n \, $:  \,
then  $ \, a \, b \in D_{m+n} \, $  by property  {\it (c)\/}  in \S
3.2.  Similarly, we have  $ \, [a,b] \in D_{m+n-1} \leq m+n-1 \, $ 
because of property  {\it (d)\/}  in \S 3.2.  The claim follows.
                                                 \par
   {\it (c)} \, By part  {\it (a)\/}  we have  $ \, \kappa(\a_n)
\leq \partial(\a_n) = n \, $.  Moreover, by definition  $ \,
\delta_2(\a_n) = \sum_{k=1}^{n-1} \a_k \otimes Q^k_{n-k}(\a_*)
\, $,  thus  $ \, \delta_n(\a_n) = (\delta_{n-1} \otimes \delta_1)
\big( \delta_2(\a_n) \big) = \sum_{k=1}^{n-1} \delta_{n-1}(\a_k)
\otimes \delta_1 \big( Q^k_{n-k} (\a_*) \big) \, $  by coassociativity.
Since  $ \, \delta_\ell(\a_m) = 0 \, $  for  $ \, \ell > m \, $,  $ \,
Q^{n-1}_1(\a_*) = n \, \a_1 \, $  and  $ \, \delta_1(\a_1) = \a_1 \, $,
\, we have  $ \, \delta_n(\a_n) = \delta_{n-1}(\a_{n-1}) \otimes (n
\, \a_1) \, $,  \, thus by induction  $ \; \delta_n(\a_n) = n! \,
{\a_1}^{\!\otimes n} \, $  ($ \, \not= 0 \, $),  \, whence  $ \,
\kappa(\a_n) = n \, $.  But also  $ \, \delta_n({\a_1}^{\!n}) = n!
\, {\a_1}^{\!\otimes n} \, $.  Thus  $ \, \delta_n(\boldalpha_n) =
\delta_n(\a_n) - \delta_n({\a_1}^{\!n}) = 0 \, $  for  $ \, n > 1 \, $.
                                                  \par
   Clearly  $ \, \kappa(\boldalpha_1) = 1 \, $.  For the general case,
for all  $ \, \ell \geq 2 \, $  we have
  $$ \, \delta_{\ell-1}(\a_\ell) \; = \; (\delta_{\ell-2} \otimes
\delta_1) \big( \delta_2(\a_\ell) \big) \; = \; {\textstyle
\sum\nolimits_{k=1}^{\ell-1}} \; \delta_{\ell-2}(\a_k)
\otimes \delta_1 \big( Q^k_{\ell-1-k}(\a_*) \big)  $$
which by the previous analysis gives
$ \; \delta_{\ell-1}(\a_\ell) \, = \, \delta_{\ell-2}(\a_{\ell-2})
\otimes \big( {(\ell-1)} \, \a_2 + {{\ell-1} \choose 2} \, {\a_1}^{\!2}
\big) + \delta_{\ell-2} (\a_{\ell-1}) \otimes \ell \, \a_1 \, = \,
{(\ell-1)}! \cdot {\a_1}^{\!\otimes (\ell-2)} \otimes \left( \a_2
+ {{\;\ell-1\,} \over {\,2\,}} \cdot {\a_1}^{\!2} \right) \, +
\, \ell \cdot \delta_{\ell-2} (\a_{\ell-1}) \otimes \a_1 \; $.
%
%
Iterating we get, for all  $ \, \ell \geq 2 \, $  (with
$ \; {{-1} \choose 2} := 0 \, $,  \, and changing indices)
  $$  \delta_{\ell-1}(\a_\ell) \; = \; {\textstyle
\sum\nolimits_{m=1}^{\ell-1} {{\,\ell\,!\,} \over {m+1}}}
\cdot {\a_1}^{\!\otimes (m-1)} \otimes \left( \a_2 + {\textstyle
{{\,m-1\,} \over {\,2\,}}} \cdot {\a_1}^{\!2} \right)
\otimes {\a_1}^{\!\otimes (\ell-1-m)} \; .  $$
   \indent   On the other hand, we have also  $ \; \delta_{\ell-1}
\big({\a_1}^{\!\ell}\big) \, = \, \sum\nolimits_{m=1}^{\ell-1}
{{\,\ell\,!\,} \over {\,2\,}} \cdot {\a_1}^{\!\otimes (m-1)}
\otimes {\a_1}^{\!2} \otimes {\a_1}^{\! \otimes (\ell-1-m)} \, $.
Therefore, for  $ \, \delta_{n-1}(\boldalpha_n) = \delta_{n-1}(\a_n)
- \delta_{n-1}({\a_1}^{\!n}) \, $  (for all  $ \, n \in \N_\nu \, $,
$ \, n \geq 2 \, $) the outcome is
  $$  \hbox{ $ \eqalign{
   \delta_{n-1}(\boldalpha_n) \;  &  = \; {\textstyle \sum_{m=1}^{n-1}
{{\,n!\,} \over {\,m+1\,}}} \cdot {\a_1}^{\!\otimes (m-1)} \otimes
\big( \a_2 - {\a_1}^{\!2} \big) \otimes {\a_1}^{\!\otimes (n-1-m)}
\; =  \cr
   {}  &  = \; {\textstyle \sum_{m=1}^{n-1} {{\,n!\,} \over {\,m+1\,}}}
\cdot {\boldalpha_1}^{\!\otimes (m-1)} \otimes \boldalpha_2 \otimes
{\boldalpha_1}^{\!\otimes (n-1-m)}  \quad ;  \cr } $ }   \eqno (4.2)  $$
in particular  $ \, \delta_{n-1}(\boldalpha_n) \not= 0 \, $,  \,
whence  $ \, \boldalpha_n \not\in D_{n-2} \, $  and so
$ \, \kappa(\boldalpha_n) = n-1 \, $,  \, q.e.d.
                                                  \par
   {\it (d)} \, Let  $ \, r \not= 1 \not= s \, $.  From
{\it (b)--(c)\/} we get  $ \, \kappa \big([\boldalpha_r,
\boldalpha_s]\big) < \kappa(\boldalpha_r) + \kappa(\boldalpha_s)
= r+s-2 \, $.  In addition, we prove that  $ \; \delta_{r+s-3}
\big([\boldalpha_r,\boldalpha_s]\big) \not= 0 \, $,  \, yielding
{\it (d)}.  Property  {\it(d)\/}  in \S 3.2  gives
  $$  \displaylines{
   \delta_{r+s-3}\big([\boldalpha_r,\boldalpha_s]\big) \; =
\hskip-9pt  \sum_{\Sb \Lambda \cup Y = \{1,\dots,r+s-3\}  \\
\Lambda \cap Y \not= \emptyset  \endSb}  \hskip-26pt
\big[ \delta_\Lambda(\boldalpha_r), \delta_Y(\boldalpha_s) \big]
\; =  \hskip-21pt
\sum_{\Sb \Lambda \cup Y = \{1,\dots,r+s-3\} \\
\Lambda \cap Y \not= \emptyset, \, |\Lambda| = r-1, \,
|Y| = s-1 \endSb}  \hskip-38pt  \big[ j_\Lambda \big( \delta_{r-1}
(\boldalpha_r) \big), j_Y \big( \delta_{s-1}(\boldalpha_s) \big)
\big] \, .  \cr }  $$
Using (4.2) in the form  $ \; \delta_{\ell-1}(\a_\ell) =
\sum_{m=1}^{\ell-1} {{\,\ell\,!\,} \over 2} \cdot \boldalpha_2
\otimes {\boldalpha_1}^{\!\otimes (\ell-2)} + \boldalpha_1
\otimes \eta_\ell \; $  (for some  $ \, \eta_\ell \in \calH \, $),
\, and counting how many  $ \Lambda $'s  and  $ Y $'s  exist with
$ \, 1 \in \Lambda \, $  and  $ \, \{1,2\} \subseteq Y \, $,  \,
and   --- conversely ---   how many of them exist with  $ \,
\{1,2\} \subseteq \Lambda \, $  and  $ \, 1 \in Y \, $,  \,
we argue
  $$  \delta_{r+s-3}\big([\boldalpha_r,\boldalpha_s]\big) \; =
\; c_{r,s} \cdot [\boldalpha_2,\boldalpha_1] \otimes \boldalpha_2
\otimes {\boldalpha_1}^{\! \otimes (r+s-5)} \, + \, \boldalpha_1
\otimes \varphi_1 \, + \, \boldalpha_2 \otimes \varphi_2 \, + \,
[\boldalpha_2,\boldalpha_1] \otimes \boldalpha_1 \otimes \psi  $$
for some  $ \, \varphi_1, \varphi_2 \in \calH^{\otimes (r+s-4)}
\, $,  $ \, \psi \in \calH^{\otimes (r+s-5)} \, $,  \, and with
  $$  {\textstyle c_{r,s} = {r! \over 2} \cdot {s! \over 3} \cdot
{{r+s-5} \choose {r-2}} - {s! \over 2} \cdot {r! \over 3} \cdot
{{s+r-5} \choose {s-2}} = {{\,2\,} \over {\,3\,}} \, {r \choose 2}
{s \choose 2} (s-r) (r+s-5)! \not= 0 \, .}  $$
In particular  $ \; \delta_{r+s-3}\big([\boldalpha_r,\boldalpha_s]\big)
\, = \, c_{r,s} \cdot [\boldalpha_2,\boldalpha_1] \otimes \boldalpha_2
\otimes {\boldalpha_1}^{\! \otimes (r+s-5)} \, + \, \text{\sl l.i.t.}
\, $,  \; where  ``{\sl l.i.t.}''  stands for some further  {\sl
terms\/}  which are  {\sl linearly independent\/}  of  $ \,
[\boldalpha_2, \boldalpha_1] \otimes \boldalpha_2 \otimes
{\boldalpha_1}^{\! \otimes (r+s-5)} \, $  and  $ \, c_{r,s}
\not= 0 \, $.  Then  $ \; \delta_{r+s-3} \big( [\boldalpha_r,
\boldalpha_s] \big) \not= 0 \, $,  \; q.e.d.
                                              \par
   Finally, if  $ \, r > 1 = s \, $  (and similarly if  $ \, r
= 1 < s \, $)  things are simpler.  Indeed, again  {\it (b)\/}
and  {\it (c)\/}  together give  $ \, \kappa \big([\boldalpha_r,
\boldalpha_1]\big) < \kappa(\boldalpha_r) + \kappa(\boldalpha_1)
= (r-1) + 1 = r \, $,  \, and we prove that  $ \; \delta_{r-1}
\big( [\boldalpha_r, \boldalpha_1] \big) \not= 0 \, $.  Like
before,  property  {\it(d)\/}  in \S 3.2 gives (since  $ \,
\delta_1(\boldalpha_1) = \boldalpha_1 \, $)
  $$  \displaylines{
   \delta_{r-1}\big([\boldalpha_r,\boldalpha_1]\big) \; =
\hskip-15pt  \sum_{\Sb \Lambda \cup Y = \{1,2,\dots,r-1\} \\
\Lambda \cap Y \not= \emptyset, \, |\Lambda| = r-1, \, |Y| = 1
\endSb}  \hskip-35pt  \big[ \delta_\Lambda(\boldalpha_r),
\delta_Y(\boldalpha_1) \big]  =  {\textstyle \sum\limits_{k=1}^{r-1}}
\, \Big[ \delta_{r-1}(\boldalpha_r), 1^{\otimes (k-1)} \otimes
\boldalpha_1 \otimes 1^{\otimes (r-1-k)} \Big] =  \cr
   {} \hfill   = \; {\textstyle \sum\limits_{m=1}^{r-1}} \, {{\,r!\,}
\over {\,m+1\,}} \cdot {\boldalpha_1}^{\!\otimes (m-1)} \otimes
[\boldalpha_2, \boldalpha_1] \otimes {\boldalpha_1}^{\!\otimes
(n-1-m)} \; \not= \; 0  \cr }  $$
                                              \par
   {\it (e)} \, We perform induction upon  $ d(b) \, $:  \, the case
$ \, d(b) < 2 \, $  is dealt with in parts  {\it (c)\/}  and  {\it
(d)\/},  thus we assume  $ \, d(b) \geq 2 \, $,  \, so  $ \, b =
\big[ b', x_\ell \big] \, $  for some  $ \, \ell \in \N_\nu \, $
and some  $ \, b' \in B_\nu \, $  with  $ \, d(b') = d(b) - 1 \, $;
\, then  $ \, \tau(\boldalpha_b) = \tau\big([\boldalpha_{b'},
\boldalpha_\ell]\big) = \tau(\boldalpha_{b'}) + \tau(\boldalpha_\ell)
- 1 \, $,  \, directly from definitions.  Moreover  $ \,
\tau(\boldalpha_\ell) = \kappa\,(\boldalpha_\ell) \, $
by part  {\it (c)},  \, and  $ \, \tau(\boldalpha_{b'}) =
\kappa\,(\boldalpha_{b'}) \, $  by inductive assumption.
                                       \par
   From  {\it (b)\/}  we have  $ \, \kappa(\boldalpha_b) =
\kappa\big([\boldalpha_{b'},\boldalpha_\ell]\big) \leq
\kappa(\boldalpha_{b'}) + \kappa(\boldalpha_\ell) - 1  =
\tau(\boldalpha_{b'}) + \tau(\boldalpha_\ell) - 1 = \tau(\boldalpha_b)
\, $,  \, i.{} e.{}  $ \, \kappa(\boldalpha_b) \leq \tau(\boldalpha_b)
\, $;  \, we must prove the converse, for which it is enough to show
  $$  \delta_{\tau(\boldalpha_b)}(\boldalpha_b) \; = \; c_b \cdot
[\, \cdots [\, [ \undersetbrace{d(b)+1}\to{\boldalpha_1,\boldalpha_2],
\boldalpha_2], \dots, \boldalpha_2} \,] \otimes \boldalpha_2 \otimes
{\boldalpha_1}^{\!\otimes (\tau(\boldalpha_b)-2)} \; + \;
\text{\sl l.i.t.}   \eqno (4.3)  $$
for some  $ \, c_b \in \Bbbk \setminus \{0\} \, $,  \, where
``{\sl l.i.t.}''  means the same as before.
                                      \par
   Since  $ \, \tau(\boldalpha_b) = \tau \big( [\boldalpha_{b'},
\boldalpha_\ell] \big) = \tau(\boldalpha_{b'}) + \ell - 2 \, $,
\, using property  {\it (d)\/}  in \S 3.2 we get
  $$  \displaylines{
   {} \quad   \delta_{\tau(\boldalpha_b)}(\boldalpha_b) \; = \;
\delta_{\tau(\boldalpha_b)} \big([\boldalpha_{b'},\boldalpha_\ell]\big)
\; = \; {\textstyle \sum_{\Sb  \Lambda \cup Y = \{1, \dots,
\tau(\boldalpha_b)\}  \\
   \Lambda \cap Y \not= \emptyset  \endSb}}
\big[ \delta_\Lambda(\boldalpha_{b'}), \delta_Y(\boldalpha_\ell) \big]
\; = \;   \hfill  \cr
   {} \hfill   = \; {\textstyle \sum_{\Sb  \Lambda \cup Y = \{1, \dots,
\tau(\boldalpha_b)\} \, , \; \Lambda \cap Y \not= \emptyset  \\
   |\Lambda| = \tau(\boldalpha_{b'}) \, , \, |Y| = \ell-1  \endSb}}
\big[ j_\Lambda\big(\delta_{\tau(\boldalpha_{b'})}(\boldalpha_{b'})\big),
j_Y\big(\delta_{\ell-1}(\boldalpha_\ell)\big) \big] \; =  \cr
%
%
   {} \hfill   =  \hskip-79pt \sum_{ \Sb  {} \hskip71pt \Lambda \cup Y
= \{1, \dots, \tau(\boldalpha_b)\} \, , \; \Lambda \cap Y \not=
\emptyset  \\
   {} \hskip71pt |\Lambda| = \tau(\boldalpha_{b'}) \, ,
\quad  |Y| = \ell-1  \endSb}  \hskip-69pt
\Big[ j_\Lambda \big( c_{b'} \, [\, \cdots [
\undersetbrace{d(b')+1}\to{\boldalpha_1,\boldalpha_2],
\dots, \boldalpha_2} \,] \otimes \boldalpha_2 \otimes
{\boldalpha_1}^{\!\otimes (\tau(\boldalpha_{b'}\!)-2)} \big),
\, j_Y \big( \textstyle{{{\,\ell\,!\,} \over {\,2\,}}}
\, \boldalpha_2 \otimes {\boldalpha_1}^{\!\otimes (\ell-2)}
\big) \Big] + \, \text{\sl l.i.t.} =  \cr
   {} \hfill   = \; c_{b'} \cdot {\textstyle {{\,\ell\,!\,} \over {\,2\,}}
\cdot {{\tau(\boldalpha_b)-2\,} \choose {\ell-2}}} \cdot [\,[\, \cdots
[[ \undersetbrace{d(b')+1+1 \, = \, d(b)+1} \to{\boldalpha_2,
\boldalpha_1], \boldalpha_2], \dots, \boldalpha_2], \boldalpha_2}
\,\big] \otimes \boldalpha_2 \otimes {\boldalpha_1}^{\!\otimes
(\tau(\boldalpha_b)-2)} \; + \; \text{\sl l.i.t.}  \cr }  $$
(using induction about  $ \boldalpha_{b'} $);  this proves (4.3)
with  $ \, c_b =  c_{b'} \cdot {{\,\ell\,!\,} \over {\,2\,}}
\cdot \Big( \! {{\tau(\boldalpha_b)-2} \atop {\ell-2}} \Big)
\not= 0 \, $.
                                                 \par
   Thus (4.3) holds, yielding  $ \, \delta_{\tau(\boldalpha_b)}
(\boldalpha_b) \not= 0 \, $,  \, hence  $ \, \kappa(\boldalpha_b)
\geq \tau(\boldalpha_b) \, $,  \, q.e.d.
                                                 \par
   {\it (f)} \, The case  $ \, \ell = 1 \, $  is proved by part
{\it (e)},  so we can assume  $ \, \ell > 1 \, $.  By part  {\it
(b)\/}  and the case  $ \, \ell = 1 \, $  we have  $ \, \kappa\,
(\boldalpha_{b_1} \boldalpha_{b_2} \cdots \boldalpha_{b_\ell})
\leq \sum_{i=1}^\ell \kappa\,(\boldalpha_{b_i}) = \sum_{i=1}^\ell
\tau(\boldalpha_{b_i}) = \tau(\boldalpha_{b_1} \boldalpha_{b_2}
\cdots \boldalpha_{b_\ell}) \, $;  \, so we must only prove the
converse inequality.  We begin with  $ \, \ell = 2 \, $  and
$ \, d(b_1) = d(b_2) = 0 \, $,  \, so  $ \, \boldalpha_{b_1}
= \boldalpha_r \, $,  $ \, \boldalpha_{b_2} = \boldalpha_s
\, $,  \, for some  $ \, r $,  $ s \in \N_\nu \, $.
                                                 \par
   If  $ \, r = s = 1 \, $  then  $ \, \kappa(\boldalpha_r) =
\kappa(\boldalpha_s) = \kappa(\boldalpha_1) = 1 \, $,  \, by
part  {\it (c)}.  Then
  $$  \delta_2(\boldalpha_1 \, \boldalpha_1) \; =
\delta_2(\a_1 \, \a_1) \; = \; {(\id - \epsilon)}^{\otimes 2}
\Delta\big({\a_1}^{\!2}\big) \; = \; 2 \cdot \a_1 \otimes \a_1 \;
= \; 2 \cdot \boldalpha_1 \otimes \boldalpha_1 \; \not= \; 0  $$
so that  $ \, \kappa(\boldalpha_1 \, \boldalpha_1) \geq 2
= \kappa(\boldalpha_1) + \kappa(\boldalpha_1) \, $,  \, hence
$ \, \kappa(\boldalpha_1 \, \boldalpha_1) = \kappa(\boldalpha_1)
+ \kappa(\boldalpha_1) \, $,  \, q.e.d.
                                                \par
   If  $ \, r > 1 = s \, $  (and similarly if  $ \, r = 1 < s
\, $)  then  $ \, \kappa(\boldalpha_r) = r-1 \, $,  $ \,
\kappa(\boldalpha_s) = \kappa(\boldalpha_1) = 1 \, $,  \,
by part  {\it (c)}.  Then  property  {\it(d)\/}  in \S 3.2 gives
  $$  \displaylines{
   \delta_r(\boldalpha_r \, \boldalpha_1) \; = \; {\textstyle
\sum_{\Sb
   \Lambda \cup Y = \{1,\dots,r\}  \\
   |\Lambda| = r-1 \, , \; |Y| = 1  \endSb}}  \hskip1pt
\delta_\Lambda(\boldalpha_r) \; \delta_Y(\boldalpha_1) \; = \;
{\textstyle \sum_{m=1}^r \sum_{k<m} {{\,r!\,} \over {\,m+1\,}}}
\, \times   \hfill  \cr
   \hfill   \times \big( {\boldalpha_1}^{\!\otimes
(k-1)} \otimes 1 \otimes {\boldalpha_1}^{\!\otimes (m-1-k)} \otimes
\boldalpha_2 \otimes {\boldalpha_1}^{\!\otimes (r-1-m)} \big)
\times \big( 1^{\otimes (k-1)} \otimes \boldalpha_1 \otimes
1^{\otimes (r-k)} \big) \; +  \quad {}  \cr
   + {\textstyle \sum\limits_{m=1}^r \sum\limits_{k>m}
{{\,r!\,} \over {\,m+1\,}}} \, \big( {\boldalpha_1}^{\!\otimes
(m-1)} \otimes \boldalpha_2 \otimes {\boldalpha_1}^{\!\otimes
(k-1-m)} \otimes 1 \otimes {\boldalpha_1}^{\!\otimes (r-1-k)}
\big) \times \big( 1^{\otimes (k-1)} \otimes \boldalpha_1
\otimes 1^{\otimes (r-k)} \big) =  \cr
%
%
%
%
   {} \hfill   = \; {\textstyle \sum\nolimits_{m=1}^r {{\,r!\,} \over
{\,m+1\,}}} \cdot {\boldalpha_1}^{\!\otimes (m-1)} \otimes \boldalpha_2
\otimes {\boldalpha_1}^{\!\otimes (r-1-m)} \; \not= \; 0  \cr }  $$
so that  $ \, \kappa(\boldalpha_r \, \boldalpha_1) \geq r
= \kappa(\boldalpha_r) + \kappa(\boldalpha_1) \, $,  \, hence
$ \, \kappa(\boldalpha_r \, \boldalpha_1) = \kappa(\boldalpha_r)
+ \kappa(\boldalpha_1) \, $,  \, q.e.d.
                                                \par
   Finally let  $ \, r, s > 1 \, $  (and  $ \, r \not= s \, $).
Then  $ \, \kappa(\boldalpha_r) = r-1 \, $,  $ \,
\kappa(\boldalpha_s) = s-1 \, $,  \, by part  {\it (c)\/};
then property  {\it(d)\/}  in \S 3.2 gives
  $$  \delta_{r+s-2}\big(\boldalpha_r \, \boldalpha_s \big) \;
=  \hskip-7pt  \sum_{\Sb \Lambda \cup Y = \{1,\dots,r+s-2\}  \\
|\Lambda| = r-1 \, , \; |Y| = s-1  \endSb}  \hskip-21pt
\delta_\Lambda(\boldalpha_r) \cdot \delta_Y(\boldalpha_s) \;
=  \hskip-13pt  \sum_{\Sb \Lambda \cup Y = \{1,\dots,r+s-2\} \\
|\Lambda| = r-1 \, , \; |Y| = s-1  \endSb}  \hskip-21pt
j_\Lambda\big( \delta_{r-1}(\boldalpha_r) \big) \cdot
j_Y \big( \delta_{s-1}(\boldalpha_s) \big) \, .  $$
Using (4.2) in the form  $ \; \delta_{t-1}(\a_t) =
\sum_{m=1}^{t-1} {{\,t\,!\,} \over 2} \cdot \boldalpha_2
\otimes {\boldalpha_1}^{\!\otimes (t-2)} + \boldalpha_1
\otimes \eta_t \; $  (for some  $ \, \eta_t \in \calH \, $
and  $ \, t \in \{r,s\} \, $)  and counting how many
$ \Lambda $'s  and  $ Y $'s  exist with  $ \, 1 \in \Lambda
\, $  and  $ \, 2 \in Y \, $  and viceversa   --- actually, it
is a matter of counting  $ (r-2,s-2) $-shuffles  ---   we argue
  $$  \delta_{r+s-2}\big(\boldalpha_r \, \boldalpha_s\big) \;
= \; e_{r,s} \cdot \boldalpha_2 \otimes \boldalpha_2 \otimes
{\boldalpha_1}^{\! \otimes (r+s-4)} \, + \, \boldalpha_1
\otimes \varphi  $$
for some  $ \, \varphi \in \calH^{\otimes (r+s-3)} \, $  with
$ \; e_{r,s} = {\textstyle {{\,r!\,} \over {\,2\,}} \cdot
{{\,s!\,} \over {\,2\,}} \cdot \Big( {{r+s-4} \choose {r-2}} +
{{s+r-4} \choose {s-2}} \Big)} = {\textstyle {{\, r! \, s! \,}
\over {\,2\,}} \cdot {{r+s-4} \choose {r-2}}} \not= 0 \, $.
In particular  $ \; \delta_{r+s-2} \big( \boldalpha_r \,
\boldalpha_s \big) \, = \, e_{r,s} \cdot \boldalpha_2 \otimes
\boldalpha_2 \otimes {\boldalpha_1}^{\! \otimes (r+s-4)} \, +
\, \text{\sl l.i.t.} \, $,  \; where  ``{\sl l.i.t.}''  stands
again for some further  {\sl terms\/}  which are  {\sl linearly
independent\/}  of  $ \, \boldalpha_2 \otimes \boldalpha_2 \otimes
{\boldalpha_1}^{\! \otimes (r+s-4)} \, $  and  $ \, e_{r,s}
\not= 0 \, $.  Then  $ \; \delta_{r+s-2} \big( \boldalpha_r \,
\boldalpha_s \big) \not= 0 \, $,  \, so  $ \, \kappa(\boldalpha_r
\, \boldalpha_1) \geq r+s-2 = \kappa(\boldalpha_r) +
\kappa(\boldalpha_1) \, $,  \, q.e.d.
                                                \par
   Now let again  $ \, \ell = 2 \, $  but  $ \, d(b_1), d(b_2) >
0 \, $.  Set  $ \, \kappa_i := \kappa(\boldalpha_{b_i}) \, $  for
$ \, i = 1, 2 \, $.  Applying (4.3) to  $ \, b = b_1 \, $  and
$ \, b = b_2 \, $  (and reminding  $ \, \tau \equiv \kappa \, $)
gives
  $$  \displaylines{
   \delta_{\kappa_1 + \kappa_2}(\boldalpha_{b_1} \, \boldalpha_{b_2})
\; =  \hskip-4pt  {\textstyle \sum\limits_{\Lambda \cup Y = \{1, \dots,
\kappa_1 + \kappa_2\}}}  \hskip-23pt
\delta_\Lambda(\boldalpha_{b_1}) \, \delta_Y(\boldalpha_{b_2})
\; =  \hskip-4pt  {\textstyle \sum\limits_{\Sb
      \Lambda \cup Y = \{1, \dots, \kappa_1 + \kappa_2\}  \\
      |\Lambda| = \kappa_1, \, |Y| = \kappa_2  \endSb}}  \hskip-23pt
j_\Lambda\big(\delta_{\kappa_1}(\boldalpha_{b_1})\big) \,
j_Y\big(\delta_{\kappa_2}(\boldalpha_{b_2})\big) \; =  \cr
   = {\textstyle \sum\limits_{\Sb  \Lambda \cup Y = \{1, \dots,
\kappa_1 + \kappa_2\}  \\
       |\Lambda| = \kappa_1, \, |Y| = \kappa_2  \endSb}}  \hskip-9pt
j_\Lambda \big(\, c_{b_1} \cdot [\, \cdots [\, [ \undersetbrace{d(b_1)+1}
\to{\boldalpha_1,\boldalpha_2],\boldalpha_2], \dots, \boldalpha_2} \,]
\otimes \boldalpha_2 \otimes {\boldalpha_1}^{\!\otimes (\kappa_1 - 2)}
\; + \; \hbox{\sl l.i.t.} \,\big) \, \times   \hfill  \cr
   {} \hfill   \times \, j_Y \big(\, c_{b_2} \cdot [\, \cdots [\, [
\undersetbrace{d(b_2)+1} \to{\boldalpha_1,\boldalpha_2],\boldalpha_2],
\dots, \boldalpha_2} \,] \otimes \boldalpha_2 \otimes {\boldalpha_1}^{\!
\otimes (\kappa_2 - 2)} \; + \; \hbox{\sl l.i.t.} \,\big) \; =  \cr
   = 2 \, c_{b_1} c_{b_2} \! {\textstyle {{\kappa_1 +
\kappa_2 - 4 \,} \choose {\kappa_1 - 2}}} \cdot
[\, \cdots [ \undersetbrace{d(b_1)+1} \to
{\boldalpha_1, \! \boldalpha_2], \dots, \! \boldalpha_2} \,]
\otimes [\, \cdots [ \undersetbrace{d(b_2)+1} \to {\boldalpha_1,
\! \boldalpha_2], \dots, \! \boldalpha_2} \,] \otimes
{\boldalpha_2}^{\! \otimes 2} \otimes {\boldalpha_1}^{\!\!
\otimes (\kappa_1 + \kappa_2 - 4)} + \, \hbox{\sl l.i.t.}  \cr
%
%
 }  $$
which proves the claim for  $ \, \ell = 2 \, $.  In addition, we can
take this last result as the basis of induction (on  $ \ell \, $)
to prove the following: for all  $ \, \underline{b} := (b_1, \dots,
b_\ell) \in {B_\nu}^{\!\ell} \, $,  \, one has
  $$  \delta_{|\underline{\kappa}|}\left( {\textstyle
\prod_{i=1}^\ell} \boldalpha_{b_i} \!\right) \, = \, c_{\underline{b}}
\cdot \bigg( \!\! {\textstyle \bigotimes_{i=1}^\ell} \, [\, \cdots
[\, [\undersetbrace{d(b_i)+1} \to {\boldalpha_1,\boldalpha_2],
\boldalpha_2], \dots, \boldalpha_2} \,] \! \bigg) \otimes \,
{\boldalpha_2}^{\!\otimes \ell} \otimes \, {\boldalpha_1}^{\!
\otimes (\! |\underline{\kappa}| - 2\,\ell )} \, + \;
\hbox{\sl l.i.t.}   \hfill \hskip11pt (4.4)  $$
for some  $ \, c_{\underline{b}} \in \Bbbk \setminus \{0\} \, $,
\, with  $ \, |\underline{\kappa}| := \sum_{i=1}^\ell \kappa_i \, $
and  $ \, \kappa_i := \kappa(\boldalpha_{b_i}) \, $  ($ \, i = 1,
\dots, \ell \, $).  The induction step, from  $ \ell \, $  to
$ (\ell + 1) $,  amounts to compute  (with  $ \, \kappa_{\ell
+ 1} := \kappa (\boldalpha_{b_{\ell + 1}}) \, $)
  $$  \displaylines{
   \delta_{|\underline{\kappa}| + \kappa_{\ell + 1}} \big(
\boldalpha_{b_1} \cdots \boldalpha_{b_\ell} \cdot \boldalpha_{b_{\ell
+ 1}} \big) \;\; =  {\textstyle \sum\limits_{\Lambda \cup Y = \, \{1,
\dots,|\underline{\kappa}| + \kappa_{\ell + 1}\}}}  \hskip-11pt
\delta_\Lambda(\boldalpha_{b_1} \cdots \boldalpha_{b_\ell})
\, \delta_Y(\boldalpha_{b_{\ell + 1}}) \; =   \hfill  \cr
   {} \hfill   = {\textstyle \sum\limits_{\Sb
     \Lambda \cup Y = \, \{1, \dots, |\underline{\kappa}|
+ \kappa_{\ell + 1}\}  \\
     |\Lambda| = |\underline{\kappa}|, \,
|Y| = \kappa_{\ell + 1}  \endSb}}  \hskip-11pt
j_\Lambda\big(\delta_{|\underline{\kappa}|}
(\boldalpha_{b_1} \cdots \boldalpha_{b_\ell})\big) \cdot
j_Y \big( \delta_{\kappa_{\ell + 1}} (\boldalpha_{b_{\ell
+ 1}}) \big) \; =  \cr
 }  $$
  $$  \displaylines{
   =  \hskip-7pt  {\textstyle \sum\limits_{\Sb
       \Lambda \cup Y = \{1, \dots, |\underline{\kappa}|
+ \kappa_{\ell + 1}\}  \\
       |\Lambda| = |\underline{\kappa}|, \, |Y| =
\kappa_{\ell + 1} \endSb}}  \hskip-11pt
j_\Lambda \Big( c_{\underline{b}} \, \cdot \Big( {\textstyle
\bigotimes_{i=1}^\ell} \; [\, \cdots [\, [ \undersetbrace{d(b_i)+1}
\to {\boldalpha_1,\boldalpha_2], \boldalpha_2], \dots, \boldalpha_2}
\,] \Big) \otimes {\boldalpha_2}^{\!\otimes \ell} \otimes
{\boldalpha_1}^{\!\otimes (|\underline{\kappa}| - 2\,\ell)}
\; + \, \hbox{\sl l.{}i.{}t.}  \Big) \times   \hfill  \cr
   {} \hfill   \times j_Y \Big( c_{b_{\ell + 1}} \cdot [\, \cdots [\,
[ \undersetbrace{d(b_{\ell + 1})+1} \to{\boldalpha_1,\boldalpha_2],
\boldalpha_2], \dots, \boldalpha_2} \,] \otimes \boldalpha_2
\otimes {\boldalpha_1}^{\!\otimes (\kappa_{\ell + 1} - 2)}
\; + \; \hbox{\sl l.i.t.} \,\Big) \; =  \cr
   = \; c_{\underline{b}} \, c_{b_{\ell + 1}} \cdot (\ell + 1)
%
%
\, \Big(\! {\textstyle {{|\underline{\kappa}| + \kappa_{\ell + 1} -
2 \, (\ell + 1)} \atop {|\underline{\kappa}| - 2 \, \ell}}} \!\Big)
\cdot
\Big( {\textstyle \bigotimes_{i=1}^\ell} \; [\, \cdots [\, [
\undersetbrace{d(b_i)+1} \to {\boldalpha_1,\boldalpha_2],
\boldalpha_2], \dots, \boldalpha_2} \,] \Big) \otimes   \hfill  \cr
   {} \hfill   \otimes [\, \cdots [\, [ \undersetbrace{d(b_{\ell + 1})
+1} \to {\boldalpha_1, \boldalpha_2], \boldalpha_2], \dots,
\boldalpha_2} \,] \otimes {\boldalpha_2}^{\!\otimes (\ell + 1)} \otimes
{\boldalpha_1}^{\!\otimes (|\underline{\kappa}| + \kappa_{\ell + 1}
- 2 \, (\ell + 1))} \; + \, \hbox{\sl l.i.t.}  \cr }  $$
which proves (4.4) for  $ \, (\underline{b} \, , b_{\ell + 1}) \, $
with  $ \, c_{(\underline{b} \, , b_{\ell + 1})} = c_{\underline{b}}
\, c_{b_{\ell + 1}} \cdot (\ell + 1) \Big( {{|\underline{\kappa}| +
\kappa_{\ell + 1} - 2 \, (\ell + 1)} \atop {|\underline{\kappa}| -
2 \, \ell}} \Big) \, \not= \, 0 \, $.  Finally, (4.4) yields  $ \,
\delta_{|\underline{\kappa}|}(\boldalpha_{b_1} \hskip-0,7pt \cdots
\boldalpha_{b_\ell}) \not= 0 \, $,  \, so
       \hbox{$ \kappa(\boldalpha_{b_1} \cdots \boldalpha_{b_\ell})
\hskip-0,9pt \geq \hskip-1,5pt \kappa(\boldalpha_{b_1}) + \cdots
+ \hskip-0,3pt \kappa(\boldalpha_{b_\ell}) $,  q.e.d.}
                                         \par
   {\it (g)} \, Part  {\it (d)\/}  proves the claim for  $ \, d(b_1)
= d(b_2) = 0 \, $,  \, that is  $ \, b_1, b_2 \in {\{x_n\}}_{n \in \N}
\, $.  Moreover, when  $ \, b_2 = x_n \in {\{x_m\}}_{m \in \N_\nu}
\, $  we can replicate the proof of part  {\it (d)\/}  to show
that  $ \, \kappa \big( [\boldalpha_{b_1},\boldalpha_{b_2}]
\big) = \kappa \big( [\boldalpha_{b_1},\boldalpha_n] \big)
= \partial_-\big([\boldalpha_{b_1},\boldalpha_n]\big) -
d\big([\boldalpha_{b_1},\boldalpha_n]\big) \, $:  \, but
the latter is exactly  $ \, \tau \big([\boldalpha_{b_1},
\boldalpha_{b_2}]\big) \, $,  \, q.e.d.  Everything is
similar if  $ \, b_1 = x_n \in {\{x_m\}}_{m \in \N_\nu} \, $.
                                         \par
   Now let  $ \, b_1, b_2 \in B_\nu \setminus {\{x_n\}}_{n \in \N_\nu}
\, $.  Then  {\it (b)}  gives  $ \, \kappa\big([\boldalpha_{b_1},
\boldalpha_{b_2}]\big) \leq \kappa\,(\boldalpha_{b_1}) + \kappa\,
(\boldalpha_{b_2}) - 1 = \tau \big( [\boldalpha_{b_1},
\boldalpha_{b_2}] \big) $.  Applying (4.3) to  $ \, b =
b_1 \, $  and  $ \, b = b_2 \, $  we get, for  $ \, \kappa_i
:= \kappa(\boldalpha_{b_i}) $  ($ \, i = 1, 2 \, $)
  $$  \displaylines{
   \delta_{\kappa_1 + \kappa_2 - 1} \big( [\boldalpha_{b_1},
\boldalpha_{b_2}] \big) \; =  \hskip-4pt  {\textstyle \sum\limits_{\Sb
\Lambda \cup Y = \{ 1, \dots, \kappa_1 + \kappa_2 - 1 \}  \\
           \Lambda \cap Y \not= \emptyset  \endSb}}  \hskip-23pt
\big[ \delta_\Lambda(\boldalpha_{b_1}), \delta_Y(\boldalpha_{b_2})
\big] \; =   \hfill  \cr
%
%
   = {\textstyle \sum\limits_{\Sb  \Lambda \cup Y = \{1, \dots,
\kappa_1 + \kappa_2\}  \\
       |\Lambda| = \kappa_1, \, |Y| = \kappa_2  \endSb}}
\big[\, j_\Lambda \big(\, c_{b_1} \cdot [\, \cdots [\,
[ \undersetbrace{d(b_1)+1} \to{\boldalpha_2,\boldalpha_1],
\boldalpha_2], \dots, \boldalpha_2} \,] \otimes \boldalpha_2
\otimes {\boldalpha_1}^{\!\otimes (\kappa_1 - 2)} \; + \;
\hbox{\sl l.i.t.} \,\big) \, \times   \hfill  \cr
   {} \hfill   \times \, j_Y \big(\, c_{b_2} \cdot [\, \cdots [\, [
\undersetbrace{d(b_2)+1} \to{\boldalpha_2,\boldalpha_1],\boldalpha_2],
\dots, \boldalpha_2} \,] \otimes \boldalpha_2 \otimes {\boldalpha_1}^{\!
\otimes (\kappa_2 - 2)} \; + \, \hbox{\sl l.i.t.} \,\big) \big]
\; =  \cr
   = 2 \, c_{b_1}  c_{b_2} {\textstyle {{\kappa_1 + \kappa_2
- 4 \,} \choose {\kappa_1 - 2}}} \big[ [\, \cdots
[ \undersetbrace{d(b_1)+1} \to {\boldalpha_2, \! \boldalpha_1],
\dots, \! \boldalpha_2} \,] ,  [\, \cdots [ \undersetbrace{d(b_2)+1}
\to {\boldalpha_2, \! \boldalpha_1], \dots, \! \boldalpha_2} \,] \big]
\! \otimes {\boldalpha_2}^{\! \otimes 2} \otimes {\boldalpha_1}^{\!\!
\otimes (\kappa_1 + \kappa_2 - 4 \!)} + \, \hbox{\sl l.i.t.}  \cr
%
%
 }  $$
(note that  $ \, d(b_i) \geq 1 \, $  because  $ b_i \not\in \big\{\,
x_n \,\big|\, n \in \N_\nu \,\big\} \, $  for  $ \, i = 1, 2 \, $).
In particular this means  $ \, \delta_{\kappa_1 + \kappa_2 - 1}
\big( [\boldalpha_{b_1}, \boldalpha_{b_2}] \big) \not= 0 \, $,
\, thus  $ \, \kappa\big([\boldalpha_{b_1},\boldalpha_{b_2}]\big)
\geq \kappa\,(\boldalpha_{b_1}) + \kappa\,(\boldalpha_{b_2}) - 1
= \tau\big([\boldalpha_{b_1},\boldalpha_{b_2}]\big) \, $.   \qed
\enddemo

\vskip7pt

\proclaim{Lemma 4.3} \, Let  $ V $  be a  $ \Bbbk $--vector  space,
and  $ \, \psi \in \text{\sl Hom}_{\,\Bbbk}\big(V, V \wedge V) \, $.
Let  $ \L(V) $  be the free Lie algebra over  $ V $,  \, and  $ \,
\psi_{d\L} \in \text{\sl Hom}_{\,\Bbbk}\big(\L(V), \L(V) \wedge \L(V)
\big) \, $  the unique extension of  $ \, \psi $  from  $ V $  to
$ \L(V) $  by derivations, i.e.~such that  $ \, \psi_{d\L}{\big|}_V
= \psi \, $  and  $ \; \psi_{d\L}\big([x,y]\big) =
\big[ x \otimes 1 + 1 \otimes x,
 \allowbreak
\, \psi_{d\L}(y) \big] + \big[ \psi_{d\L}(x),
\, y \otimes 1 + 1 \otimes y \, \big]
= x.\psi_{d\L}(y) - y.\psi_{d\L}(x) \; $
in the  $ \L(V) $--module  $ \L(V) \wedge \L(V) \, $,
$ \, \forall \, x, y \in \L(V) \, $.
   \hbox{Let  $ \, K := \text{\sl Ker}\,(\psi) \, $:  then
$ \, \text{\sl Ker}\,\big(\psi_{d\L}\big) = \L(K) \, $,  \,
the free Lie algebra over  $ K \, $.}
\endproclaim

\demo{Proof} Standard, by universal arguments (for a direct
proof see [Ga2], Lemma 10.15).   \qed
\enddemo

\vskip7pt

\proclaim{Lemma 4.4}  The Lie cobracket  $ \delta $  of
$ \, U(\L_\nu) $  preserves  $ \tau $.  That is, for each  $ \,
\vartheta \in U(\L_\nu) \, $  in the expansion  $ \, \delta_2
(\vartheta) = \sum_{\underline{b}_1, \underline{b}_{\,2} \in \Bbb{B}}
c_{\underline{b}_1,\underline{b}_{\,2}} \, \boldalpha_{\underline{b}_1}
\otimes \boldalpha_{\underline{b}_{\,2}} \, $  (w.r.t.~the basis  $ \,
\Bbb{B} \otimes \Bbb{B} \, $,  \, where  $ \Bbb{B} $  is a PBW basis
as in \S 1.1 w.r.t.~some total order of  $ B_\nu $)  we have
$ \, \tau \big( \underline{\hat{b}}_{\,1} \big) + \tau \big(
\underline{\hat{b}}_{\,2} \big) = \tau(\vartheta) \, $  for some
$ \underline{\hat{b}}_{\,1} $,  $ \underline{\hat{b}}_{\,2} $  with
$ \, c_{\underline{\hat{b}}_1, \underline{\hat{b}}_{\,2}} \not= 0
\, $,  \, so
   \hbox{$ \, \tau \big( \delta(\vartheta) \big) := \max \big\{
\tau(\underline{b}_{\,1}) + \tau(\underline{b}_{\,2}) \;\big|\;
c_{\underline{b}_1,\underline{b}_{\,2}} \not= 0 \big\} =
\tau(\vartheta) \, $  if  $ \, \delta(\vartheta) \not= 0 \, $.}
\endproclaim

\demo{Proof}  It follows from Proposition 4.1 that  $ \, \tau \big(
\delta(\vartheta) \big) \leq \tau(\vartheta) \, $;  \, so  $ \; \delta
\, \colon \, U(\L_\nu) \longrightarrow {U(\L_\nu)}^{\otimes 2} \; $
is a morphism of filtered algebras, hence it naturally induces a
morphism of graded algebras  $ \; \overline{\delta} \, \colon \,
G_{\underline{\varTheta}}\big(U(\L_\nu)\big) \! \llongrightarrow
{G_{\underline{\varTheta}}\big(U(\L_\nu)\big)}^{\otimes 2} $.  Thus
proving the claim is equivalent to showing that  $ \, \text{\sl
Ker}\,\big(\overline{\delta}\,\big) = G_{\underline{\varTheta}
\, \cap {Ker}(\delta)} \big( \text{\sl Ker}\,(\delta) \big) =:
\overline{\text{\sl Ker}\,(\delta)} \, $,  \, the latter being
embedded into  $ G_{\underline{\varTheta}}\big(U(\L_\nu)\big) \, $.
                                       \par
   By construction,  $ \, \tau(x \, y - y \, x) = \tau\big([x,y]\big)
< \tau(x) + \tau(y) \, $  for  $ \, x $,  $ y \in U(\L_\nu) \, $,
\, so  $ G_{\underline{\varTheta}}\big(U(\L_\nu)\big) $  is
commutative: indeed, it is clearly isomorphic   --- as an algebra
---   to  $ S(V_\nu) $,  the symmetric algebra over  $ V_\nu \, $.
Moreover,  $ \delta $  acts as a derivation, that is  $ \; \delta(x\,y)
= \delta(x) \, \Delta(y) + \Delta(x) \, \delta(y) \; $  (for all
$ \, x $,  $ y \in U(\L_\nu) \, $),  \, thus the same holds for
$ \overline{\delta} $  too.  Like in Lemma 4.3, since  $ \,
G_{\underline{\varTheta}}\big(U(\L_\nu)\big) \, $  is generated
by  $ \; G_{\underline{\varTheta} \cap \L_\nu}(\L_\nu) =:
\overline{\L_\nu} \; $  it follows that  $ \, \text{\sl Ker}\,
\big(\overline{\delta}\,\big) \, $  is the free (associative
sub)algebra over  $ \, \text{\sl Ker}\, \Big( \overline{\delta}\,
{\big|}_{\overline{\L_\nu}} \Big) \, $,  \, in short  $ \; \text{\sl
Ker}\,\big(\overline{\delta}\,\big) = \Big\langle \text{\sl Ker}\,
\Big( \overline{\delta}\,{\big|}_{\overline{\L_\nu}} \Big) \Big\rangle
\; $.  Now, by definition  $ \; \delta(x_n) = \sum_{\ell=1}^{n-1}
(\ell+1) \, x_\ell \wedge x_{n-\ell} \; $  (cf.~Theorem 2.1) is
$ \tau $--{\,}homogeneous,  of  $ \tau $--{\,}degree  equal to
$ \, \tau(x_n) = n - 1 \, $.  As  $ \delta $  also enjoys  $ \,
\delta\big([x,y]\big) = \big[x \otimes 1 + 1 \otimes x, \delta(y)
\big] + \big[\delta(x), y \otimes 1 + 1 \otimes y \,\big] \, $
(for  $ \, x $,  $ y \in \L_\nu \, $)  we have that  $ \,
\delta{\big|}_{\L_\nu} \, $  is even  $ \tau $--{\,}homogeneous,
i.e.~such that  $ \, \tau\big(\delta(z)\big) = \tau(z) \, $,  \,
for any\break
$ \tau $--homogeneous
  $ \, z \in \L_\nu \, $  such that  $ \, \delta(z) \not= 0
\, $;  \, this implies that the induced map  $ \, \overline{\delta}
\,{\big|}_{\overline{\L_\nu}} \, $  enjoys  $ \; \overline{\delta}
\,{\big|}_{\overline{\L_\nu}} \big( \overline{\vartheta} \,\big)
= \overline{0} \iff \delta(\vartheta) = 0 \; $  for any  $ \,
\vartheta \in \L_\nu \, $,  \, whence  $ \, \text{\sl Ker}\,
\Big( \overline{\delta}\,{\big|}_{\overline{\L_\nu}} \Big) \,
= \, \overline{\text{\sl Ker}\, \big(\delta{\big|}_{\L_\nu}\big)}
\, $.  On the upshot we get  $ \; \text{\sl Ker}\,\big(
\overline{\delta}\,\big) = \Big\langle \text{\sl Ker}\, \Big(
\overline{\delta}\,{\big|}_{\overline{\L_\nu}} \Big) \Big\rangle
= \Big\langle \overline{\text{\sl Ker}\, \big(\delta{\big|}_{\L_\nu}
\big)} \, \Big\rangle = \overline{\text{\sl Ker}\,(\delta)} \; $,
\; q.e.d.   \qed  
\enddemo

\vskip7pt

\proclaim{Proposition 4.2}  $ \, \underline{D} = \underline{\varTheta}
\, $,  \, that is  $ \, D_n = \varTheta_n \, $  for all  $ \, n \in \N
\, $,  or  $ \, \kappa = \tau \, $.  Therefore, given any total order
$ \preceq $  in  $ B_\nu \, $,  \, the set  $ \, \calA_{\leq n} =
\calA \cap \varTheta_n = \calA \cap D_n \, $  of ordered monomials
  $$  \calA_{\leq n} \; = \; \Big\{\, \boldalpha_{\underline{b}} =
\boldalpha_{b_1} \cdots \boldalpha_{b_k} \,\Big|\; k \in \N \, ,
\; b_1, \dots, b_k \in B_\nu \, , \; b_1 \preceq \cdots \preceq
b_k \, ,  \; \tau(\underline{b}\,) \leq n \,\Big\}  $$
is a\/  $ \Bbbk $--basis  of  $ \, D_n \, $,  \, and  $ \, \calA_n
:= \big(\, \calA_{\leq n} \! \mod D_{n-1} \big) \, $  is a\/
$ \Bbbk $--basis  of  $ \, D_n \big/ D_{n-1} \; $  ($ \,
\forall \; n \in \N \, $).
\endproclaim

\demo{Proof} Both claims about the  $ \calA_{\leq n} $'s
and  $ \calA_n $'s  are equivalent to  $ \, \underline{D}
= \underline{\varTheta} \, $.  Also,  $ \, \calA_n :=
\big(\, \calA_{\leq n} \! \mod D_{n-1} \big) = \big(\,
\calA_{\leq n} \setminus \calA_{\leq n-1} \! \mod D_{n-1}
\big) \, $,  \, with  $ \, \calA_{\leq n} \setminus
\calA_{\leq n-1} = \big\{\, \boldalpha_{\underline{b}} \in
\calA \;\big|\; \tau(\underline{b}\,) = n \,\big\} \, $.
                                             \par
   By  Lemma 4.2{\it (f)\/}  we have  $ \, \calA_{\leq n} =
\calA \bigcap \varTheta_n \subseteq \calA \bigcap D_n \subseteq
D_n \, $;  \, since  $ \calA $  is a basis,  $ \calA_{\leq n} $  is
linearly independent and is a  $ \Bbbk $--basis  of  $ \varTheta_n $
(by definition): so  $ \, \varTheta_n \subseteq D_n \, $  for all
$ \, n \in \N \, $.
 \vskip3pt
   $ \underline{n=0} \, $:  \, By definition  $ \, D_0 := \hbox{\sl
Ker}(\delta_1) = \Bbbk \cdot 1_{\scriptscriptstyle \calH} =:
\varTheta_0 \, $,  \, spanned by  $ \,\calA_{\leq 0} =
\{1_{\scriptscriptstyle \calH}\} \, $,  \, q.e.d.
 \vskip3pt
   $ \underline{n=1} \, $:  \, Let  $ \, \eta' \in D_1 := \hbox{\sl
Ker}(\delta_2) \, $.  Let  $ \Bbb{B} $  be a PBW-basis of  $ \,
{\calH_\h}^{\!\vee} = U(\L_\nu) \, $  as in Lemma 4.4; expanding
$ \eta' $  w.r.t.~$ \calA $  we have  $ \, \eta' =
\sum_{\boldalpha_{\underline{b}} \in \calA} c_{\underline{b}}
\, \boldalpha_{\underline{b}} = \sum_{\underline{b} \in \Bbb{B}}
c_{\underline{b}} \, \boldalpha_{\underline{b}} \, $.  Then  $ \,
\eta := \eta' - \sum_{\tau(\underline{b}\,) \leq 1} c_{\underline{b}}
\, \boldalpha_{\underline{b}} = \sum_{\tau(\underline{b}\,) > 1}
c_{\underline{b}} \, \boldalpha_{\underline{b}} \in D_1 \, $,
           \hbox{\, since  $ \, \boldalpha_{\underline{b}} \in
\calA_1 \subseteq \varTheta_1 \subseteq D_1 \, $  for  $ \,
\tau(\underline{b}\,) \leq 1 \, $.}
                                              \par
   Now,  $ \, \boldalpha_1 := \a_1 \, $  and  $ \, \boldalpha_s :=
\a_s - {\a_1}^{\!s} = \h \, \big( \x_s + \h^{s-1} {\x_1}^{\!s} \big)
\, $  for all  $ \, s \in \N_\nu \setminus \{1\} \, $  yield
  $$  \eta \; = \,  {\textstyle \sum_{ \hskip-5pt \Sb  \underline{b}
\in \Bbb{B}  \\   \tau(\underline{b}\,) > 1  \endSb}} \hskip-3pt
c_{\underline{b}} \, \boldalpha_{\underline{b}}  \; = \;
{\textstyle \sum_{ \hskip-5pt \Sb  \underline{b}
\in \Bbb{B}  \\   \tau(\underline{b}\,) > 1  \endSb}} \hskip-3pt
\h^{\,g(\underline{b}\,)} \, c_{\underline{b}} \, \big(
\x_{\underline{b}} + \h \, \chi_{\underline{b}} \big)
\; \in \;  {\calH_\h}^{\!\vee}  $$
for some  $ \, \chi_{\underline{b}} \in {\calH_\h}^{\!\vee} \, $:
\, hereafter we set  $ \, g(\underline{b}\,) := k \, $  for each
$ \, \underline{b} = b_1 \cdots b_k \in \Bbb{B} \, $  (i.e.~$ \,
g(\underline{b}\,) \, $  is the degree of  $ \underline{b} $  as
a monomial in the  $ b_i $'s).  If  $ \, \eta \not= 0 \, $,  \,
let  $ \, g_0 := \min \big\{\, g(\underline{b}\,) \,\big|\,
\tau(\underline{b}\,) > 1 \, , \, c_{\underline{b}} \not= 0 \,\big\}
\, $;  \, then  $ \, g_0 > 0 \, $,  $ \, \eta_+ := \h^{-g_0} \,
\eta \in {\calH_\h}^{\!\vee} \setminus \h \, {\calH_\h}^{\!\vee}
\, $  and
  $$  0 \;\; \not= \;\; \overline{\,\eta_+} \;\; = \, {\textstyle
\sum\limits_{g(\underline{b}\,) = g_0}} \, c_{\underline{b}}
\, \overline{\x_{\underline{b}}} \;\; = \, {\textstyle
\sum\limits_{g(\underline{b}\,) = g_0}} \, c_{\underline{b}}
\, x_{\underline{b}} \;\; \in \;\, {\calH_\h}^{\!\vee} \Big/
\h \, {\calH_\h}^{\!\vee} \; = \; U(\L_\nu) \; .  $$
Now  $ \, \delta_2(\eta) = 0 \, $  yields  $ \, \delta_2
\big( \overline{\,\eta_+} \,\big) = 0 \, $,  \, thus  $ \,
\sum_{g(\underline{b}\,) = g_0} \, c_{\underline{b}} \,
x_{\underline{b}} = \overline{\,\eta_+} \in P\big(U(\L_\nu)\big)
= \L_\nu \, $;  \, therefore all PBW monomials occurring in the
last sum do belong to  $ B_\nu $  (and  $ \, g_0 = 1 \, $).  In
addition,  $ \, \delta_2(\eta) = 0 \, $  also implies  $ \,
\delta_2(\eta_+) = 0 \, $  which yields also  $ \, \delta \big(
\overline{\,\eta_+} \,\big) = 0 \, $  for the Lie cobracket
$ \delta $  of  $ \L_\nu $  arising as semiclassical limit of
$ \Delta_{{\calH_\h}^{\!\vee}} $  (see Theorem 2.1); therefore
$ \; \overline{\,\eta_+} \, = \sum_{b \in B_\nu} c_b \, x_b \; $  is
an element of  $ \L_\nu $  killed by the Lie cobracket  $ \delta $,
i.e.~$ \, \overline{\,\eta_+} \in \text{\sl Ker}\,(\delta) \, $.
                                             \par
   Now we apply Lemma 4.3 to  $ \, V = V_\nu \, $,  $ \, \L(V) =
\L(V_\nu) =: \L_\nu \, $  and  $ \, \psi = \delta{\big|}_{V_\nu}
\, $,  \, so that  $ \, \psi_{d\L} = \delta \, $.  By the formulas
for  $ \delta $  in Theorem 2.1 we get  $ \, K := \text{\sl Ker}
\,(\psi) = \text{\sl Ker}\,\Big( \delta{\big|}_{V_\nu} \Big) =
\text{\sl Span}\,\big(\{x_1,x_2\}\big) \, $,  \, hence  $ \, \L(K)
= \L\big(\text{\sl Span}\,\big(\{x_1,x_2\}\big)\big) = \text{\sl
Span}\,\Big( \big\{\, x_b \,\big|\, b \in B_\nu \, ; \, \tau(b)
= 1 \,\big\} \Big) \, $,  \, thus eventually (via Theorem 2.1)
$ \, \text{\sl Ker}\,(\delta) = \L(K) = \text{\sl Span}\,\big(
\big\{\, x_b \,\big|\, b \in B_\nu \, ; \, \tau(b) = 1 \,\big\}
\big) \, $.
                                             \par
   As  $ \, \overline{\,\eta_+} \in \text{\sl Ker}\,(\delta) =
\text{\sl Span}\,\big( \big\{\, x_b \,\big|\, b \in B_\nu \, ; \,
\tau(b) = 1 \,\big\} \big) $,  we have  $ \, \overline{\,\eta_+} =
%
%
 \! \sum_{b \in B_\nu , \tau(b)=1} \hskip-1pt c_b \,
x_b \, $;  \, but  $ \, c_b = 0 \, $  whenever  $ \, \tau(b) \leq 1 \, $,
\, by construction of  $ \eta \, $:  thus  $ \, \overline{\,\eta_+} = 0
\, $,  \, a contradiction.  The outcome is  $ \, \eta = 0 \, $,  \,
whence finally  $ \, \eta' \in \varTheta_1 \, $,  \, q.e.d.
 \vskip3pt
   $ \underline{n>1} \, $:  \, We must show that  $ \, D_n = \varTheta_n
\, $,  \, while assuming by induction that  $ \, D_m = \varTheta_m \, $
for all  $ \, m < n \, $.  Let  $ \, \eta = \sum_{\underline{b} \in
\Bbb{B}} c_{\underline{b}} \, \boldalpha_{\underline{b}} \in D_n \, $;
\, then  $ \, \tau(\eta) = \max \big\{\, \tau(\underline{b}\,) \,\big|\,
c_{\underline{b}} \not= 0 \,\big\} \, $.  If  $ \, \delta_2(\eta) =
0 \, $  then  $ \, \eta \in D_1 = \varTheta_1 \, $  by the previous
analysis, and we're done.  Otherwise,  $ \, \delta_2(\eta) \not= 0
\, $  and  $ \, \tau\big(\delta_2(\eta)\big) = \tau(\eta) \, $  by
Lemma 4.4.  On the other hand, since  $ \underline{D} $  is a Hopf
algebra filtration we have  $ \, \delta_2(\eta) \in \sum_{\Sb  r+s=n
\\   r, s > 0  \endSb} D_r \otimes D_s = \sum_{\Sb  r+s=n  \\   r, s
> 0  \endSb} \varTheta_r \otimes \varTheta_s \, $,  \, thanks to
the induction; but then  $ \, \tau\big(\delta_2(\eta)\big) \leq n \, $,
\, by definition of  $ \tau $.  Thus  $ \, \tau(\eta) = \tau \big(
\delta_2(\eta) \big) \leq n \, $,  \, which means  $ \, \eta \in
\varTheta_n \, $.   \qed
\enddemo

\vskip7pt

\proclaim{Theorem 4.1} \, For any  $ \, b \in B_\nu \, $  set
$ \; \widehat{\boldalpha}_b := \h^{\,\kappa(\boldalpha_b)} \,
\boldalpha_b = \h^{\,\tau(b)} \, \boldalpha_b \; $.
                                 \hfill\break
   \indent   (a) \, The set of ordered monomials
  $$  \widehat{\Cal{A}}_{\leq n} \; := \; \Big\{\,
\widehat{\boldalpha}_{\underline{b}} := \widehat{\boldalpha}_{b_1}
\cdots \widehat{\boldalpha}_{b_k} \;\Big|\;  k \in \N \, , \, b_1,
\dots, b_k \in B \, , \, b_1 \preceq \cdots \preceq b_k \, , \,
\kappa\,(\boldalpha_{\underline{b}}\,) = \tau(\underline{b}\,)
\leq n \,\Big\}  $$
is a  $ \, \Bbbk[\h\,] $--basis  of  $ \, D'_n = D_n \big(
{\calH_\h}^{\!\prime} \big) = \h^n D_n \, $.  So  $ \; \widehat{\Cal{A}}
:= \bigcup_{n \in \N} \widehat{\Cal{A}}_{\leq n} \; $  is a  $ \,
\Bbbk[\h\,] $--basis  of  $ \, {\calH_\h}^{\!\prime} \, $.
                                 \hfill\break
   \indent   (b)  $ \; \displaystyle{ {\calH_\h}^{\!\prime}  \;
= \;  \Bbbk[\h\,] \, \Big\langle {\big\{\, \widehat{\boldalpha}_b
\,\big\}}_{b \in B_\nu} \Big\rangle \! \Bigg/ \!\! \bigg( \! \Big\{\,
\big[ \widehat{\boldalpha}_{b_1}, \widehat{\boldalpha}_{b_2} \big]
- \h \, \widehat{\boldalpha}_{[b_1,b_2]} \;\Big|\; \forall \; b_1,
b_2 \in B_\nu \,\Big\} \bigg) } \, $.
                                 \hfill\break
  \indent   (c) \,  $ {\calH_\h}^{\!\prime} $  is a graded Hopf
$ \; \Bbbk[\h\,] $--subalgebra  of  $ \, {\calH_\h} \, $.
                                        \hfill\break
  \indent   (d) \,  $ {\calH_\h}^{\!\prime}{\Big|}_{\h=0}
:= {\calH_\h}^{\!\prime} \Big/ \h \, {\calH_\h}^{\!\prime}
= \widetilde{\calH} = F \big[ {\varGamma_{\!\L_\nu}
\phantom{|}}^{\hskip-8pt \star} \big] \, $,  \, where
$ \, {\varGamma_{\!\L_\nu}\phantom{|}}^{\hskip-8pt\star} \, $
is a connected Poisson algebraic group with cotangent Lie bialgebra
isomorphic to  $ \L_\nu $  (as a Lie algebra) with the graded Lie
bialgebra structure given by  $ \, \delta(x_n) = (n-2) \, x_{n-1}
\wedge x_1 \, $  (for all  $ \, n \in \N_\nu $).  Indeed,  $ \,
{\calH_\h}^{\!\prime}{\Big|}_{\h=0} $  is the free Poisson
(commutative) algebra over  $ \N_\nu \, $,  generated by all
the\/  $ \bar{\boldalpha}_n := \widehat{\boldalpha}_n{\big|}_{\h=0} $
($ \, n \in \N_\nu \, $)  with Hopf structure given (for all  $ \,
n \in \N_\nu $)  by
  $$  \displaylines{
   \Delta\big(\bar{\boldalpha}_n\big)  \; = \;\,
\bar{\boldalpha}_n \otimes 1 \, + \, 1 \otimes
\bar{\boldalpha}_n \, + \, {\textstyle \sum_{k=2}^{n-1} {n \choose k}}
\, \bar{\boldalpha}_k \otimes {\bar{\boldalpha}_1}^{\,n-k} \,
+ \, {\textstyle \sum_{k=1}^{n-1}} \, (k+1) \,
{\bar{\boldalpha}_1}^{\,k} \otimes \bar{\boldalpha}_{n-k}  \cr
   S\big(\bar{\boldalpha}_n\big)  \; = \;\,
- \, \bar{\boldalpha}_n \, - \, {\textstyle \sum_{k=2}^{n-1}
{n \choose k}} \, S\big(\bar{\boldalpha}_k\big) \,
{\bar{\boldalpha}_1}^{\,n-k} \, - \, {\textstyle
\sum_{k=1}^{n-1}} \, (k+1) \, S \big( \bar{\boldalpha}_1
\big)^k \, \bar{\boldalpha}_{n-k} \; ,  \qquad
\epsilon\big(\bar{\boldalpha}_n\big) \, = \, 0  \; .  \cr }  $$
Thus  $ \, {\calH_\h}^{\!\prime}{\Big|}_{\h=0} $  is the polynomial
algebra  $ \, \Bbbk \big[ {\{\, \eta_b \,\}}_{b \in B_\nu} \big] \, $
generated by a set of indeterminates  $ \, {\{\, \eta_b \,\}}_{b \in
B_\nu} \, $  in bijection with  $ B_\nu \, $,  so  $ \; {\varGamma_{\!
\L_\nu} \phantom{|}}^{\hskip-8pt \star} \cong \Bbb{A}_\Bbbk^{B_\nu}
\, $
%
%
as algebraic varieties.
                                        \hfill\break
  \indent   Finally,  $ \, {\calH_\h}^{\!\prime}{\Big|}_{\h=0} =
F \big[ {\varGamma_{\!\L_\nu} \phantom{|}}^{\hskip-8pt \star} \big]
= \Bbbk \big[ {\{\, \eta_b \,\}}_{b \in B_\nu} \big] \, $  is a
{\sl graded Poisson Hopf algebra}  w.r.t.~the grading
$ \, \partial(\bar{\boldalpha}_n) = n \, $  (inherited from
$ {\calH_\h}^{\!\prime} $)  and w.r.t.~the grading induced
from  $ \, \kappa = \tau \, $  (on  $ \calH $),  and a  {\sl
graded algebra}  w.r.t.~the  ({\sl polynomial})  grading  $ \,
d(\bar{\boldalpha}_n) = 1 \, $  (for all  $ \, n \in \N_+ $).
                                        \hfill\break
  \indent   (e) \,  The analogues of statements (a)--(d) hold with
$ \calK $  instead of  $ \, \calH \, $,  \, with  $ X^+ $  instead
of  $ X $  for all  $ \, X = \L_\nu, B_\nu, \N_\nu \, $,  and with
$ {\varGamma_{\!\L_\nu^+}\phantom{|}}^{\hskip-8pt \star} $  instead
of  $ {\varGamma_{\!\L_\nu}\phantom{|}}^{\hskip-8pt \star} \, $.
\endproclaim

\demo{Proof} \, {\it (a)} \, This follows from Proposition 4.2 and
the definition of  $ {\calH_\h}^{\!\prime} $  in \S 4.2.
                                              \par
   {\it (b)} \, This is a direct consequence of claim  {\it (a)\/}
and  Lemma 4.2{\it (g)}.
                                              \par
   {\it (c)} \, Thanks to claims  {\it (a)\/}  and  {\it (b)},  we
can look at  $ {\calH_\h}^{\!\prime} $  as a Poisson algebra, whose
Poisson bracket is given by  $ \, \{x,y\,\} := \h^{-1} [x,y] = \h^{-1}
(x \, y - y \, x) \, $  (for all  $ \, x $,  $ y \in {\calH_\h}^{\!
\prime} \, $);  then  $ {\calH_\h}^{\!\prime} $  itself is the
free associative Poisson algebra generated by  $ \big\{\,
\widehat{\boldalpha}_n \,\big|\, n \in \N \,\big\} $.  Clearly
$ \Delta $  is a Poisson map, therefore it is enough to prove that
$ \, \Delta\big(\widehat{\boldalpha}_n\big) \in {\calH_\h}^{\!\prime}
\otimes {\calH_\h}^{\!\prime} \, $  for all  $ \, n \in \N_+ \, $.
This is clear for  $ \boldalpha_1 $  and  $ \boldalpha_2 $  which
are primitive; as for  $ \, n > 2 \, $,  \, we have, like in
Proposition 4.1,
  $$  \hbox{ $ \eqalign{
   \hskip-2pt   \Delta \big( \widehat{\boldalpha}_n  &  \big)  \;
= \;  {\textstyle \sum_{k=2}^n} \, \h^{k-1} \boldalpha_k \otimes
\h^{n-k} Q^k_{n-k}(\a_*) \, + \, {\textstyle \sum_{k=0}^{n-1}} \,
\h^k {\boldalpha_1}^{\!k} \otimes \h^{n-k-1} Z^k_{n-k}(\boldalpha_*)
\; =   \hfill  \cr
   {} \hfill   &  \hskip-3pt = \;  {\textstyle \sum_{k=2}^n} \,
\widehat{\boldalpha}_k \otimes \h^{n-k} Q^k_{n-k}(\a_*) \, + \,
{\textstyle \sum_{k=0}^{n-1}} \, {\widehat{\boldalpha}_1}^{\,k} \otimes
\h^{n-k-1} Z^k_{n-k}(\boldalpha_*) \; \in \; {\calH_\h}^{\!\prime}
\otimes {\calH_\h}^{\!\prime}  \cr } $ }   \hfill (4.5)  $$
thanks to Lemma 4.1 (with notations used therein).  In addition,
$ \, S\big({\calH_\h}^{\!\prime}\big) \subseteq {\calH_\h}^{\!\prime}
               \, $  also\break
 \noindent
follows by induction from (4.5) because Hopf algebra
axioms along with (4.5) give
  $$  S\big(\widehat{\boldalpha}_n\big)  \; = \;
- \widehat{\boldalpha}_n \, - \, {\textstyle \sum_{k=2}^{n-1}} \,
S\big(\widehat{\boldalpha}_k\big) \, \h^{n-k} Q^k_{n-k}(\a_*) \, - \,
{\textstyle \sum_{k=1}^{n-1}} \, S \big( {\widehat{\boldalpha}_1}^{\,k}
\big) \, \h^{n-k-1} Z^k_{n-k}(\boldalpha_*) \; \in \; {\calH_\h}^{\!
\prime}  $$
for all  $ \, n \in \N_\nu \, $  (using induction).  The claim
follows.
                                              \par
   {\it (d)} \, Thanks to  {\it (a)\/}  and  {\it (b)},  $ \,
{\calH_\h}^{\!\prime} {\Big|}_{\h=0} $  is a polynomial
$ \Bbbk $--algebra  as claimed, over the set of indeterminates
$ \Big\{ \bar{\boldalpha}_b := \widehat{\boldalpha}_b {\big|}_{\h=0}
\, \big(\! \in \! {\calH_\h}^{\!\prime}{\big|}_{\h=0} \big) \Big\}_{b
\in B_\nu} \, $.  Furthermore, in the proof of  {\it (c)\/}  we noticed
that  $ {\calH_\h}^{\!\prime} $  is also the free Poisson algebra
generated by  $ \big\{\, \widehat{\boldalpha}_n \,\big|\, n \in \N
\,\big\} $;  therefore  $ \, {\calH_\h}^{\!\prime}{\Big|}_{\h=0} $
is the free commutative Poisson algebra generated by  $ \big\{
\bar{\boldalpha}_n := \check{\boldalpha}_{x_n} {\big|}_{\h=0}
\,\big\}_{n \in \N} \, $.  Then formula (4.5)   --- for all 
$ \, n \in \N_\nu \, $  ---   describes uniquely the Hopf
structure of  $ {\calH_\h}^{\!\prime} $,  hence the formula
it yields at  $ \, \h = 0 \, $  will describe the Hopf
structure of  $ \, {\calH_\h}^{\!\prime}{\big|}_{\h=0} $.
                                              \par
   Expanding  $ \, \h^{n-k} Q^k_{n-k}(\a_*) \, $  in (4.5) w.r.t.~the
basis  $ \widehat{\calA} $  in  {\it (a)\/}  we find a sum of terms
of  $ \tau $--degree  less or equal than  $ (n-k) $,  and the sole
one achieving equality is  $ \, {\widehat{\boldalpha}_1}^{\,n-k} \, $,
\, which occurs with coefficient  $ {n \choose k} $:  \, similarly,
when expanding  $ \, \h^{n-k-1} Z^k_{n-k}(\boldalpha_*) \, $  in (4.5)
w.r.t.~$ \widehat{\calA} $  all summands have  $ \tau $--degree  less
or equal than  $ (n-k-1) $,  and equality holds only for  $ \,
\widehat{\boldalpha}_{n-k} \, $,  \, whose coefficient is  $ \,
(k+1) \, $.  Therefore for some  $ \, \boldsymbol{\eta} \in
{\calH_\h}^{\!\prime} {\big|}_{\h=0} \; $  we have
  $$  \Delta\big(\widehat{\boldalpha}_n\big) \; = \;
{\textstyle \sum_{k=2}^n} \, \widehat{\boldalpha}_k \otimes
{\textstyle {n \choose k}} \, {\widehat{\boldalpha}_1}^{\,n-k}
\, + \, {\textstyle \sum_{k=0}^{n-1}} \, (k+1) \,
{\widehat{\boldalpha}_1}^{\,k} \otimes \widehat{\boldalpha}_{n-k}
\, + \, \h \; \boldsymbol{\eta} \; ;  $$
this yields the formula for  $ \Delta $,  from which the formula
for  $ S $  follows too as usual.
                                              \par
   Finally, let  $ \, \varGamma := \text{\sl Spec}\,\big(
{\calH_\h}^{\!\prime}{\big|}_{\h=0} \big) \, $  be the
algebraic Poisson group such that  $ \, F\big[\varGamma\big]
= {\calH_\h}^{\!\prime}{\big|}_{\h=0} \, $,  \, and let  $ \,
\boldsymbol{\gamma}_\nu := \text{\sl coLie}\,(\varGamma) \, $
be its cotangent Lie bialgebra.  Since  $ \, {\calH_\h}^{\!
\prime}{\big|}_{\h=0} $  is Poisson free over  $ \big\{
\bar{\boldalpha}_n \big\}_{n \in \N_\nu} \, $,  \,
as a Lie algebra  $ \boldsymbol{\gamma}_\nu $  is free over
$ \, \big\{\, d_n := \bar{\boldalpha}_n \mod \germ^2
\,\big\}_{n \in \N_\nu} \, $  (where  $ \, \germ :=
J_{{\calH_\h}^{\!\prime}{|}_{\h=0}} \, $),  so  $ \,
\boldsymbol{\gamma}_\nu \cong \L_\nu \, $,  via  $ \,
d_n \mapsto x_n (n \in \N_+) \, $  as a Lie algebra.
The Lie cobracket is
  $$  \displaylines{
   {} \;   \delta_{\boldsymbol{\gamma}_\nu}\big(d_n\big)  \,
= \,  (\Delta - \Delta^{\text{op}})\big(\bar{\boldalpha}_n\big)
\mod \germ_\otimes  \; = \;\,  {\textstyle \sum_{k=2}^{n-1}
{n \choose k}} \, \bar{\boldalpha}_k \wedge {\bar{\boldalpha}_1}^{\,n-k}
\, +   \hfill   \cr
   + \, {\textstyle \sum_{k=1}^{n-1}} \, (k+1) \,
{\bar{\boldalpha}_1}^{\,k} \wedge \bar{\boldalpha}_{n-k}
\mod \germ_\otimes  \; = \; {\textstyle {n \choose n-1}}
\, \bar{\boldalpha}_{n-1} \wedge \bar{\boldalpha}_1 \, +
\, 2 \, \bar{\boldalpha}_1 \wedge \bar{\boldalpha}_{n-1}
\mod \germ_\otimes  \; = \;  \cr
   {} \hfill   = \;  (n-2) \, \bar{\boldalpha}_{n-1} \wedge
\bar{\boldalpha}_1  \mod \germ_\otimes  \; = \; (n-2) \,
d_{n-1} \wedge d_1  \; \in \;  \boldsymbol{\gamma} \otimes
\boldsymbol{\gamma}  \cr }  $$
where  $ \, \germ_\otimes := \Big( \germ^2 \otimes {\calH_\h}^{\!
\prime}{|}_{\h=0} + \germ \otimes \germ + {\calH_\h}^{\!\prime}
{|}_{\h=0} \otimes \germ^2 \Big) \, $,  \, whence  $ \, \varGamma
= {\varGamma_{\!\L_\nu}\phantom{|}}^{\hskip-8pt \star} \, $  as
claimed in  {\it (d)}.
                                              \par
   Finally, the statements about gradings of  $ \, {\calH_\h}^{\!\prime}
{\Big|}_{\h=0} = F \big[ {\varGamma_{\!\L_\nu} \phantom{|}}^{\hskip-8pt
\star} \big] \, $  hold by construction.
                                              \par
   {\it (e)} \, This should be clear from the whole discussion, since
all arguments apply again   --- {\sl mutatis mutandis\/} ---   when
starting with  $ \calK $  instead of  $ \, \calH \, $;  \, we leave
details to the reader.   \qed
\enddemo

\vskip1,5truecm

\centerline {\bf \S \; 5 \ Drinfeld's deformation
$ \, \big( {\calH_\h}^{\!\prime} \big)^\vee \, $. }

\vskip10pt

  {\bf 5.1 The goal.} \, Like in \S 3.1, there is a second step in
the crystal duality principle which builds another deformation basing
upon the Rees deformation  $ {\calH_\h}^{\!\prime} $.  This will be again
a Hopf  $ \Bbbk[\h] $--algebra,  namely  $ \big( {\calH_\h}^{\!\prime}
\big)^{\!\vee} $,  which specializes to  $ \calH $  for  $ \, \h = 1 \, $
and for  $ \, \h = 0 \, $  instead specializes to  $ U(\gerk_-) $,  for
some Lie bialgebra  $ \gerk_- \, $.  In other words,  $ \, \big(
{\calH_\h}^{\!\prime} \big)^{\!\vee}\Big|_{\h=1} \!\! = \calH \, $
and  $ \, \big( {\calH_\h}^{\!\prime} \big)^{\!\vee}\Big|_{\h=0} \!\!
= U(\gerk_-) \, $,  \; the latter meaning that  $ \big( {\calH_\h}^{\!
\prime} \big)^{\!\vee} $  is a  {\sl quantized universal enveloping
algebra\/}  (QUEA in the sequel).  Thus  $ \big( {\calH_\h}^{\!\prime}
\big)^{\!\vee} $  is a  {\sl quantization\/}  of  $ U(\gerk_-) \, $,
and the quantum symmetry  $ \calH $  is a deformation of the
classical Poisson symmetry  $ U(\gerk_-) \, $.
                                           \par
   The general theory describes explicitly the relationship between
$ \gerk_- $  and  $ {\varGamma_{\!\L_\nu} \phantom{|}}^{\hskip-8pt
\star} $  in \S 4, which is  $ \, \gerk_- = \boldsymbol{\gamma}_\nu
:= \text{\sl coLie}\,\big( {\varGamma_{\!\L_\nu}\phantom{|}}^{\hskip-8pt
\star} \big) \cong \L_\nu  $  (with the structure in  Theorem 4.1{\it
(d)\/}),  \, the cotangent Lie bialgebra of  $ {\varGamma_{\!\L_\nu}
\phantom{|}}^{\hskip-8pt \star} \, $.  Thus, from this and \S 4 we
see that the quantum symmetry encoded by  $ \calH $  is (also)
intermediate between the two classical,
          \hbox{Poisson symmetries ruled by
$ {\varGamma_{\!\L_\nu}\phantom{|}}^{\hskip-8pt \star} $
and  $ \boldsymbol{\gamma}_\nu \, $.}
                                       \par
   In this section I describe explicitly  $ \big( {\calH_\h}^{\!\prime}
\big)^{\!\vee} $  and its semiclassical limit  $ U(\gerk_-) \, $,  \,
hence  $ \gerk_- $  itself too.  This provides a direct proof of the
above mentioned results on  $ \big( {\calH_\h}^{\!\prime} \big)^{\!
\vee} $  and  $ \gerk_- \, $.

\vskip7pt

  {\bf 5.2 Drinfeld's algebra  $ \big( {\calH_\h}^{\!\prime} \big)^{\!
\vee} $.} \, Let  $ \, J^{\,\prime} := J_{{\calH_\h}^{\!\prime}} \, $,
\; and define
  $$  \big( {\calH_\h}^{\!\prime} \big)^{\!\vee} \; := \; {\textstyle
\sum_{n \in \N}} \, \h^{-n} {J^{\,\prime}}^n \; = \; {\textstyle
\sum_{n \in \N}} \, \big( \h^{-1} {J^{\,\prime}} \big)^n  \quad
\big( \subseteq \calH(\h) \,\big) .   \eqno (5.1)  $$
   \indent   Now I describe  $ \big( {\calH_\h}^{\!\prime}
\big)^{\!\vee} $  and its specializations at  $ \, \h = 1 \, $
and  $ \, \h = 0 \, $.  The main step is

\vskip7pt

\proclaim{Theorem 5.1} \, For any  $ \, b \in B_\nu \, $  set
$ \; \check{\boldalpha}_b := \h^{\,\kappa(\boldalpha_b)-1}
\, \boldalpha_b = \h^{\,\tau(b)-1} \, \boldalpha_b = \h^{-1}
\, \widehat{\boldalpha}_b \; $.
                                 \hfill\break
   \indent   (a)  $ \; \displaystyle{ {\big( {\calH_\h}^{\!\prime}
\big)}^{\!\vee}  \; = \;  \Bbbk[\h\,] \, \Big\langle {\big\{\,
\check{\boldalpha}_b \,\big\}}_{b \in B_\nu} \Big\rangle \!
\Bigg/ \!\! \bigg( \! \Big\{\, \big[ \check{\boldalpha}_{b_1},
\check{\boldalpha}_{b_2} \big] - \, \check{\boldalpha}_{[b_1,b_2]}
\;\Big|\; \forall \; b_1, b_2 \in B_\nu \,\Big\} \bigg) } \, $.
                                 \hfill\break
  \indent   (b) \,  $ {\big( {\calH_\h}^{\!\prime} \big)}^{\!\vee} $
is a graded Hopf  $ \; \Bbbk[\h\,] $--subalgebra  of  $ \, \calH_\h
\, $.
                                        \hfill\break
  \indent   (c) \,  $ {\big( {\calH_\h}^{\!\prime} \big)}^{\!\vee}
{\Big|}_{\h=0} := {\big( {\calH_\h}^{\!\prime} \big)}^{\!\vee} \Big/
\h \, {\big( {\calH_\h}^{\!\prime} \big)}^{\!\vee} \cong U \big( \L_\nu
\big) \; $  as co-Poisson Hopf algebra, where  $ \, \L_\nu \, $  bears
the Lie bialgebra structure given by  $ \, \delta(x_n) = (n-2) \,
x_{n-1} \wedge x_1 \, $  (for all  $ \, n \in \N_\nu $).
                                        \hfill\break
   \indent   Finally, the grading  $ \, d $  given by  $ \, d(x_n)
:= 1 \;\, (n \in \N_+) \, $  makes  $ \, {\big( {\calH_\h}^{\!\prime}
\big)}^{\!\vee}{\Big|}_{\h=0} \! = U(\L_\nu) \, $  into a  {\sl graded
co-Poisson Hopf algebra},  and the grading  $ \, \partial $  given
by  $ \, \partial(x_n) := n \;\, (n \in \N_+) \, $  makes  $ \, {\big(
{\calH_\h}^{\!\prime} \big)}^{\!\vee}{\Big|}_{\h=0} \! = U(\L_\nu)
\, $     into a  {\sl graded Hopf algebra}  and  $ \L_\nu $  into
a  {\sl graded Lie bialgebra.}
                                        \hfill\break
  \indent   (d) \,  The analogues of statements (a)--(c) hold
with  $ \calK \, $,  $ \L_\nu^+ $,  $ B_\nu^+ $  and  $ \N_\nu^+ $
respectively instead of  $ \, \calH \, $,  $ \L_\nu^+ $,  $ B_\nu $
and  $ \N_\nu^+ \, $.
\endproclaim

\demo{Proof} \, {\it (a)} \, This follows from  Theorem 4.1{\it
(b)\/}  and the very definition of  $ {\big( {\calH_\h}^{\!\prime}
\big)}^{\!\vee} $  in \S 5.2.
                                              \par
   {\it (b)} \, This is a direct consequence of claim  {\it (a)\/}
and  Theorem 4.1{\it (c)}.
                                              \par
   {\it (c)} \, It follows from claim  {\it (a)\/}  that mapping
$ \, \check{\boldalpha}_b{\big|}_{\h=0} \mapsto b \, $  ($ \forall
\, b \in B_\nu \, $)  yields a well-defined algebra isomorphism
$ \, \Phi \, \colon \, {\big( {\calH_\h}^{\!\prime}
\big)}^{\!\vee}{\Big|}_{\h=0} {\buildrel \cong \over
{\lhook\joinrel\relbar\joinrel\relbar\joinrel\twoheadrightarrow\,}}
U\big(\L_\nu) \, $.  In addition, when expanding  $ \, \h^{n-k}
Q^k_{n-k}(\a_*) \, $  in (4.5) w.r.t.~the basis  $ \calA $  (see
Proposition 4.2) we find a sum of terms of  $ \tau $--degree  less
or equal than  $ (n-k) $,  and equality is achieved only for  $ \,
\boldalpha_1^{\,n-k} \, $,  \, which occurs with coefficient
$ {n \choose k} $:  \, similarly, the expansion of  $ \, \h^{n-k-1}
Z^k_{n-k}(\boldalpha_*) \, $  in (4.5) yields a sum of terms whose
$ \tau $--degree  is less or equal than  $ (n-k-1) $,  with equality
only for  $ \, \boldalpha_{n-k} \, $,  \, whose coefficient is  $ \,
(k+1) \, $.   Thus using the relation  $ \, \widehat{\boldalpha}_s
= \h \, \check{\boldalpha}_s \, $  ($ \, s \in \N_+ \, $)  we get
  $$  \displaylines{
   \Delta \big( \check{\boldalpha}_n \big)  \, = \,
\check{\boldalpha}_n \otimes 1  \, + \,  1 \otimes \check{\boldalpha}_n
\, + \,  {\textstyle \sum_{k=2}^{n-1}} \, \check{\boldalpha}_k \otimes
\h^{n-k} Q^k_{n-k}(\a_*)  \, + \, {\textstyle \sum_{k=1}^{n-1}} \,
{\check{\boldalpha}_1}^{\,k} \otimes \h^{n-1} Z^k_{n-k}(\boldalpha_*)
\, =   \hfill {}  \cr
   = \,  \check{\boldalpha}_n \otimes 1  \, + \,  1 \otimes
\check{\boldalpha}_n  \, + \, {\textstyle \sum_{k=2}^{n-1}} \,
\h^{n-k} \, \check{\boldalpha}_k \otimes {\textstyle {n \choose k}} \,
{\check{\boldalpha}_1}^{\,n-k}  \, + \, {\textstyle \sum_{k=1}^{n-1}}
\, \h^k \, (k+1) \, {\check{\boldalpha}_1}^{\,k} \otimes
\check{\boldalpha}_{n-k}  \, + \,  \h^2 \; \boldsymbol{\eta} \, =  \cr
   {} \hfill   = \, \check{\boldalpha}_n \otimes 1  \, + \,  1 \otimes
\check{\boldalpha}_n  \, + \,  \h \, \big( n \, \check{\boldalpha}_{n-1}
\otimes \check{\boldalpha}_1  \, + \, 2 \, \check{\boldalpha}_1 \otimes
\check{\boldalpha}_{n-1} \big)  \, + \,  \h^2 \; \boldsymbol{\chi}
\cr }  $$
for some  $ \, \boldsymbol{\eta}, \boldsymbol{\chi} \in {\big(
{\calH_\h}^{\!\prime} \big)}^{\!\vee} \otimes {\big( {\calH_\h}^{\!
\prime} \big)}^{\!\vee} \, $.  It follows that  $ \, \Delta
\big( \check{\boldalpha}_n{\big|}_{\h=0} \big)  \, = \,
\check{\boldalpha}_n{\big|}_{\h=0} \otimes 1  \, + \,
1 \otimes \check{\boldalpha}_n{\big|}_{\h=0} \, $  for
all  $ \, n \in\N_\nu \, $.  Similarly we have  $ \, S \big(
\check{\boldalpha}_n{\big|}_{\h=0} \big) = - \check{\boldalpha}_n
{\big|}_{\h=0} \, $  and  $ \, \epsilon\big(\check{\boldalpha}_n
{\big|}_{\h=0}\big) = 0 \, $  for all  $ \, n \in \N_\nu \, $,
\, thus  $ \Phi $  {\sl is an isomorphism of Hopf algebras\/}  too.
In addition, the Poisson cobracket of  $ {\big( {\calH_\h}^{\!\prime}
\big)}^{\!\vee}{\Big|}_{\h=0} $  inherited from  $ {\big( {\calH_\h}^{\!
\prime} \big)}^{\!\vee} $  is given by
  $$  \displaylines{
   \delta\big(\check{\boldalpha}_n{\big|}_{\h=0}\big)
\, = \,  \Big( \h^{-1} (\Delta - \Delta^{\text{op}}) \big(
\check{\boldalpha}_n \big) \Big)  \mod \h \, {\big( {\calH_\h}^{\!
\prime} \big)}^{\!\vee} \otimes {\big( {\calH_\h}^{\!\prime}
\big)}^{\!\vee}  \; =   \hfill  \cr
   {} \hfill   = \;\,  \big( n \, \check{\boldalpha}_{n-1} \wedge
\check{\boldalpha}_1  \, + \, 2 \, \check{\boldalpha}_1 \wedge
\check{\boldalpha}_{n-1} \big)  \mod \h \, {\big( {\calH_\h}^{\!
\prime} \big)}^{\!\vee} \! \otimes {\big( {\calH_\h}^{\!\prime}
\big)}^{\!\vee}  \, = \;  (n-2) \, \check{\boldalpha}_{n-1}
{\big|}_{\h=0} \wedge \check{\boldalpha}_1{\big|}_{\h=0}  \cr }  $$
hence  $ \Phi $  is also an isomorphism of  {\sl co-Poisson\/}
Hopf algebras, as claimed.
                                              \par
   The statements on gradings of  $ \; {\big( {\calH_\h}^{\!\prime}
\big)}^{\!\vee}{\Big|}_{\h=0} = U(\L_\nu) \; $  should be clear by
construction.
%
%
 \eject
   {\it (d)} \, This should be clear from the whole discussion, as
all arguments apply again   --- {\sl mutatis mutandis\/} ---   when
starting with  $ \calK $  instead of  $ \, \calH \, $;  \, details
are left to the reader.   \qed
\enddemo

\vskip7pt

   {\bf 5.3 Specialization limits.} \, So far, Theorem 4.1{\it
(d)\/}  and  Theorem 5.1{\it (c)\/}  prove the following
specialization results for  $ {\calH_\h}^{\!\prime} $  and
$ {\big( {\calH_\h}^{\!\prime} \big)}^{\!\vee} $  respectively:
  $$  {\calH_\h}^{\!\prime} \;{\buildrel \h \rightarrow 0
\over \llongrightarrow}\; F \big[ {\varGamma_{\!\L_\nu}
\phantom{|}}^{\hskip-8pt\star} \big] \quad ,  \qquad \qquad
{\big( {\calH_\h}^{\!\prime} \big)}^{\!\vee} \;{\buildrel
\h \rightarrow 0 \over \llongrightarrow}\; U(\L_\nu)  $$
as graded Poisson or co-Poisson Hopf  $ \Bbbk $--algebras.
In addition,  Theorem 4.1{\it (b)\/}  implies that  $ \;
{\calH_\h}^{\!\prime} \,{\buildrel \h \rightarrow 1 \over
\llongrightarrow}\, \calH' = \calH \; $  as graded Hopf
$ \Bbbk $--algebras.  Indeed, by  Theorem 4.1{\it (b)\/}
$ \, \calH $  (or even  $ \calH_\h $)  embeds as an algebra
into  $ {\calH_\h}^{\!\prime} $,  via  $ \; \boldalpha_n
\mapsto \widehat\boldalpha_n \; $  (for all  $ \, n \in
\N_\nu \, $):  \, then
  $$  [\boldalpha_n, \boldalpha_m] \, \mapsto \, \big[
\widehat\boldalpha_n, \widehat\boldalpha_m \big]  \, =
\,  \h \, \widehat{\boldalpha}_{[x_n,x_m]} \, \equiv \,
\widehat{\boldalpha}_{[x_n,x_m]} \mod (\h \! - \! 1) \,
{\calH_\h}^{\!\prime}   \eqno \big(\, \forall \; n, m
\in \N_\nu \big)  $$
thus, thanks to the presentation of  $ {\calH_\h}^{\!\prime} $
in  Theorem 4.1{\it (b)\/},  $ \calH $  is isomorphic
to  $ \; {\calH_\h}^{\!\prime}{\Big|}_{\h=1} := \;
{\calH_\h}^{\!\prime} \Big/ (\h \! - \! 1) \, {\calH_\h}^{\!\prime}
\; = \;  \Bbbk \big\langle {\widehat\boldalpha}_1{\big|}_{\h=1},
{\widehat\boldalpha}_2{\big|}_{\h=1}, \dots, {\widehat\boldalpha}_n
{\big|}_{\h=1}, \ldots \big\rangle \, $,  \, as a  $ \Bbbk $--algebra,
via  $ \, \boldalpha_n \mapsto {\widehat\boldalpha}_n{\big|}_{\h=1}
\, $.  Moreover, the Hopf structure of  $ {\calH_\h}^{\!\prime}
{\Big|}_{\h=1} $  is given by
  $$  \Delta \big( \widehat{\boldalpha}_n{\big|}_{\h=1} \big)
=  {\textstyle \sum_{k=2}^n} \, \widehat{\boldalpha}_k \otimes
\h^{n-k} Q^k_{n-k}(\a_*)  +  {\textstyle \sum_{k=0}^{n-1}} \,
{\widehat{\boldalpha}_1}^{\,k} \otimes \h^{n-1} Z^k_{n-k}
(\boldalpha_*) \hskip-6pt  \mod (\h-1) {\calH_\h}^{\!\prime}
\otimes {\calH_\h}^{\!\prime} \, .  $$
   \indent   Now,  $ \, Q^k_{n-k}(\a_*) = Q^k_{n-k}(\boldalpha_*
+ {\boldalpha_1}^*) = {\Cal Q}^k_{n-k}(\boldalpha_*) \, $  for
some polynomial  $ {\Cal Q}^k_{n-k}(\boldalpha_*) $  in the
$ \boldalpha_i $'s{\,};  let  $ \, {\Cal Q}^k_{n-k}(\boldalpha_*)
= \sum_s {\Cal T}^{s,k}_{n-k}(\boldalpha_*) \, $  be the splitting
of  $ {\Cal Q}^k_{n-k} $  into  $ \tau $--homogeneous  summands
(i.e., each  $ {\Cal T}^{s,k}_{n-k}(\boldalpha_*) $  is a
homogeneous polynomial of  $ \tau $--degree  $ s \, $):  then
  $$  \h^{n-k} Q^k_{n-k}(\a_*)  \, = \,  \h^{n-k} {\Cal Q}^k_{n-k}
(\boldalpha_*)  \, = \,  \h^{n-k} {\textstyle \sum_s}
{\Cal T}^{s,k}_{n-k}(\boldalpha_*)  \, = \,  {\textstyle \sum_s}
\h^{n-k-s} {\Cal T}^{s,k}_{n-k}(\widehat{\boldalpha}_*)  $$
with  $ \, n-k-s > 0 \, $  for all  $ s $  (by construction).  Since
clearly  $ \, \h^{n-k-s}{\Cal T}^{s,k}_{n-k}(\widehat{\boldalpha}_*)
\equiv {\Cal T}^{s,k}_{n-k}(\widehat{\boldalpha}_*) \mod (\h-1) \,
{\calH_\h}^{\!\prime} \, $,  \, we find  $ \; \h^{n-k} \, Q^k_{n-k}
(\a_*) \, = \, \h^{n-k} \, {\Cal Q}^k_{n-k} (\boldalpha_*) \,
= \, {\textstyle \sum_s} \h^{n-k-s} \, {\Cal T}^{s,k}_{n-k}
(\widehat{\boldalpha}_*) \, \equiv \, {\textstyle \sum_s}
{\Cal T}^{s,k}_{n-k}(\widehat{\boldalpha}_*)  \mod (\h-1) \,
{\calH_\h}^{\!\prime} \, = \, {\Cal Q}^k_{n-k}(\widehat{\boldalpha}_*)
\, $,  \, for all  $ k $  and  $ n \, $.  Similarly we argue that
$ \, \h^{n-1} Z^k_{n-k}(\boldalpha_*) \equiv Z^k_{n-k}
(\widehat{\boldalpha}_*) \mod (\h-1) \, {\calH_\h}^{\!\prime} \, $,
\, for all  $ k $  and  $ n \, $.  The outcome is that
  $$  \displaylines{
   \Delta \big( \widehat{\boldalpha}_n{\big|}_{\h=1} \big)  \! = \!
{\textstyle \sum_{k=2}^n} \, \widehat{\boldalpha}_k \otimes \h^{n-k}
\! {\Cal Q}^k_{n-k}(\boldalpha_*)  +  {\textstyle \sum_{k=0}^n} \,
{\widehat{\boldalpha}_1}^{\,k} \otimes \h^{n-1} Z^k_{n-k}(\boldalpha_*)
\hskip-4pt  \mod \hskip-2pt (\h-1) \, {\calH_\h}^{\!\prime} \otimes
{\calH_\h}^{\!\prime}  =  \cr
   {} \hfill   = \,  {\textstyle \sum_{k=2}^{n-1}} \,
\widehat{\boldalpha}_k \otimes {\Cal Q}^k_{n-k}(\widehat{\boldalpha}_*)
+  {\textstyle \sum_{k=0}^{n-1}} \, {\widehat{\boldalpha}_1}^{\,k}
\otimes Z^k_{n-k}(\widehat{\boldalpha}_*) \mod (\h-1) \,
{\calH_\h}^{\!\prime} \otimes {\calH_\h}^{\!\prime} \, .  \cr }  $$
   \indent   On the other hand, we have  $ \; \Delta(\boldalpha_n)
\, = \,  {\textstyle \sum_{k=2}^n} \, \boldalpha_k \otimes
{\Cal Q}^k_{n-k}(\boldalpha_*) \, + \, {\textstyle \sum_{k=0}^{n-1}}
\, \boldalpha_1^{\,k} \otimes Z^k_{n-k}(\boldalpha_*) \, $
in  $ \calH $.  Thus the  {\sl graded algebra\/}  isomorphism 
$ \, \Psi \, \colon \, \calH \,{\buildrel \cong \over
{\lhook\joinrel\relbar\joinrel\relbar\joinrel\twoheadrightarrow}}\,
{\calH_\h}^{\!\prime}{\Big|}_{\h=1} \, $  given by  $ \,
\boldalpha_n \mapsto {\widehat\boldalpha}_n{\big|}_{\h=1} \, $
preserves the coproduct too.  Similarly,  $ \Psi $  respects
the antipode and the counit, hence it is a graded Hopf algebra
isomorphism.  In a nutshell, we have (as graded Hopf
$ \Bbbk $--algebras)  $ \; {\calH_\h}^{\!\prime}
\;{\buildrel {\h \rightarrow 1} \over \llongrightarrow}\;
\calH' = \calH \; $.  Similarly, Theorem 5.1 implies that
$ \, {\big( {\calH_\h}^{\!\prime} \big)}^{\!\vee} \,{\buildrel
\h \rightarrow 1 \over \llongrightarrow}\, \calH \, $  as graded
Hopf  $ \Bbbk $--algebras.  Indeed,  Theorem 5.1{\it (a)\/}  shows
that  $ \, {\big( {\calH_\h}^{\!\prime} \big)}^{\!\vee} \cong
\Bbbk[\h\,] \otimes_\Bbbk U(\L_\nu) \, $  {\sl as graded associative
algebras},  via  $ \, \check\boldalpha_n \mapsto x_n \, $  ($ \, n
\in \N_\nu \, $),  \, in particular  $ {\big( {\calH_\h}^{\!\prime}
\big)}^{\!\vee} $  is the free associative  $ \Bbbk[\h\,] $--algebra
over  $ \big\{ \check{\boldalpha}_n \big\}_{n \in \N_\nu} $;  \, then
specialization yields a graded algebra isomorphism
  $$  \Omega \, \colon \, {\big( {\calH_\h}^{\!\prime} \big)}^{\!\vee}
{\Big|}_{\h=1} := \, {\big( {\calH_\h}^{\!\prime} \big)}^{\!\vee} \!
\Big/ (\h \! - \! 1) \, {\big( {\calH_\h}^{\!\prime} \big)}^{\!\vee}
\,{\buildrel \cong \over
{\lhook\joinrel\relbar\joinrel\relbar\joinrel\twoheadrightarrow}}\,
\calH \; ,  \qquad  {\check\boldalpha}_n{\big|}_{\h=1} \mapsto
\boldalpha_n \;\; .  $$
As for the Hopf structure, in  $ {\big( {\calH_\h}^{\!\prime}
\big)}^{\!\vee}{\Big|}_{\h=1} $  it is given by
  $$  \Delta \big( \check{\boldalpha}_n{\big|}_{\h=1} \big)  =
{\textstyle \sum_{k=2}^n} \, \check{\boldalpha}_k{\big|}_{\h=1}
\otimes \h^{n-k} {\Cal Q}^k_{n-k}(\boldalpha_*){\big|}_{\h=1}  +
{\textstyle \sum_{k=0}^{n-1}} \, {\check{\boldalpha}_1}^{\,k}{\big|}_{\h=1}
\otimes \h^{n-2} Z^k_{n-k}(\boldalpha_*){\big|}_{\h=1} \; .  $$
   \indent   As before, split  $ {\Cal Q}^k_{n-k}(\boldalpha_*) $  as
$ \, {\Cal Q}^k_{n-k}(\boldalpha_*) = \sum_s {\Cal T}^{s,k}_{n-k}
(\boldalpha_*) \, $,  \, and split each  $ {\Cal T}^{s,k}_{n-k}
(\widehat{\boldalpha}_*) $  into homogeneous components w.r.t.~the
total degree in the  $ \widehat{\boldalpha}_i $'s,  say  $ \,
{\Cal T}^{s,k}_{n-k}(\widehat{\boldalpha}_*) = {\textstyle \sum_r}
{\Cal Y}^{s,k}_{r,n}(\widehat{\boldalpha}_*) \, $:  \, then  $ \,
\h^{n-k-s} {\Cal T}^{s,k}_{n-k}(\widehat{\boldalpha}_*) = \h^{n-k-s}
{\textstyle \sum_r} {\Cal Y}^{s,k}_{r,n}(\widehat{\boldalpha}_*)
= {\textstyle \sum_r} \h^{n-k-s+r} {\Cal Y}^{s,k}_{r,n}
(\check{\boldalpha}_*) \, $,  \, because  $ \, \widehat{\boldalpha}_*
= \h \, \check{\boldalpha}_* \, $.  As  $ \, \h^{n-k-s+r}
{\Cal Y}^{s,k}_{r,n}(\check{\boldalpha}_*) \equiv {\Cal Y}^{s,k}_{r,n}
(\check{\boldalpha}_*) \! \mod (\h-1) \, {\big( {\calH_\h}^{\!\prime}
\big)}^{\!\vee} \, $,  \, we eventually get
  $$  \h^{n-k} {\Cal Q}^k_{n-k}(\boldalpha_*)  =  {\textstyle \sum_{s,r}}
\h^{n-k-s+r} {\Cal Y}^{s,k}_{r,n}(\check{\boldalpha}_*)  \equiv
{\textstyle \sum_{s,r}} {\Cal Y}^{s,k}_{r,n}(\check{\boldalpha}_*)
\hskip-5pt \mod \hskip-2pt (\h-1) \, {\big( {\calH_\h}^{\!\prime}
\big)}^{\!\vee}  \! =  Q^k_{n-k}(\a_*)  $$
for all  $ k $  and  $ n \, $.  Similarly  $ \, \h^{n-1} Z^k_{n-k}
(\boldalpha_*) \equiv Z^k_{n-k}(\boldalpha_*) \hskip-1pt \mod (\h-1)
\, {\big( {\calH_\h}^{\!\prime} \big)}^{\!\vee} $  ($ \forall
\, k \, $,  $ n $).  Thus
  $$  \displaylines{
   \Delta \big( \check{\boldalpha}_n{\big|}_{\h=1} \big)  =
{\textstyle \sum_{k=2}^n} \, \check{\boldalpha}_k{\big|}_{\h=1}
\otimes \h^{n-k} {\Cal Q}^k_{n-k}(\boldalpha_*){\big|}_{\h=1}  +
{\textstyle \sum_{k=0}^{n-1}} \, {\check{\boldalpha}_1}^{\,k}{\big|}_{\h=1}
\otimes \h^{n-2} Z^k_{n-k}(\boldalpha_*){\big|}_{\h=1}  =  \cr
   {} \hfill   = \;  {\textstyle \sum_{k=2}^n} \,
\check{\boldalpha}_k{\big|}_{\h=1} \otimes {\Cal Q}^k_{n-k}
(\boldalpha_*){\big|}_{\h=1}  \, + \,  {\textstyle \sum_{k=0}^{n-1}}
\, {\check{\boldalpha}_1}^{\,k}{\big|}_{\h=1} \otimes Z^k_{n-k}
(\boldalpha_*){\big|}_{\h=1}  \; .  \cr }  $$
   \indent   On the other hand, one has  $ \; \Delta(\boldalpha_n)
\, = \,  {\textstyle \sum_{k=2}^n} \, \boldalpha_k \otimes
{\Cal Q}^k_{n-k}(\boldalpha_*) \, + \, {\textstyle \sum_{k=0}^{n-1}} \,
\boldalpha_1^{\,k} \otimes Z^k_{n-k}(\boldalpha_*) \, $  in  $ \calH $,
thus the algebra isomorphism  $ \, \Omega \, \colon \, {\big(
{\calH_\h}^{\!\prime}}\big)^{\!\vee}{\Big|}_{\h=1} \, {\buildrel \cong
\over {\lhook\joinrel\relbar\joinrel\relbar\joinrel\twoheadrightarrow}}
\, \calH \, $  given by  $ \, {\widehat\boldalpha}_n{\big|}_{\h=1}
\mapsto \boldalpha_n \, $  also preserves the coproduct; similarly,
it also respects the antipode and the counit, hence it is a graded
Hopf algebra isomorphism.  In a nutshell, we have (as graded Hopf
$ \Bbbk $--algebras)  $ \; {\big({\calH_\h}^{\!\prime}\big)}^{\!\vee}
\;{\buildrel {\h \rightarrow 1} \over \llongrightarrow}\; \calH \; $.
Therefore we have filled in the bottom part of the diagram
($ \maltese $) in the Introduction, for  $ \, H = \calH \;
(:= \calH_\nu) \, $,  \, because  $ \, \calH' := \cup_{n \in \N}
D_n = \calH \, $  by \S 4.2: namely,
  $$  F \big[ {\varGamma_{\!\L_\nu}\phantom{|}}^{\hskip-8pt
\star} \big] \underset{{\calH_\h}^{\!\prime}}
\to {\overset{0 \leftarrow \h \rightarrow 1} \to
{\longleftarrow\joinrel\relbar\joinrel\relbar\joinrel\llongrightarrow}}
\, \calH \underset{{({\calH_\h}^{\!\prime})}^{\!\vee}}  \to
{\overset{1 \leftarrow \h \rightarrow 0} \to
{\longleftarrow\joinrel\relbar\joinrel\relbar\joinrel\llongrightarrow}}
U(\L_\nu)  $$
where now in right-hand side  $ \L_\nu $  is given the Lie bialgebra
structure of Theorems 4.1 and 5.1, and  $ {\varGamma_{\!\L_\nu}
\phantom{|}}^{\hskip-8pt \star} $  is the corresponding dual
Poisson group mentioned in Theorem 4.1.   

\vskip1,5truecm

\centerline {\bf \S \; 6 \ Summary and generalizations. }

\vskip10pt

  {\bf 6.1 Summary.} \, The analysis in \S\S 2--5 yields a complete
description of the  {\sl non-trivial\/}  deformations of  $ \calH $
--- namely the Rees deformations  $ {\calH_\h}^{\!\vee} $  and
$ {\calH_\h}^{\!\prime} $  and the Drinfeld's deformations
$ \big( {\calH_\h}^{\!\vee} \big)^{\!\prime} $  and  $ \big(
{\calH_\h}^{\!\prime} \big)^{\!\vee} $  ---   built out of the
{\sl trivial\/}  deformation  $ \calH_\h \, $.  In particular
  $$  \gerg_-^{\,\times} = \big(\L_\nu, \delta_\bullet\big) \;\, ,
\qquad  G_- = {G_{\!\L_\nu}\phantom{|}}^{\hskip-8pt \star} \;\, ,
\qquad \quad  G_+ = {\varGamma_{\!\L_\nu} \phantom{|}}^{\hskip-8pt
\star} \;\, ,  \qquad  \gerg_+^{\,\times} = \big(\L_\nu,\delta_*\big)
\eqno (6.1)  $$
(with notation of $ (\maltese) \, $)  where  $ \, \delta_\bullet \, $
and  $ \, \delta_* \, $  denote the Lie cobracket on  $ \L_\nu $
defined respectively in Theorem 2.1 and in Theorems 4.1 and 5.1.
Next result shows that the four objects in (6.1) are really
different, though they share some common features:

\vskip7pt

\proclaim{Theorem 6.1}
                                         \hfill\break
   \indent   (a)  $ \, \big( {\calH_\h}^{\!\vee} \big)' \cong
{\calH_\h}^{\!\prime} \, $  as Poisson  $ \, \Bbbk[\h\,] $--algebras,
but  $ \, \big( {\calH_\h}^{\!\vee} \big)' \not\cong {\calH_\h}^{\!
\prime} \, $  as Hopf  $ \, \Bbbk[\h\,] $--algebras.
                                         \hfill\break
   \indent   (b)  $ \, \big(\L_\nu, \delta_\bullet\big) \cong
\big(\L_\nu, \delta_*\big) \, $  as Lie algebras, but  $ \,
\big(\L_\nu, \delta_\bullet\big) \not\cong \big(\L_\nu,
\delta_*\big) \, $  as Lie bialgebras.
                                         \hfill\break
   \indent   (c)  $ \, {G_{\!\L_\nu}\phantom{|}}^{\hskip-8pt \star}
\cong {\varGamma_{\!\L_\nu}\phantom{|}}^{\hskip-8pt \star} \, $
as (algebraic) Poisson varieties, but  $ \, {G_{\!\L_\nu}
\phantom{|}}^{\hskip-8pt \star} \not\cong {\varGamma_{\!\L_\nu}
\phantom{|}}^{\hskip-8pt \star} \, $  as (algebraic) groups.
                                        \hfill\break
  \indent   (d) \,  The analogues of statements (a)--(c) hold
with  $ \calK $  and  $ \L_\nu^+ $  instead of  $ \, \calH $
and  $ \L_\nu \, $.   
\endproclaim  

\demo{Proof}  It follows from  Theorem 3.1{\it (a)\/}  that
$ {\big( {\calH_\h}^{\!\vee} \big)}' $  can be seen as a Poisson
Hopf algebra, with Poisson bracket given by  $ \, \{x,y\,\} := \h^{-1}
[x,y] = \h^{-1} (x \, y - y \, x) \, $  (for all  $ \, x $,  $ y
\in {\big( {\calH_\h}^{\!\vee} \big)}' \, $);  then  $ {\big(
{\calH_\h}^{\!\vee} \big)}' $  is the free Poisson algebra generated
by  $ \Big\{\, \widetilde{\b}_{x_n} \! = \widetilde{\x}_n = \a_n
\,\Big|\, n \in \N \,\Big\} \, $;  \, since  $ \, \a_n = \boldalpha_n +
(1 - \delta_{1,n}) \, {\boldalpha_1}^{\!n} \, $  and  $ \, \boldalpha_n
= \a_n - (1 - \delta_{1,n}) \, {\a_1}^{\!n} \, $  ($ \, n \in \N_+
\, $)  it is also (freely) Poisson-generated by  $ \big\{
\boldalpha_n \,\big|\, n \in \N \,\big\} $.  We also saw that
$ {\calH_\h}^{\!\prime} $  is the free Poisson algebra over  $ \,
\big\{\, \widehat{\boldalpha}_n \,\big|\, n \in \N \,\big\} \, $;
\, thus mapping  $ \, \boldalpha_n \mapsto \widehat{\boldalpha}_n
\, $  ($ \, \forall \, n \in \N \, $)  does define a unique Poisson
algebra isomorphism  $ \; \Phi \, \colon \, {\big( {\calH_\h}^{\!\vee}
\big)}' \,{\buildrel \cong \over \longrightarrow}\, {\calH_\h}^{\!
\prime} \, $,  \, given by  $ \; \widetilde{\boldalpha}_b := \h^{-d(b)}
\boldalpha_b \mapsto \widehat{\boldalpha}_b \, $,  \; for all  $ \,
b \in B_\nu \, $.  This proves the first half of  {\it (a)},  and
then also (taking semiclassical limits and spectra) of  {\it (c)\/}.
                                               \par
   The group structure of either  $ {G_{\!\L_\nu}
\phantom{|}}^{\hskip-8pt \star} $  or  $ {\varGamma_{\!\L_\nu}
\phantom{|}}^{\hskip-8pt \star} $  yields a Lie cobracket onto
the cotangent space at the unit point of the above, isomorphic
Poisson varieties: this cotangent space identifies with
$ \L_\nu $,  and the two cobrackets are given respectively
by  $ \, \delta_\bullet(x_n) = \sum_{\ell=1}^{n-1} (\ell+1)
\, x_\ell \wedge x_{n-\ell} \, $  for  $ {G_{\!\L_\nu}
\phantom{|}}^{\hskip-8pt \star} $  (by Theorem 3.1) and
by  $ \, \delta_*(x_n) = (n-2) \, x_{n-1} \wedge x_1 \, $  for
$ {\varGamma_{\!\L_\nu}\phantom{|}}^{\hskip-8pt \star} $  (by
Theorem 4.1), for all  $ \, n \in \N_\nu \, $.  It follows that
$ \, \text{\sl Ker}\,(\delta_\bullet) = \{0\} \not= \text{\sl Ker}
\,(\delta_*) \, $,  \, which implies that the two Lie coalgebra
structures on  $ \L_\nu $  are not isomorphic.  This proves  {\it
(b)},  and also means that  $ \, {G_{\!\L_\nu}\phantom{|}}^{\hskip-8pt
\star} \not\cong {\varGamma_{\!\L_\nu}\phantom{|}}^{\hskip-8pt \star}
\, $  as (algebraic) groups, hence  $ \, F \big[ {G_{\!\L_\nu}
\phantom{|}}^{\hskip-8pt \star} \big] \not\cong F \big[ {\varGamma_{\!
\L_\nu} \phantom{|}}^{\hskip-8pt \star} \big] \, $  as Hopf
$ \Bbbk $--algebras,  and so  $ \, \big( {\calH_\h}^{\!\vee}
\big)' \not\cong {\calH_\h}^{\!\prime} \, $  as Hopf
$ \Bbbk[\h\,] $--algebras,  which ends the proof of
{\it (c)\/}  and  {\it (a)\/}  too.
                                           \par
   Finally, claim  {\it (d)\/}  should be clear: one applies
the like arguments  {\sl mutatis mutandis},  and everything
follows as before.   \qed
\enddemo

\vskip7pt

   {\bf 6.2 Generalizations.} \, Plenty of features of  $ \, \calH
= \calH^\dif \, $  are shared by a whole bunch of  {\sl graded\/}
Hopf algebras, which usually arose in connection with some physical
problem or some (co)homological topic and all bear a nice combinatorial
content; essentially, most of them can be described as ``formal series''
over indexing sets   --- replacing  $ \N $  ---   of various
(combinatorial) nature: planar trees (with or without labels),
forests, graphs, Feynman diagrams, etc.  Besides the ice-breaking
examples in physics provided by Connes and Kreimer (cf.~[CK1--3]),
which are all commutative or cocommutative Hopf algebras, other
non-commutative non-cocommutative examples (like the one of
$ \calH^\dif $)  are introduced in [BF], roughly through a
``disabelianization process'' applied to the commutative Hopf algebras
of Connes and Kreimer.  A very general analysis and wealth of examples
in this context is due to Foissy (see [Fo1--3]), who also makes an
interesting study of  $ \delta_\bullet $--maps  and of the functor
$ \, H \mapsto H' \, $  ($ H $  a Hopf  $ \Bbbk $--algebra).  Other
examples, issued out of topological motivations, can be found in
the works of Loday et al.: see e.g.~[LR], and references therein.
                                             \par
   When performing the like analysis, as we did for  $ \calH $,
for a graded Hopf algebra  $ H $  of the afore mentioned type,
{\sl the arguments used for  $ \calH $  apply essentially the
same, up to minor changes, and give much the same results.}
To give an example, the Hopf algebras considered by Foissy
are non-commutative polynomial, say  $ \, H = \Bbbk \big\langle
\{x_i\}_{i \in \Cal{I}} \big\rangle \, $  for some index set
$ \Cal{I} \, $:  \, then one finds  $ \, {\calH_\h}^{\!\vee}
{\big|}_{\h=0} = U(\gerg_-) = U(\L_{\Cal{I}}) \, $  where
$ \L_{\Cal{I}} $  is the free Lie algebra over  $ \Cal{I} \, $.
                                             \par
   This opens the way to apply the methods presented in this paper
to all these graded Hopf algebras, of great interest for their
applications in mathematical physics or in topology (or whatever);
the simplest case of  $ \calH^\dif $  plays the role of a toy
model which realizes a clear and faithful pattern for many
common features of all Hopf algebras of this kind.

\vskip1,3truecm

\vskip15pt

\enddocument